\documentclass[12pt]{amsart}

\usepackage{xy, xypic,amsfonts}
\usepackage{color} 
\usepackage{graphicx}
\usepackage{psfrag}  

\usepackage{amsfonts, amssymb, amsmath, amsthm, eucal, latexsym, array,
pstricks}
\usepackage{verbatim}
\usepackage{fullpage}

\def\g{\mathfrak{g}}  
\def\q{\mathfrak{q}}  
\def\l{\mathfrak{l}}  
\def\z{{\mathfrak{z}}}  
\def\a{{\mathfrak{a}}}  %
\def\b{\mathfrak{b}}  
\def\t{{\mathfrak{t}}}  
\def\r{{\mathfrak{r}}}  
\def\h{\mathfrak{h}}  
\def\p{\mathfrak{p}}  
\def\m{\mathfrak{m}}  
\def\n{\mathfrak{n}}  
\def\s{\mathfrak{s}}  

\def\C{\mathbb{C}}      

\def\reg{{\rm reg}}            

\def\red{\mathrm{red}}     
\def\kg#1{\mathrm{k}_{#1}} 

\def\varp#1{\varphi_{#1}} 
\def\var#1{\varphi_{#1}}            
\def\hr#1{\theta_{#1}}        
\def\cal#1{\mathcal{#1}}

\def\ep{\varepsilon}

\newtheorem{thm}{Theorem}[section]
\newtheorem{lm}[thm]{Lemma}
\newtheorem{cl}[thm]{Corollary}
\newtheorem{prop}[thm]{Proposition}

\theoremstyle{remark}
\newtheorem{ex}[thm]{Example}
\newtheorem{rmk}[thm]{Remark}

\theoremstyle{definition}
\newtheorem{df}[thm]{Definition}

\newcommand{\gt}{\mathfrak}
\newcommand{\cp}{\mathbb C}

\newcommand{\id}{{\rm id}}
\newcommand{\ind}{{\rm ind}\,}

\newcommand{\rk}{\mathrm{rk\,}}

\newcommand{\Lie}{\mathrm{Lie\,}}

\newcommand{\Ann}{\mathrm{Ann}}
\newcommand{\ad}{\mathrm{ad}}

\renewcommand{\le}{\leqslant}
\renewcommand{\ge}{\geqslant}

\newfam\eusfam%
\font\euszw=eusm10 scaled 1200%
\font\eusac=eusm7 scaled 1200%
\font\eusacc=eusm7 scaled 1000%
\textfont\eusfam=\euszw\scriptfont\eusfam=\eusac%
\scriptscriptfont\eusfam=\eusacc%
\title[Coadjoint orbits of reductive type]
{Coadjoint orbits of reductive type of seaweed Lie algebras}

\author[A.~Moreau]{Anne Moreau}
\address{Anne Moreau, Laboratoire de Math\'ematiques et Applications, Universit\'e de Poitiers, France}
\email{anne.moreau@math.univ-poitiers.fr}

\author[O.~Yakimova]{Oksana Yakimova}
\address{Oksana Yakimova, Emmy-Noether-Zentrum, Department Mathematik, Universit\"at Erlangen-N\"urnberg, Germany}
\email{yakimova@mccme.ru}

\date{August 21, 2011}  


\newfam\eusfam%
\font\euszw=eusm10 scaled 1200%
\font\eusac=eusm7 scaled 1200%
\font\eusacc=eusm7 scaled 1000%
\textfont\eusfam=\euszw\scriptfont\eusfam=\eusac%
\scriptscriptfont\eusfam=\eusacc%
\begin{document}

\begin{abstract}

A connected algebraic group $Q$ defined over a field of characteristic zero is {\it quasi-reductive}
if there is an element of ${\gt q}^*$ of reductive type, that is such that the quotient
of its stabiliser by the centre of $Q$ is a reductive subgroup of $GL(\gt q)$.
Such groups appear in harmonic analysis when unitary representations are studied.
In particular, over the field of real numbers they turn out to be the groups with discrete series
and  their irreducible unitary square integrable representations are
parameterised by coadjoint
orbits of reductive type.
Due to results of M.~Duflo, coadjoint representation of a quasi-reductive $Q$ possesses a
so called {\it maximal reductive stabiliser} and knowing this subgroup,
defined up to a conjugation in $Q$, one can describe all
coadjoint  orbits of reductive type.

In this paper, we consider quasi-reductive parabolic subalgebras of simple complex
Lie algebras as well as seaweed subalgebras  of $\gt{gl}_n(\cp)$
and describe the classes of their maximal reductive stabilisers.
\end{abstract}
\maketitle


%
\section{Introduction}
\label{sec:Intro}

\subsection{}
Suppose that $\Bbbk$ is a field of characteristic zero and
$Q$ a connected algebraic (or Lie) 
group defined over $\Bbbk$.
Let  $\q=\Lie Q$ be the Lie algebra of $Q$. Let $Z$ denote the centre of $Q$.
A linear function $\gamma\in\q^*$ is said to be of
{\it reductive type}  if the quotient
$Q_\gamma/Z$ of the stabiliser $Q_\gamma\subset Q$ for the coadjoint action is
a reductive subgroup of $GL(\gt q^*)$ (or, what is equivalent, of
$GL(\q)$). Whenever it makes sense, we will also say that $\gamma$ is
of {\it compact type} if $Q_\gamma/Z$ is compact.
Further, $Q$ and $\q$ are said to be \emph{quasi-reductive}
if there is $\gamma\in\q^*$ of reductive type.

The notions go back to M.~Duflo, who
initiated the study of such Lie algebras because of applications
in harmonic analysis, see~\cite{Du}.
In order to explain his (and our) motivation let us assume for a while that
$\Bbbk=\mathbb R$.

A classical problem is to describe the unitary dual $\widehat{Q}$ of $Q$,
i.e.,~the equivalence classes of unitary irreducible representations of $Q$.
In this context, the coadjoint orbits play a fundamental r{\^o}le.
When $Q$ is a connected, simply connected, nilpotent Lie group,
Kirillov's ``orbit method''~\cite{Kir} provides a bijection between
$\widehat{Q}$ and $\q^*/Q$.
Great efforts have been made to extend the orbit method to arbitrary groups by Kostant, Duflo, Vogan and many others.
In case of square integrable representations the extension is particularly successful.

Let $(\pi, \mathcal{H})$ be an irreducible unitary representation of $Q$.
Then central elements of $Q$ act on $\mathcal{H}$ as scalar multiplications.
Hence for $v,w\in\mathcal{H}$ and $z\in Z$, $q\in Q$ the norm $|\left<v,\pi(qz)w\right>|^2$ does not depend on $z$.
Therefore $|\left<v,\pi(qz)w\right>|^2$ is a function on $Q/Z$ and the following definition makes sense:

\begin{df}\label{sq-int-mod-Z} An irreducible unitary representation
$(\pi, \mathcal{H})$ is said to be {\it square integrable modulo $Z$} if
there exists a non-zero matrix coefficient such that the integral of
the function $q\mapsto |\left<v,\pi(qz)w\right>|^2$ 
over a left invariant Haar measure on $Q/Z$ is finite.
\end{df}

Let $\widehat{Q}_2$ denote the set of the equivalence classes of irreducible square integrable modulo $Z$ representations.
In case of a nilpotent Lie group $N$ we have a result of C.~Moore and J.~Wolf~\cite{Joe1} stating that $\widehat{N}_2\ne\varnothing$ if and only if
there is $\gamma\in\gt n^*$ of compact type, where $\n=\Lie N$. This also can be formulated as
$\gamma$ is of reductive type or $N_\gamma=Z$.
In case of a reductive group $G$,  Harish-Chandra's theorems \cite{HC1},~\cite{HC2} give a bijection between $\widehat{G}_2$ and certain forms of compact type.
Square integrable representations of unimodular groups were described by
Ahn~\cite{Ahn1}, \cite{Ahn2}, and the existence of such representations leads to a rather restrictive conditions on 
the group.
For an arbitrary Lie group $Q$, description of $\widehat{Q}_2$
was obtained by Duflo~\cite[Chap.~III, Th\'eor\`eme 14]{Du}.

If $\widehat{Q}_2\ne\varnothing$, then $Q$ is quasi-reductive.
In the other direction, having $\gamma\in\q^*$ of compact type, there is
a way to  construct a square integrable unitary representation \cite[Chap.~III,~\S\S 14 and 16]{Du}.

In this paper we consider complex Lie groups and
linear forms $\gamma$ of reductive type.
If $\gamma\in\gt q^*$ is of compact type, then $\gamma$ extends to a function on
$\q{\otimes}\cp$ and the extension is of reductive type.
Going in the other direction is more difficult. One has to
find a suitable real form of $\gt q{\otimes}\cp$ and prove that
the stabiliser in it is compact.
We do not address these problems here.

\subsection{}
The coadjoint 
action of a quasi-reductive linear Lie group 
has quite remarkable properties.
In Section~\ref{S:gen}, we recall relevant results
of \cite[Section 3]{DKT} and 
give independent
proofs of them under the assumption that the centre $Z$ of $Q$
consists of semisimple elements of $Q$.
These results can be also traced back to~\cite{Du}.
One if them is that there is a unique coadjoint orbit  $Q\gamma$ of reductive type
such that the restriction of $\gamma$ to its stabiliser $\q_{\gamma}$ is zero.
Thus, the stabilisers of such $\gamma$ are all conjugate by the group $Q$.
Moreover, if $\beta\in\gt q^*$ is of reductive type, then
$Q_\beta$ is contained in $Q_\gamma$ up to conjugation.
%
Therefore the following  definition  is justified.

\begin{df}
Assume that $Z$ consists of semisimple elements. 
If $\gamma$ is a linear form of reductive type such that the restriction of $\gamma$ to $\q_{\gamma}$ is zero,
then its stabiliser $Q_\gamma$ is called a \emph{maximal reductive stabiliser},
MRS for short.
Let  $M_*(\q)$ denote the Lie algebra of a MRS.
The reader shall keep in mind that both MRS and $M_*(\q)$ are defined up to conjugation.
\end{df}

\noindent
We call a quasi-reductive Lie algebra $\gt q$ {\it strongly} quasi-reductive if
$Z$ consists of semisimple elements.

It also follows from Duflo's results~\cite{Du} (\cite[Th\'eor\`eme 3.6.2]{DKT} or see Corollary~\ref{red-orbits} here),
that the coadjoint orbits of $Q$ of reductive type are parametrised  by the
coadjoint orbits of reductive type of MRS.
The latter are in one-to-one correspondence with the semisimple
adjoint orbits of MRS, that is the closed adjoint orbits of MRS.
Thus, the set $\gt q^*_{\rm red}$ of linear forms of reductive type
has an affine geometric quotient with respect to the action of
$Q$ and $\gt q^*_{\rm red}/Q\cong M_*(\gt q)/\!\!/ {\rm MRS}$.

\subsection{Digression to positive characteristic}
It is worth mentioning that finite-characteristic analogues of quasi-reductive
Lie algebras have some nice properties.
Let $\mathbb F$ be an algebraically closed field of characteristic $p>0$
and $Q$ an algebraic group defined over $\mathbb Z$ such that
$Q(\cp)$ is strongly quasi-reductive.
Reducing scalars modulo $p$ one may view $\gt q$ as a Lie algebra
over $\mathbb F$. Identifying $\gt q$ with the $Q$-invariant
derivations of $\mathbb F[Q]$ and taking the $p$'th power of a derivation,
we get a $p$-operation: $x\to x^{[p]}$ on $\gt q$ turning it into
a restricted Lie algebra $\widetilde{\gt q}$.
Since $Q(\cp)$ is strongly quasi-reductive, for almost all $p$
there is $\gamma\in\widetilde{\gt q}^*$ such that
$\widetilde{\gt q}_\gamma$ is a toral subalgebra of $\widetilde{\gt q}$.
This implies that $\widetilde{\gt q}$ satisfies a
Kac-Weisfeiler conjecture (the so called KW1 conjecture)
on
the maximal dimension of irreducible
$\widetilde{\gt q}$-modules, \cite[Section~4]{PS}.

\subsection{Parabolic subalgebras}
Classification of quasi-reductive Lie algebras seems to be a wild problem.
(Easy examples of them are reductive algebras and Abelian algebras.)
However, it is reasonable to look on specific subalgebras of semisimple
Lie algebras, in particular on parabolic subalgebras.

By a classical result of Dynkin, a maximal subgroup of a semisimple
group is either reductive or a parabolic subgroup.
When studying branching rules one naturally is tempted to restrict first
to a maximal subgroup. This is one of the instances where parabolic subgroups come into play in
harmonic analysis. Of particular interest are branchings of
square integrable representations with finite multiplicities.
In that case a subgroup 
must have an irreducible square integrable representation,
in other words, it must be quasi-reductive.

More generally, the algebras of seaweed type,
first introduced by V.~Dergachev and A.~Kirillov \cite{DK}
in the $\mathfrak{sl}_n$ case, 
form a very interesting class of non-reductive subalgebras
of semisimple Lie algebras.
They naturally extend both the classes of parabolic subalgebras and of Levi subalgebras.
A {\it seaweed subalgebra}, or {\it biparabolic subalgebra},~\cite{Jo1}, of a semisimple Lie algebra is the intersection of
two parabolic subalgebras whose sum is the total Lie algebra.
They have been intensively studied these last years,
see~e.g.~\cite{P1},~\cite{Dv},~\cite{TY2},~\cite{P3}, \cite{Jo1},~\cite{Jo2},~\cite{BM}.

Observe that the centre of a biparabolic subalgebra always consists of semisimple elements of the total Lie algebra.
Therefore if a biparabolic subalgebra is quasi-reductive, it is
strongly quasi-reductive (Definition~\ref{d:qua}).
The biparabolic subalgebras in $\mathfrak{sl}_n$ and $\mathfrak{sp}_{2n}$ are quasi-reductive
by a result of D.~Panyushev~\cite{P3}.
The classification of quasi-reductive parabolic subalgebras of reductive Lie algebras has been
recently completed in~\cite[Section 5] {DKT} and~\cite{BM}.

\subsection{Description of the paper}
In this paper, we focus on the conjugation class of MRS in quasi-reductive parabolic subalgebras of a reductive Lie algebra.
This subject of study was suggested to the first author 
by M.~Duflo and was the original motivation of~\cite{BM}.

The problem clearly reduces to the case of quasi-reductive parabolic subalgebras of simple Lie algebras.
Let $\gt q\subset\gt g$ be a parabolic subalgebra of a simple Lie algebra $\gt g$
or a seaweed subalgebra in $\gt g=\gt{sl}_n$ (or $\gt{gl}_n$).
In this note, we describe $M_*(\gt q)$ and also
specify an embedding $M_*(\gt q)\subset\gt q$.
This allows us to get back MRS and $\gt q^*_{\rm red}$ in the following way.
Set
$$
\Upsilon:=(\gt q^*)^{M_*(\gt q)}\cap{\rm Ann}(M_*(\gt q))
 =\{\xi\in\gt q^*\mid \xi([\gt q,M_*(\gt q)]=0, \xi(M_*(\gt q))=0\}.
$$
Clearly, if $Q_\gamma$ is a maximal reductive stabiliser, then
$\gamma\in Q\Upsilon$.
Since $Q$-orbits of the maximal dimension form an open subset of
$Q\Upsilon$, for generic $\xi\in\Upsilon$ we have
$\gt q_\xi=M_*(\gt q)$. Thus generic $\xi\in \Upsilon$
is of nilpotent and reductive type at the same time (see Section~\ref{S:gen} for the definitions).
By Proposition~\ref{MRS-first}({\sf i}), $Q_\xi$ is also a
maximal reductive stabiliser.
Identifying a
$Q_\xi$-invariant complement of ${\rm Ann}(\gt q_\xi)\subset\gt q^*$
in $\gt q^*$  with $\gt q_\xi^*$ we get
$(\gt q^*)_{\rm red}=Q(\xi+(\gt q_\xi^*)_{\rm red})$
(see proof of Pro\-po\-si\-tion~\ref{MRS-first}({\sf ii})).
We do not describe explicitly the component group of MRS.
However, note that for seaweed subalgebras in $\gt{gl}_n$,
the maximal reductive stabilisers are connected.

Roughly speaking, we have two methods for calculating  $M_*(\gt q)$.
The first uses root system of $\gt g$ and related objects, e.g. Kostant's cascade
(see Section~\ref{S:sea} for more details). It allows us
to prove the crucial ``additivity'' property (Theorem~\ref{t:add}).
This property assures that it suffices to consider parabolic subalgebras with simple semisimple
parts and is particularly useful in the exceptional case.

The second method uses algebraic Levi decomposition
$\gt q=\gt l\ltimes\gt n$ and allows us to ``cut'' an Abelian ideal in
$\gt n$ (Lemma~\ref{ideal}).  In principle,  this reduction can be applied
to any algebraic linear Lie algebra such that its centre consists of semisimple
elements and at the end establish whether the algebra is quasi-reductive or not.
In particular the method nicely works for
$\mathbb Z$-graded Lie algebras with $2$ or $3$ graded components
(Lemmas~\ref{2-grad}, \ref{3-grad}).
Similar gradings were used before by Panyushev~\cite{P1}, \cite{P3}
to calculate the index of a biparabolic subalgebra and its generic stabiliser.
This approach is used in Section~\ref{S:cla} to deal with the classical Lie algebras.
%

While Kostant's cascade is strongly used in~\cite{BM} to get generic reductive stabilisers,
here we do not use it directly to construct MRS.
But one can observe relations between the both approaches.
In fact, the description using roots system is handier to write down an algorithm, for instance if one wants to
use the computer program~\texttt{GAP} to get MRS (cf. Section~\ref{ss:cas}).

For seaweed subalgebras in $\gt{gl}_n$ our answer is given in terms of the meander graph,
used in \cite{DK} to express the index of $\gt q$.
Recall that a  seaweed  (or biparabolic) subalgebra $\gt q$ is an intersection of two complementary parabolic subalgebras. 
Hence, up to conjugation, it is defined by two
compositions $\bar a,\bar b$ of $n$. To these objects one can attach
a certain graph $\Gamma=\Gamma(\bar a|\bar b)$ with $n$ vertices and
at most $n$ edges.
For example $\Gamma(5,2,2|2,4,3)$ has $\bar a$-edges
$(1,5)$, $(2,3)$, $(6,7)$, $(8,9)$ and $\bar b$-edges
$(1,2)$, $(3,6)$, $(4,5)$, $(7,9)$  (for a picture see Section~\ref{A}).
A cycle of $\Gamma$ is said to be maximal if it does not lie inside any
other cycles in the planar embedding of $\Gamma$.
To a maximal cycle one can adjust a number, $r$, its dimension,
and a subgroup $GL_{r}\subset GL_n$. Our result states
that MRS of $\gt q$ is equal to the product of $GL_r$ over all maximal cycles in $\Gamma$.
This 
confirms a prediction of Duflo.
We furthermore describe an explicit embedding of MRS.

Finally, the exceptional Lie algebras  are dealt with in
Section~\ref{S:exc} and the results are stated in
Tables~\ref{tE6},~\ref{tE7},~\ref{tE8},~\ref{tF4}.

\smallskip

It is worth noticing that our paper yields alternative proofs of some results of~\cite{BM}.
Namely, each time we explicitly describe a (maximal) reductive stabiliser, this shows at the same time that the given Lie algebra is quasi-reductive.
This is especially interesting in some exceptional cases, where \cite{BM} merely proves the quasi-reductivity by using \texttt{GAP} (cf.~Remark~\ref{r:gap} for more details).

\bigskip

\subsection*{Acknowledgments}
We 
would like to
thank Michel Duflo for 
bringing 
quasi-reductive Lie algebras
to our attention, explanations concerning their coadjoint representations,
and his interest in this work.
It is also a pleasure 
to thank Alexander Premet for useful 
discussions on the subject of
restricted Lie algebras.

\setcounter{tocdepth}{1}
\tableofcontents

%
\section{Strongly quasi-reductive Lie algebras} \label{S:gen}                    
%


From now on, $\Bbbk=\C$.
Let $Q$ be a linear algebraic group and $\gt q=\Lie Q$ its Lie algebra.
Keep the notation of the introduction.
Set $\z:=\Lie Z$. Recall that a linear form $\gamma\in \q^*$
is of  \emph{reductive type} if
$Q_\gamma/Z$ is a reductive Lie subgroup of $GL(\q)$.
Let $\q^{*}_{\red}\subset\gt q^*$
denote the set of linear forms of reductive type.
The group $Q$ acts on $\q^{*}_{\red}$ and we denote by $\q^{*}_{\red}/Q$ the
set of coadjoint orbits of reductive type.
Recall also that $\q$ is called {\it quasi-reductive} if it has linear forms
of reductive type.
Most results of this section are due to Duflo et al.\ and are contained in~\cite[Section 3]{DKT}.
We give here independent proofs for the convenience of the reader.
For our purpose, the following definition will be very useful as well.

\begin{df} \label{d:qua}
If $\q$ is quasi-reductive and $\z$ consists of semisimple elements of $\q$,
then $\q$ is said to be  {\it strongly quasi-reductive}.
\end{df}

In this section we concentrate on strongly quasi-reductive Lie algebras.
These algebras can be characterised by the property that there exists a linear form
$\gamma\in\gt q^*$ such that $\gt q_\gamma$ is a reductive subalgebra of $\gt q$.
However, the reader can keep in mind that most of results
stated in this section are still true for an arbitrary $\q$ (\cite[Section 3]{DKT}).
The statements have then to be slightly modified accordingly.

\begin{df}\label{gen-stab}
Suppose that $Q$ acts on an irreducible affine variety $Y$.
Then a subgroup $Q_y$ (with $y\in Y$) is called a
{\it generic stabiliser} of this action if there is an open subset
$U\subset Y$ such that $Q_y$ and $Q_w$ are conjugate in $Q$ for
all $w\in U$.
\end{df}

By a deep result of Richardson \cite{R},
a generic stabiliser exists for any action of a reductive  algebraic group
on a smooth affine variety.
We will say that a Lie algebra of a generic stabiliser is a generic stabiliser as well.
By \cite[Prop.~4.1]{R}, if a generic stabiliser exists on the Lie algebra level,
it also exists on the group level.

In case of a coadjoint representation, the following lemma
can be deduced from
\cite[III]{Du2} and \cite{Du0}. We give here a different 
proof applicable in a more general setting.

\begin{lm}\label{exist-gen-stab} Suppose that $Q$ acts on a linear space
$V$ and there is $v\in V$ such that $Q_v$ is reductive.
Then the action of $Q$ on $V$ has a generic stabiliser.
In particular, if $\gt q$ is strongly quasi-reductive,
then the coadjoint action of $Q$ possesses  a generic stabiliser $T$
such that $T$ is reductive and $\gt t=\Lie T$ is Abelian.
\end{lm}
 \begin{proof}
 This can be considered as a version of the Luna slice theorem
 \cite{Luna}, see also \cite[proof of Prop.~1.1]{gib}.
Since $Q_v$ is reductive,  there
is a $Q_v$-stable complement of T$_v(Qv)$ in $V$, say $N_v$.
Let us consider the associated fibre bundle
$X_v:=Q\ast_{Q_v}N_v$ and  recall that it is the (geometric) quotient of
$Q\times N_v$ by the $Q_v$-action defined by
$Q_v\times Q\times N_v\to Q\times N_v$, $(s,q,n)\mapsto (qs^{-1},sn)$.
The image of $(q,n)\in Q\times N_v$ in $X_v$ is denoted by $q\ast n$.
The natural $Q$-equivariant
morphism $\psi: X_v\to V$, $\psi(q\ast n)=q(v+n)$  is \'etale in
$e\ast 0\in X_v$ by construction. It follows that
there is an open $Q$-stable neighbourhood
$U$ of $Qv$ such that for all $y\in U$ the identity component
$Q_y^{\circ}$ is conjugate to $(Q_v)_x^{\circ}$ with $x\in N_v$.
Therefore on the Lie algebra level a generic stabiliser of the
action $Q_v\times N_v\to N_v$ is also a generic stabiliser
of $Q\times V\to V$. (The statement also follows from
\cite[Prop.~3.3]{R}.)
By  \cite[Prop.~4.1]{R}, a
generic stabiliser exists also on the group level.

In case $V=\gt q^*$, we have
T$_v(Qv)\cong {\rm Ann}(\gt q_v)$ and
$N_v\cong \gt q_v^*$.
Therefore a maximal torus $\gt t$ in $\gt q_v$ is a generic stabiliser
for the coadjoint action of  $\gt q$.
Since
all stabilisers $Q_y$ are algebraic groups,
generic stabilisers in $Q$ are reductive as well.
\end{proof}

Recall that the \emph{index} of a Lie algebra
is the minimal dimension of a stabiliser in the coadjoint representation,
$$\ind\q= \min \{\dim \q_{\gamma} \mid \gamma \in \q^*\}.$$
A linear form $\gamma \in \g^*$ is \emph{regular} if
$\dim\q_{\gamma}=\ind\q$.
Let $\q^{*}_{\reg}$ denote the set of regular linear forms;
it is an open dense subset of $\q^*$.
For reductive Lie algebras, the index is equal to the rank.
For Abelian Lie algebras,
the index is equal to the dimension.
Both instances are easy examples of quasi-reductive Lie algebras.

\begin{cl} \label{tor}
Suppose that  $\gt t=\gt q_\xi$ is a reductive torus for $\xi\in\q^*$.
Then $\gt t$ is a generic stabiliser for the coadjoint action.
\end{cl}
\begin{proof}
The statement immediately
follows from the description of a slice associated with
$Q\xi$.
\end{proof}

By a classical result of Duflo and Vergne, \cite{DV},
for any Lie algebra $\gt q$ and any $\gamma\in\gt q^*_{\reg}$,
the stabiliser $\gt q_\gamma$ is Abelian, see e.g. \cite[Corollary~1.8]{P2}.

Suppose that $\gt t=\gt q_\alpha$ is reductive and
$\alpha\in\gt q^*_{\reg}$. In this case $\gt t$ is a torus.
Let $\gt h:=\gt z_{\gt q}(\gt t)$ be the centraliser of $\gt t$
in $\gt q$, i.e., the set of elements commuting with $\gt t$.
 Since $\gt q$ is algebraic, it possesses
a Levi decomposition $\gt q=\gt l\ltimes\gt n$, where $\gt n$ is the nilpotent
radical and $\gt l$ is a maximal reductive subalgebra. We may assume that
$\gt t\subset\gt l$. 
Since $Q_\alpha^{\circ}\subset Q$ is an algebraic torus, there exists 
a $\gt q$-invariant bilinear form, say $B$, on $\gt q$ such that it is non-degenerate on $\gt t$.
Note that $\gt n$ lies in the kernel of
$B$.  The set
$\gt m_0:=\gt t^{\perp}\cap\gt h$ is a subalgebra of $\gt q$
and $\gt h=\gt t{\oplus}\gt m_0$.
We have also
\begin{equation}\label{decomp}
\gt q=\gt t\oplus\gt m_0\oplus [\gt t,\gt q], \text{  where  }\
\gt m_0\oplus[\gt t,\gt q]=\gt t^{\perp}\  \text{ and }\
 \gt n\subset \gt m_0\oplus [\gt t,\gt q].
\end{equation}
This also leads to the dual decomposition of
$\gt q^*$. 
Since $\ad^*(\gt t)\alpha=0$, we see that $\alpha([\gt t,\gt q])=0$.
One can write $\alpha$ as $\alpha=s+\gamma+0$,
where $\gamma$ is zero on $\gt t$ and $[\gt t,\gt q]$
(and $s$ is zero on ${\gt m}_0$ and $[\gt t,\gt q]$).
Clearly $\gt t\subset\gt q_\gamma$.
One the most important properties of $\gamma$ is that it is of reductive type.
This is a result of \cite[Th\'eor\`eme 3.4.2]{DKT} and we also give a proof below.

\begin{ex}\label{gamma-red} Suppose that $\gt q$ is reductive.
Then the above construction gives $\gamma=0$ and $\gt q_\gamma=\gt q$.
\end{ex}

\begin{ex}\label{gamma-heis} Suppose that $\gt q=\gt l\ltimes\gt n$,
where $\gt n$ is a Heisenberg Lie algebra and $\gt l$ is reductive.
Let $\alpha\in\gt q^*_{\reg}$ with $\gt q_\alpha=\gt t$ be of reductive type and
$\gt u$ the centre of $\gt n$.
Then $\alpha|_{\gt u}\ne 0$.
Consider the function $\tilde\alpha=\alpha|_{\gt u}$
on $\gt u$. Let $\gt v\subset\gt n$ be the kernel of
$\alpha|_{\gt n}$. Note that $[\gt t,\gt v]\subset\gt v$
and $\gt n=\gt v{\oplus}\gt u$. Moreover $[\gt n,\gt v]\subset\gt u$.
Since $[\gt q,\gt v]\subset \gt v{\oplus}\gt u$,
the normaliser $\hat{\gt l}$ of $\gt v$ in $\gt q$ is a subalgebra of
dimension at least $\dim\gt l+1$ and its intersection with $\gt n$
is equal to $\gt u$. Therefore $\hat{\gt l}/\gt u\cong\gt l$.
Replacing
$\gt l$ by a conjugate subalgebra of $\gt q$
(a reductive part of  $\hat{\gt l}$) we may assume that
$\gt v$ is also $\gt l$-stable.
Since  $\dim\gt u=1$, the stabiliser $\gt l_{\tilde\alpha}$
is reductive and there is an orthogonal decomposition
$\gt l=\gt l_{\tilde\alpha}\oplus\cp \eta$ with $\eta$ being
a central (in $\gt l$) element. Note that
$\gt t$ is a maximal torus in $\gt l_{\tilde\alpha}$.
The above construction
produces a function $\gamma$ such that $\gamma|_{\gt u}=\alpha|_{\gt u}$ and
$\gamma(\gt l_{\tilde\alpha})=0$, $\gamma(\gt v)=0$. Thereby  $\gt q_\gamma=\gt l_{\tilde\alpha}$
is reductive.
\end{ex}

The general case is proved by induction on $\dim\gt q$.

\begin{lm}[\cite{DKT}]\label{unipotent-exists}
Suppose that $\gt q$ is strongly quasi-reductive,
$\gt t=\gt q_\alpha$ is a generic (reductive) stabiliser, and
$\gamma\in\gt q^*$ is obtained from $\alpha$ by means of
the decomposition~(\ref{decomp}) as above.
Then $\gt q_\gamma$ is reductive. Moreover, $\gamma(\gt q_\gamma)=0$.
\end{lm}
\begin{proof}
Note first that 
$[\gt m_0,\gt q]\subset\gt t^{\perp}$ and hence
$\gamma$ is equal to $\alpha$ on $[\gt m_0,\gt q]$. 
Thereby $\gt q_\gamma\cap\gt m_0=0$ and since $\gt t\subset\gt q_\gamma$,
we get
$\gt q_\gamma\subset\gt t\oplus[\gt t,\gt q]$. In particular,
$\gamma(\gt q_\gamma)=0$.
Recall also that $\gt n\subset \gt t^{\perp}$ and hence
$\alpha|_{\gt n}=\gamma|_{\gt n}$.

We prove that $\gt q_\gamma$ is reductive by induction on
$\dim\gt q$. Let $A$ be a bilinear form on
$\gt n{\times}\gt q$ defined by
$A(\eta,\xi):=\alpha([\eta,\xi])=\gamma([\eta,\xi])$.
Set $\ker A:=\{\xi\in\gt q\mid A(\gt n,\xi)=0\}$ and
$\gt n(A):=\gt n\cap \ker A$. Let $\gt w\subset \gt n$ be
a complement of $\gt n(A)$ in $\gt n$. Then $A$ is non-degenerate on
$\gt w{\times}\gt w$.
Since $\gt q_\alpha\cap \gt n=0$, there is a subspace
$\tilde{\gt v}\subset\gt q$ of dimension $\dim\gt n(A)$ such that the pairing
$A(\tilde{\gt v},\gt n(A))$ is non-degenerate.
Set $\gt v:=\{v\in\tilde{\gt v}{\oplus}\gt w\mid A(v,\gt w)=0\}$.
Then $\dim\gt v=\dim\tilde{\gt v}$ and $\gt v\cap\gt w=0$,
because $A$ is non-degenerate on
$\gt w{\times}\gt w$. Since also $A(\gt n(A),\gt w)=0$,
we get that $A$ is non-degenerate on  $\gt n(A){\times}\gt v$.

Finally in $\ker A$ we fix a decomposition
$\ker A=\gt n(A)\oplus\gt s$, where
$\gt s=\gt s(\alpha)$ or $\gt s=\gt s(\gamma)$ has the property that
$\alpha([\gt s,\gt v])=0$ or $\gamma([\gt s,\gt v])=0$, respectively.
Let us choose a basis
of $\gt q$ according to the inclusions
$$
\gt w\subset\gt w\oplus\gt n(A)\subset \gt n \oplus \gt s \subset
\gt n + \ker A\subset
 \gt w\oplus\ker A\oplus\gt v=\gt n \oplus \gt s \oplus \gt v.
$$
Then the matrices of the forms
$\hat\alpha(\xi_1,\xi_2)=\alpha([\xi_1,\xi_2])$ and
$\hat\gamma$ look as shown in Picture~\ref{matrices}.

\begin{figure}[htb]
{\setlength{\unitlength}{0.1in}
\begin{center}
\begin{picture}(14,14)(0,0)

\put(0,0){\line(1,0){14}}\put(0,0){\line(0,1){14}}
\put(14,0){\line(0,1){14}}\put(0,14){\line(1,0){14}}

\put(0,10){\line(1,0){4}}\put(4,10){\line(0,1){4}}
\qbezier[24](0,8)(3,8)(6,8)
\qbezier[24](6,8)(6,11)(6,14)
\put(0,2){\line(1,0){14}}\put(12,0){\line(0,1){14}}
\put(6,0){\line(0,1){8}}\put(6,8){\line(1,0){8}}
\put(4,0){\line(0,1){2}}\put(12,10){\line(1,0){2}}

\put(1,6.5){{\Huge$0$}}\put(3.5,3.5){{\Huge$0$}}
\put(5,11){{\Huge$0$}}\put(9,9){{\Huge$0$}}
\put(1.5,0.4){{\large$0$}}\put(12.5,11.5){{\large$0$}}

\put(0.5,11){{\Large$A_1$}}
\put(4,0.5){$A_2$}\put(12.05,8.5){$A_3$}\put(8,4){{\Huge$C$}}
\put(9,0.5){$0$}
\put(12.7,0.6){$\ast$}
\put(12.7,4.5){$0$}

\put(-1.5,10.65){$\left\{ {\parbox{1pt}{\vspace{6\unitlength}}}  \right.$}
\put(-1.5,4.5){$\left\{ {\parbox{1pt}{\vspace{5.7\unitlength}}}  \right.$}
\put(-1,0.4){$\left\{ {\parbox{1pt}{\vspace{1.8\unitlength}}}  \right.$}
\put(-2.5,10.5){{$\gt n$}}
\put(-2.2,4.5){{$\gt s$}}
\put(-2,0.5){{$\gt v$}}
\put(14,11.6){$\left. {\parbox{1pt}{\vspace{4\unitlength}}}  \right\}$}
\put(14,8.5){$\left. {\parbox{1pt}{\vspace{1.9\unitlength}}}  \right\}$}
\put(15.2,11.5){{$\gt w$}}
\put(15.1,8.5){{$\gt n(A)$}}
\put(13.6,5.5){$\left. {\parbox{1pt}{\vspace{7.9\unitlength}}}  \right\}$}
\put(15.3,5.4){{$\ker A$}}
\end{picture}

\text{Here matrices $A_i$ are non-degenerate.}
\end{center} }
\caption{}\label{matrices}
\end{figure}

\noindent
Since $\gt q_\gamma=\ker\hat\gamma$, Picture~\ref{matrices}
tells us that $\gt q_\gamma\subset\ker A$ and
$\gt q_\gamma\cap\gt n(A)=0$. Moreover,
$\gt q_\gamma$ coincides with the kernel of
a skew-symmetric form defined by the matrix $C=C(\gamma)$
in $\gt s=\gt s(\gamma)$. The same holds for $\gt q_\alpha$ if we replace
$C(\gamma)$ with $C(\alpha)$.
In order to formalise this we set
$$
\gt q(A):=\{\xi\in\gt q \mid A(\,.\,,\xi)\in\cp\alpha|_{\gt n}\}/\ker\alpha\cap\gt n(A).
$$
Note that $\gt q(A)$ is a Lie algebra,
$\gt t\subset\gt q(A)$, and $B$ canonically induces an
invariant bilinear form on  $\gt q(A)$
preserving the condition $\gt t^{\perp}\cap\gt t=0$.
Both $\alpha$ and $\gamma$ can be restricted to $\gt q(A)$,
we keep the same letters for all the restrictions.

Equality $\dim\gt q(A)=\dim\gt q$ is possible only in two cases:
$\gt q$ is reductive, this is treated in   Examples~\ref{gamma-red},
or $\dim\gt n=1$, which is a particular case of Example~\ref{gamma-heis}.
In the following we assume that $\dim\gt q(A)<\dim\gt q$.

The centre of $\gt n$ has zero intersection with $\gt q_\alpha$ and
is obviously contained in $\gt n(A)$.
Hence there exists $\eta\in\gt n(A)$ such that
$\alpha(\eta)\ne 0$.
Let $\bar\eta=\eta+\ker\alpha\cap\gt n(A)$
be an element of $\gt q(A)$. Next
$\dim\{A(\,.\,,\xi)\mid \xi\in\gt q\}=\dim\gt n$ and
therefore this set is isomorphic to $\gt n^*$.
Hence there are elements $l\in \gt q$  multiplying
$\alpha|_{\gt n}$ by a non-zero constant.
We also consider $l$ as an element of $\gt q(A)$.
Then
$$
\gamma([l,\bar\eta])=\gamma([l,\eta+\ker\alpha\cap\gt n(A)])=
\ad^*(l)\gamma(\eta+\ker\alpha\cap\gt n(A))=c\gamma(\eta)=c\alpha(\eta)\ne 0,
$$
because $c$ is a non-zero constant.
Therefore $\bar\eta,\,l\not\in\gt q(A)_\gamma$ and
one can conclude that $\gt q(A)_\gamma=\gt q_\gamma$.
By the same reason, $\gt q(A)_\alpha=\gt t$.
In $\gt q(A)^*$, $\gamma$ is obtained from $\alpha$ by the same procedure as
in $\gt q^*$.  By the inductive hypothesis
$\gt q_\gamma$ is reductive.
\end{proof}

 Linear functions $\gamma\in\gt q^*$ such that  $\gamma(\gt q_\gamma)=0$ are said to be of {\it nilpotent type}. The name can be justified by the following
observation.

\begin{lm} \label{nilp-0}
Suppose that ${\rm char}\,\Bbbk=0$ and
$\gamma\in\gt q^*$ is of nilpotent type.
Then $\overline{Q\gamma}$ contains zero.
\end{lm}
\begin{proof}
Recall that $\ad(\gt q)^*\gamma=\Ann(\gt q_\gamma)$.
Since $\gamma(\gt q_\gamma)=0$, there is $\xi\in\gt q$ such that
$\ad(\xi)\gamma=\gamma$ or better $\ad(\xi)\gamma=-\gamma$.
Take $t\in\mathbb N\subset\Bbbk$ and consider
$$
\begin{array}{l}
\exp(t \ad (\xi))\gamma=\gamma+t\ad(\xi)\gamma+t^2\ad(\xi)^2\gamma/2+\ldots+
t^k\ad(\xi)^k\gamma/k!+\ldots= \\
\qquad
\gamma-t\gamma+t^2\gamma/2-t^3\gamma/3!+\ldots+(-1)^kt^k\gamma/k!+\ldots=
\exp(-t)\gamma.
\end{array}
$$
Since $\exp(-t)$ goes to zero when $t$ goes to infinity,
$0\in\overline{Q\gamma}$.
\end{proof}

Importance of linear forms, which are  simultaneously of
reductive and nilpotent type is illustrated by the following proposition (see \cite[Th\'eor\`eme 3.6.2]{DKT}).

\begin{prop}[Duflo-Khalgui-Torasso] \label{MRS-first}
Suppose that $\gt q$ is strongly quasi-reductive.
Then
\begin{itemize}
\item[({\sf i})] there is a unique orbit $Q\gamma\subset\gt q^*$ such that
$\gamma$ is of reductive and nilpotent type;
\item[({\sf ii})] if $\gt q_\beta$ is reductive, then
$\gt q_\beta$ is conjugate under $Q$ to a subalgebra
$(\gt q_\gamma)_s$ with $s\in \gt q_\gamma^*$.
\end{itemize}
\end{prop}
\begin{proof}
Existence of $\gamma$ was shown in Lemma~\ref{unipotent-exists}.
Suppose we have two linear functions $\gamma$,
$\gamma'$ of reductive and nilpotent type.
Replacing $\gamma'$ by an element of $Q\gamma'$ we may
(and will) assume that there is a generic stabiliser
$\gt t=\gt q_\alpha$ such that
$\gt t\subset (\gt q_\gamma\cap\gt q_{\gamma'})$.
Let $\gt h=\gt z_{\gt q}(\gt t)=\gt t{\oplus}\gt m_0$ and $\gt m_0$ be as in (\ref{decomp}).
Set $\gt v:=[\gt t,\gt q]$.
We have $\gamma(\gt v)=\gamma'(\gt v)=0$  
and  
both functions can be viewed as elements of $\gt h^*$.
Since $\gt t$ is a maximal torus of $\gt q_\gamma$,
we have $\gt h\cap\gt q_{\gamma}=\gt t$.
Now let $y\in\gt h$. Then $[y,\gt v]\subset\gt v$ and
$\gamma([y,\gt v])=0$.
Thereby $\gt h_\gamma=\gt h\cap\gt q_\gamma=\gt t$.
The same holds for $\gamma'$ and both functions are zero on $\gt t$.

Let $H\subset Q$ be a connected subgroup with $\Lie H=\gt h$.
Note that $\dim H\gamma=\dim\gt h-\dim\gt t$ and $H\gamma$
is a dense open subset of
$Y:=\{\xi\in\gt h^*\mid \xi(\gt t)=0\}$. The same holds for $H\gamma'$.
Since $Y$ is irreducible, $H\gamma\cap H\gamma'\ne\varnothing$ and
the orbits coincide. Therefore $\gamma'\in Q\gamma$.

Now part ({\sf i}) is proved and we pass to ({\sf ii}).
Let $\gt t=\gt q_\alpha$ be a maximal torus of $\gt q_\beta$ and $\gamma$ a linear form of reductive and nilpotent type
constructed  from $\alpha$ by means of (\ref{decomp}).
Since $B$ is non-degenerate on $\gt t$, it is also
non-degenerate on $\gt q_\beta$ and
$\gt q=\gt q_\beta\oplus \gt w$, where $\gt w=\gt q_\beta^{\perp}$.
We have $\gt q_\alpha\cap\gt h=\gt t=\gt q_\beta\cap \gt h$
and $H\beta$, $H\alpha$ are dense open subsets in
$\beta+\gt m_0^*$, $\alpha+\gt m_0^*$, respectively.
Replacing $\alpha$ by a conjugate function we may (and will)
assume that $\alpha|_{\gt m_0}=\beta|_{\gt m_0}$.
This implies that $\beta$ and $\gamma$ are equal on
$\gt m_0\oplus\gt v$.

Let $x\in\gt q_\beta$. Then
$[x,\gt q_\beta]\subset\gt q_\beta\subset\gt t\oplus\gt v$ and
$\gamma([x,\gt q_\beta])=0$.
In addition, $[x,\gt w]\subset\gt w\subset\gt m_0\oplus\gt v$ and
$\gamma([x,\gt w])=\beta([x,\gt w])=0$.
This proves the inclusion $\gt q_\beta\subset\gt q_\gamma$.
By a similar reason, $\beta([y,\gt q_\gamma^{\perp}])=0$ for
$y\in\gt q_\gamma$.
Therefore $\gt q_\beta=(\gt q_\gamma)_s$
for $s=\beta|_{\gt q_\gamma}$.
To conclude, note that by part
({\sf i}) there is only one orbit of both reductive and nilpotent type.
\end{proof}

Proposition~\ref{MRS-first} can be also deduced from \cite[Chap.~I]{Du}.

\begin{cl}[{\cite[Th{\'e}or{\`e}me 3.6.2]{DKT}}] \label{red-orbits}
Suppose that $\gt q$ is strongly quasi-reductive.
Let $\gamma\in\gt q^*$ be both of reductive and nilpotent type.
%
Then the  coadjoint orbits of reductive type  in $\q^*$
are in bijection with the closed (co)adjoint orbits of ${Q_{\gamma}}$.
In other words, the geometric quotient
$\gt q^*_{\rm red}/Q$ exists and coincides with
$\gt q_{\gamma}/\!\!/Q_\gamma$.
\end{cl}
\begin{proof}
By the proof of Proposition~\ref{MRS-first}({\sf ii}), each $Q$-orbit of reductive type
contains a point  $\beta=\gamma+s$ with $s\in\gt q_\gamma^*$ lying in the slice to
$Q\gamma$.
Let $\psi$ and $X_{\gamma}$ be as in the proof of Lemma~\ref{exist-gen-stab},
with $V=\gt q^*$.
We get that $Q\beta$ meets $\psi(X_\gamma)$.
Each connected component of $Q_\beta$ contains an element
$g$ preserving the maximal torus $\gt t$ and hence the decomposition~(\ref{decomp}).
Therefore $gs=s$ and $g\gamma=\gamma$.
This proves that $Q_\beta=(Q_\gamma)_s$ and the result follows.
\end{proof}

Corollary~\ref{red-orbits} justifies the following definition:

\begin{df} \label{d:max}
Suppose that $\q$ is strongly quasi-reductive
and $\gamma\in\gt q^*$ is a linear form of reductive type such that
$\gamma(\q_{\gamma})=0$. Then  $Q_\gamma$
is called a \emph{maximal reductive stabiliser} of $\q$.
As an abbreviation, we will  write MRS for a maximal reductive stabiliser.
\end{df}

Note that a maximal reductive stabiliser, as well as a generic stabiliser, is defined up to conjugation. 
Let $M_*(\q)$ denote the Lie algebra of a representative of the conjugation class of a MRS 
for a strongly quasi-reductive Lie algebra $\q$. 
For convenience, we set $M_*(\q)=\varnothing$ whenever $\q$ is not
strongly quasi-reductive.
Note  also that $M_*(\q)$ is defined up to conjugation.

\begin{rmk}\label{r:max}
If $\gt q$ is strongly quasi-reductive,
then the index of $\gt q$ is equal to the rank of $M_*(\gt q)$.
\end{rmk}

In this paper we deal with $M_*(\gt q)$ and
leave the description of the MRS on the group level for further investigation.
Following examples show that this problem is not entirely trivial.

\begin{ex}\label{group-mrs-1}
Let $Q=\cp^{^\times}{\ltimes}\exp(\cp^2{\oplus}\cp)$ be
a semi-direct product of a one-dimensional torus and
a Heisenberg Lie group.
Assume that $\cp^{^\times}$ acts on $\cp^2$ with characters
$(1,1)$, hence on the derived algebra of the Heisenberg
algebra
with character $2$.
Then $\ind\gt q=0$. Here MRS is equal to
$\{1,-1\}$. Since the centre of $Q$ is trivial, MRS will stay disconnected
after taking a quotient by $Z$.
\end{ex}

\begin{ex}\label{group-mrs-2}
Consider a semi-direct product
$Q=(\cp^{^\times}{\times}SO_9(\cp))\ltimes\exp(\cp^9)$,
where the central torus of the reductive part acts on
$\cp^9$ with character $1$.
Then MRS of $\gt q$ is equal to $O_8(\cp)$ and the component group
acts non-trivially on the set of coadjoint orbits of $M_*(\gt q)$.
\end{ex}

Motivated by the assertion of Corollary~\ref{red-orbits},
one is interested in a more precise description of $M_*(\gt q)$
for (strongly) quasi-reductive $\gt q$.

\section{On quasi-reductive biparabolic subalgebras} \label{S:sea}

In this section,
$\g=\Lie G$ is a finite-dimensional semisimple complex Lie algebra.
The dual of $\g$ is identified with $\g$ through the Killing form $\kappa$ of $\g$. 
For $u \in \g$, we denote by $\varphi_u$ the corresponding element of $\g^*$.
When the orthogonality in $\g$ refers to $\kappa$, we use the symbol $\perp$.

Recall that a {\it biparabolic subalgebra} of $\g$ is defined to be the intersection of two parabolic subalgebras
whose sum is $\g$.
They are also called {\it seaweed subalgebra} 
because of their shape in the case
of $\mathfrak{sl}_n$, see Picture~\ref{sl}.

\begin{figure}[htb]
{\setlength{\unitlength}{0.1in}
\begin{center}
\begin{picture}(14,14)(0,0)

\put(0,0){\line(1,0){14}}
\put(0,0){\line(0,1){14}}
\put(14,0){\line(0,1){14}}
\put(0,14){\line(1,0){14}}

\qbezier[75](0,14)(7,7)(14,0)

\put(0,10){\line(1,0){4}}
\put(4,10){\line(0,-1){4}}
\put(4,6){\line(1,0){4}}
\put(8,0){\line(0,1){6}}

\put(5,14){\line(0,-1){5}}
\put(5,9){\line(1,0){6}}
\put(11,3){\line(0,1){6}}
\put(14,3){\line(-1,0){3}}

\put(0,10){\line(1,4){1}}
\put(1,10){\line(1,4){1}}
\put(2,10){\line(1,4){1}}
\put(3,10){\line(1,4){1}}
\put(4,10){\line(1,4){1}}

\qbezier[200](4,6)(4.5,9.5)(5,13)
\qbezier[150](4.5,6)(4.75,8.5)(5,11)

\put(5,6){\line(1,3){1}}
\put(6,6){\line(1,3){1}}
\put(7,6){\line(1,3){1}}
\put(8,6){\line(1,3){1}}

\qbezier[200](8,3)(8.75,6)(9.5,9)
\qbezier[300](8,0)(9,4.5)(10,9)
\qbezier[300](8.5,0)(9.5,4.5)(10.5,9)
\qbezier[300](9,0)(10,4.5)(11,9)
\qbezier[250](9.5,0)(10.25,4)(11,8)

\put(10,0){\line(1,6){1}}

\qbezier[150](10.5,0)(10.75,1.5)(11,3)

\put(11,0){\line(1,3){1}}
\put(12,0){\line(1,3){1}}
\put(13,0){\line(1,3){1}}

\end{picture}
\end{center}}
\caption{}\label{sl}
\end{figure}

\subsection{}

Let $\q$ be a biparabolic subalgebra of $\g$ and we assume that $\q$ is quasi-reductive.
Let $Q$ be the connected Lie subgroup of $G$ with Lie algebra $\q$.
As it has been noticed in the Introduction, the centre of $\q$ consists of semisimple elements of $\g$.
Hence $\q$ is strongly quasi-reductive and results of Section~\ref{S:gen} apply.

The Killing form enables to identify the dual of $\q$ with
a biparabolic subalgebra $\q^-$ of $\g$.
Thus, to $\q_{\red}^{*}$ and $\q_{\reg}^{*}$ correspond subspaces of $\q^-$:
\begin{eqnarray*}
    && \q_{\reg}^- := \{u\in \q^-; \ (\var{u})_{|\q} \in \q_{\reg}^* \};\\
    && \q_{\red}^- := \{u\in \q^-; \ (\var{u})_{|\q} \in \q_{\red}^* \}.
\end{eqnarray*}

By Lemma~\ref{exist-gen-stab}, there is $x$ in $\q_{\reg}^- \cap \q_{\red}^-$.
Let $\alpha$ be the restriction of $\var{x}$ to $\q$.
The stabiliser  $Q_\alpha$ is an algebraic torus in $G$ and hence 
we have $\gt g=\gt t\oplus\t^{\perp}$ for $\gt t=\gt q_\alpha$.
Let  $x_\t$ and $x_{\t^\perp}$ denote the components of $x$ in
$\t$ and $\t^{\perp}$, respectively.
%

Obviously, $\t$ is contained in the stabiliser of $\var{x_{\t}}$ in $\q$ and hence also in the stabiliser  of $\var{x_{\t^\perp}}$. 
In the notation of Lemma~\ref{unipotent-exists},
we have $\alpha=s+\gamma+0$, where $s$ and $\gamma$ are the restrictions to $\q$ 
of $\var{x_{\t}}$ and $\var{x_{\t^\perp}}$, respectively. 
Hence, by Lemma~\ref{unipotent-exists}, Proposition~\ref{MRS-first}, and Definition~\ref{d:max}, we can claim:

\begin{thm} \label{t:perp}
Suppose that $\q$ is quasi-reductive and
 let $\gt t=\gt q_{\varphi_x}$ be a generic stabiliser with $x\in\q_-$. 
Then 
the restriction of $\var{x_{\t^\perp}}$ to $\q$  is of nilpotent and reductive type 
as a linear function on $\gt q$. 
In particular, its stabiliser in $Q$ is a MRS.
\end{thm}

Whenever the generic reductive stabiliser $\t$ can be explicitly computed,
Theorem~\ref{t:perp} provides a procedure to describe $M_*(\q)$ (see Section~\ref{ss:cas}).

%
\subsection{}
%
Let  $\Pi$ denote the set of simple roots with respect to a fixed triangular decomposition
$$\mathfrak{g}=\mathfrak{n}^+ \oplus \mathfrak{h} \oplus \mathfrak{n}^-,$$
and by $\Delta$ (respectively $\Delta^+$, $\Delta^{-}$)
the corresponding root system (respectively positive root system, negative root system).
If $\pi$ is a subset of $\Pi$, denote by $\Delta_\pi$ the root subsystem of $\Delta$ generated
by $\pi$ and set $\Delta_\pi^{\pm} := \Delta_\pi \cap \Delta^{\pm}$.
For $\alpha \in \Delta$, denote by $\mathfrak{g}_{\alpha}$ the $\alpha$-root subspace of $\g$ and let $h_\alpha$ be the unique element of
$[\mathfrak{g}_{\alpha},\mathfrak{g}_{-\alpha}]$ such that $\alpha(h_{\alpha})=2$.
For each $\alpha \in \Delta$, fix $x_{\alpha} \in \g_{\alpha}$ so that the family
$\{e_{\alpha}, h_{\beta} \ ; \  \alpha \in \Delta, \beta \in \Pi \}$ is a Chevalley basis of $\g$.

We briefly recall a classical construction due to B.~Kostant.
It associates to a subset of $\Pi$ a system of strongly orthogonal positive roots in $\Delta$.
This construction is known to be very helpful to obtain regular forms on biparabolic subalgebras of $\g$.
For a recent account about the \emph{cascade} construction of Kostant, we refer to~\cite[\S1.5]{TY2} or~\cite[\S40.5]{TY1}.

Recall that two roots $\alpha$ and $\beta$ in $\Delta$ are said to be {\it strongly orthogonal} if neither $\alpha+\beta$ nor $\alpha -\beta$ is in  $\Delta$.
Let $\pi$ be a subset of $\Pi$.
The {\em cascade} $\mathcal{K}_{\pi}$ of $\pi$ is defined by induction on the cardinality of $\pi$ as follows:
\begin{itemize}
\item[(1)] $\mathcal{K}_\varnothing=\varnothing$,
\item[(2)] If $\pi_1$,\ldots,$\pi_r$ are the connected components of $\pi$, then~
$\mathcal{K}_{\pi}=\mathcal{K}_{\pi_1} \cup \cdots \cup \mathcal{K}_{\pi_r}$,
\item[(3)] If $\pi$ is connected,  then $\mathcal{K}_{\pi}=\{\pi\}
\cup \mathcal{K}_{T}$ where $T$ is the set of simple roots that are orthogonal to the highest positive root
$\theta_{\pi}$ of $\Delta_{\pi}^+$.
\end{itemize}
For  a subset $\pi\subset\Pi$, let $\cal{E}_{\pi}$  denote the set of
the highest roots $\hr{K}$, where $K$ runs over the elements of the cascade of $\pi$.
The cardinality of $\mathcal{K}_{\Pi}$  depends only on $\mathfrak{g}$;
it is independent of the choices of $\mathfrak{h}$ and $\Pi$.
Denote it by $\kg{\g}$.
The values of $\kg{\g}$
for the different types of simple Lie algebras are given in Table~\ref{Tkg};
in this table, for a real number $x$, we denote by $[x]$ the largest integer
$\le x$.

{
\begin{table}[h] 
\begin{center}
\begin{tabular}{|c|c|c|c|c|c|c|c|c|}
\hline
  & & & & & & & & \\
$\mathrm{A}_{\ell}, \ell \ge 1$  &
$\mathrm{B}_{\ell}, \ell \ge 2$  &
$\mathrm{C}_{\ell}, \ell \ge 3$  &
$\mathrm{D}_{\ell}, \ell \ge 4$  &
$\mathrm{G}_2$  &
$\mathrm{F}_4$  &
$\mathrm{E}_6$  &
$\mathrm{E}_7$  &
$\mathrm{E}_8$   \\
  & & & & & & & & \\
\hline
 & & & & & & & & \\
$\left[ \displaystyle{\frac{\ell+1}{2}} \right]$ & $\ell$ & $\ell$ &
$2 \left[ \displaystyle{\frac{\ell}{2}} \right]$ &
$2$ & $4$ & $4$ & $7$ & $8$ \\
  & & & & & & & & \\
\hline
\end{tabular}
\vspace{.1cm}
\caption{\label{Tkg} $\kg{\g}$ for the simple  Lie algebras.}
\end{center}
\end{table}}

We denote by $\mathfrak{p}_{\pi}^{+}$ the standard parabolic subalgebra of $\mathfrak{g}$, which  is the subalgebra generated by $\mathfrak{b}^+=\mathfrak{h}
\oplus \mathfrak{n}^+$ and by $\mathfrak{g}_{-\alpha}$, for $\alpha \in \pi$.
It is well-known that any parabolic subalgebra of $\g$ is conjugate to a standard one.
We denote by $\mathfrak{p}_{\pi}^{-}$ the ``opposite parabolic subalgebra'' generated by $\mathfrak{b}^- = \mathfrak{n}^- \oplus \mathfrak{h}$ and
by $\mathfrak{g}_{\alpha}$, for $\alpha \in \pi$.
Let
$\mathfrak{l}_{\pi}=\mathfrak{p}_{\pi}^{+} \cap \mathfrak{p}_{\pi}^{-}$ be the standard
Levi factor of both $\mathfrak{p}_{\pi}^{+}$
and $\mathfrak{p}_{\pi}^{-}$.
We denote by $\g_{\pi}$ its derived Lie algebra.
Let $\mathfrak{m}_{\pi}^{+}$ (respectively $\mathfrak{m}_{\pi}^{-}$) be the nilradical of $\mathfrak{p}_{\pi}^{+}$
(respectively $\mathfrak{p}_{\pi}^{-}$).
In the sequel, we will make use of the following element of $\p_{\pi}^-$:
\begin{eqnarray*}
u_{\pi}^- &=& \sum\limits_{\hr{} \in \cal{E}_{\Pi}\mbox{, }\hr{}\notin \Delta_{\pi}^{+}} e_{-\hr{}}.
\end{eqnarray*}

%
\subsection{}
%

In~\cite[Section 2]{BM}, the authors described an additivity property for the quasi-reductivity of certain parabolic subalgebras of $\g$.
The additivity property runs particularly smoothly when $\rk\g=\kg{\g}$. 
For maximal reductive stabilisers as well, we are going to prove an additivity property for certain quasi-reductive parabolic subalgebras of $\g$. 
We start by recalling some definitions and results of~\cite{BM}. 

Until the end of this section, we assume that $\g$ is simple and that the condition $\rk\g=\kg{\g}$ is satisfied. 

\begin{rmk}\label{not-simple}
If $\g=\g_1\oplus\g_2$ is not simple, then
$\gt q=\gt q_1\oplus\gt q_2$ with $\gt q_i\subset\g_i$
and $M_*(\gt q)=M_*(\gt q_1)\oplus M_*(\gt q_2)$.
Therefore for the description of $M_*(\gt q)$ it suffices to
consider only simple Lie algebras.
\end{rmk}

\begin{df} \label{d:s}
Let $\pi',\pi''$ be subsets of $\Pi$.
We say that $\pi'$ is {\em not connected to} $\pi''$ if $\alpha'$ is orthogonal to $\alpha''$,
for all $(\alpha',\alpha'')$ in $\pi'\times\pi''$.

For $\alpha \in \Delta^+$, we define an element $K_{\Pi}^+(\alpha)$ of $\mathcal{K}_\Pi$ as follows: 
if $\alpha$ does not belong to $\mathcal{E}_\Pi$, then $K_{\Pi}^+(\alpha)$
is the unique element $L$ of $\cal{K}_{\Pi}$ such that 
$\theta_{L}-\alpha$ is a positive root;
if $\alpha=\theta_K$ with $K\in\mathcal{K}_\Pi$, then we set $K_{\Pi}^+(\alpha):=K$.
\end{df}

It can be checked that $K_{\Pi}^+(\alpha)=K_{\Pi}^+(\beta)$ for simple roots
$\alpha$ and $\beta$ if and only if $\alpha$ and $\beta$ lie in the same orbit of $-w_0$, where $w_0$ is the longest element in
the Weyl group of $\g$. 
Since $\rk\g=\kg{\g}$, the element $-w_0$ is trivial, and it results from properties of $\mathcal{K}_\Pi$ 
that two subsets $\pi',\pi''$ of $\Pi$   
which are not connected to each other satisfy the condition: 
$K_{\Pi}^{+}(\alpha')\not= K_{\Pi}^{+}(\alpha'')$, $\forall\ (\alpha',\alpha'')\in \pi' \times \pi''$.  
By the properties of $\mathcal{K}_\Pi$, this condition is equivalent to the condition $(\clubsuit)$:
\begin{eqnarray*}
(\clubsuit)
&& K_{\Pi}^{+}(\alpha')\not= K_{\Pi}^{+}(\alpha'')\quad \forall\ (\alpha',\alpha'')\in \Delta_{\pi'}^+ \times \Delta_{\pi''}^+\, .
\end{eqnarray*}

Let $\pi',\pi''$ be two subsets of $\Pi$ which are not connected to each other. 
Set $\b^+:=\h\oplus\n^+$ and let $B$ be the corresponding Borel subgroup of $G$, i.e.,  $\Lie B=\b^+$. 
By a result of Kostant, for any ideal $\m$ of $\b^+$ contained in $\n^+$, 
the $B$-orbit of $\varphi_{u_\varnothing^-}$ in $\m^*$ is an open dense subset of $\m^*$. 
Hence there is $w'$ in $\l_{\pi'}$ such that the restriction of $\varp{w}$ to $\p_{\pi'}^+$  is regular, 
where $w=w'+u_{\pi'}^-$. 
Let $\t'$ be the stabiliser  in $\p_{\pi'}^+$ of 
the restriction of $\varp{w}$ to $\p_{\pi'}^+$.

\begin{lm} \label{l:add1} Keep the above notation. 
Then  $\t'$ is contained in the derived algebra 
$[\p_{\pi'}^+,\p_{\pi'}^+]=\g_{\pi'}{\oplus}\m_{\pi' \cup\pi''}^+$. 
In particular, $[\t', \p_{\pi'\cup \pi''}^+] \subset \p_{\pi'}^+$.
\end{lm}

\begin{proof}
Let  $\s'$ denote  the image of $\t'$ under the projection map from
$\p_{\pi'}^+$ to $\g_{\pi'}{\oplus}\m_{\pi'}^+$ along the centre $\z(\l_{\pi'})$ of
$\l_{\pi'}$  (and of $\p_{\pi'}^+$).

Let $x$ be an element of $\t'$. 
Write $x=x_0+x'+x^+$ with $x_0\in \z(\l_{\pi'})$, $x'\in\g_{\pi'}$, and $x^+ \in \m_{\pi'}^+$.
Since $[x^+,w']$ lies in $\m_{\pi'}^+$, the fact that $x \in \t'$ means
$[x_0,u_{\pi'}^-]+[x',w'] +[x',u_{\pi'}^-] +[x^+,u_{\pi'}^-] \in \m_{\pi'}^+$. 
From this, we first deduce that $x_0$ is in $\ker \theta$ for all 
$\theta \in \mathcal{E}_{\Pi}\setminus\Delta_{\pi'}^+$. 
On the other hand, $\theta(x_0)=0$ for any $\theta \in \mathcal{E}_{\Pi}\cap \Delta_{\pi'}^+$ 
since $x_0$ is in the centre of $\l_{\pi'}$. 
Our assumption $\rk\g=\kg{\g}$ tells us that $\mathcal{E}_{\Pi}$ is a basis of $\h^*$, so $x_0=0$. 

It remains to show that $x^{+} \in \m_{\pi'\cup \pi''}^+$.
If not, there are $\gamma \in \Delta_{\pi''}^+$, $K \in \cal{K}_{\Pi}$, and $\alpha' \in \Delta_{\pi'}^+$  such that
$\gamma - \theta_{K_{\Pi}^{+}(\gamma)} = - (\alpha' + \theta_{K})$, 
i.e.~$\theta_{K_{\Pi}^{+}(\gamma)} = \gamma + (\alpha' + \theta_{K})$ 
whence $K_{\Pi}^{+}(\alpha')=K_{\Pi}^{+}(\gamma)$.
But this contradicts condition $(\clubsuit)$.
\end{proof}

Let $\pi',\pi''$ be two subsets of $\Pi$ which are not connected to each other as before.  
Assume besides that $\p_{\pi'}^+$ and $\p_{\pi''}^+$ are both quasi-reductive and set $\pi:=\pi'\cup\pi''$. 
Then $\p_\pi^+$ is quasi-reductive by~\cite[Theorem 2.11]{BM}, 
and there is $(w_0,w',w'')\in \z(\l_{\pi' \cup\pi''}) \times \g_{\pi'} \times  \g_{\pi''}$
such that the restriction of $\var{v}$ to $\p_{\pi'}^+$ (resp.~$\p_{\pi''}^+$ and $\p_{\pi}^+$) 
is regular and of reductive type, where $v=w_0 + w' + w'' + u_{\pi}^-$. 
Let $\t'$, $\t''$, and $\t$ be the respective stabilisers 
of the restrictions of $\var{v}$ to $\p_{\pi'}^+$, $\p_{\pi''}^+$, and $\p_{\pi}^+$. 
Let also $\r'$, $\r''$, and $\r$ be the respective stabilisers
of the restrictions of $\var{v_{\t^\perp}}$ to $\p_{\pi'}^+$, $\p_{\pi''}^+$, and $\p_{\pi}^+$. 

\begin{lm} \label{l:add2}
{\sf (i)} We have $\t=\t'\oplus\t''$. 

{\sf (ii)} $\r'$ is contained in $\g_{\pi'}\oplus \m_{\pi}^+$ and 
$\r''$ is contained in $\g_{\pi''}\oplus \m_{\pi}^+$. 

\end{lm}

\begin{proof} 
{\sf (i)} By~\cite[\S2.2, equality (2)]{BM}, we have $\ind \p_\pi^+=\ind\p_{\pi'}^++\ind\p_{\p''}^+$ 
since $\rk\g=\kg{\g}$, that is, $\dim \t=\dim\t'+\dim\t''$.
By Lemma~\ref{l:add1}, 
$\t' \subset \g_{\pi'}\oplus\m_\pi^+$ and $\t'' \subset \g_{\pi''}\oplus\m_\pi^+$. 
This forces $\t=\t'\oplus\t''$ because $\t'$ and $\t''$ are both consisted of semisimple elements of $\g$. 

{\sf (ii)} By Theorem~\ref{t:perp}, $\r'$ is a maximal reductive stabiliser of $\p_{\pi'}^+$. 
Moreover, $\t'$ is a maximal torus of $\r'$ and we can write $\r'=\t' \oplus [\t',\r']$.
So, the statement follows from Lemma~\ref{l:add1},    
since $\g_{\pi'}\oplus \m_{\pi}^+$ is an ideal of $\p_{\pi'}^+$. 
The same goes for $\r''$.
\end{proof}

\begin{thm}[Additivity property] \label{t:add}
Assume that $\g$ is simple and $\rk\g=\kg{\g}$.
Let $\pi',\pi''\subset \Pi$ be two subsets which are not connected to each other and
assume that $\p_{\pi'}^+$ and $\p_{\pi''}^+$ are both quasi-reductive. 
Set $\pi:=\pi'\cup\pi''$. 
Then $\p_{\pi}^+$ is quasi-reductive and $M_*(\p_{\pi'}^+) \oplus M_*(\p_{\pi''}^+) = M_*(\p_{\pi}^+)$.
More precisely, there exist a maximal reductive stabiliser $\r'$ of $\p_{\pi'}^+$
and a maximal reductive stabiliser $\r''$ of $\p_{\pi''}^+$ such that $\r' \oplus \r''$ is a maximal reductive stabiliser of 
$\p_{\pi}^+$.
\end{thm}

\begin{proof}
We need to prove only the second part of the statement. 
Fix a triple $(w_0,w',w'') \in \z(\l_{\pi}) \times \g_{\pi'} \times  \g_{\pi''}$
such that for $v=w_0 + w' + w'' + u_{\pi}^-$,
the restriction of $\var{v}$ to $\p_{\pi'}^+$ (resp.~$\p_{\pi''}^+$ and $\p_{\pi}^+$)
is regular and of reductive type; and keep the above notations.  
Our goal is to show that $\r=\r'\oplus\r''$.

\medskip

We prove 
the inclusion $\r' \subseteq \r$; the proof of the inclusion $\r'' \subseteq \r$ is similar.
From the equality $\kappa(v_{\t'^{\perp}},[\r',\p_{\pi'}^+])=\{0\}$, we deduce that
$\kappa(v_{\t'^{\perp}},[\r',\p_{\pi}^+])=
\kappa(v_{\t'^{\perp}},[\r',\g_{\pi''}])$.
In turn, by Lemma~\ref{l:add2}{\sf (ii)}, $[\r',\g_{\pi''}]$ is contained in $\m_{\pi}^+$,
whence $\kappa(v_{\t'^{\perp}},[\r',\p_{\pi}^+]) \subset \kappa(v_{\t'^{\perp}},\m_{\pi}^+)$.
Now, $\kappa(v_{\t'^{\perp}},\m_{\pi}^+)=\{0\}$, since $\t'^{\perp} \supset \m_{\pi}^+$.
To sum up, we have obtained:
\begin{eqnarray} \label{e:add}
\kappa( v'_{\t'^{\perp}},[\r',\p_{\pi}^+] ) &=& \{0\}.
\end{eqnarray}

It is worth noting that $\kappa$ is non-degenerate on $\t'\times\t'$ (resp.~$\t''\times\t''$ and $\t\times\t$) 
and on $\t'^\perp\times\t'^\perp$ (resp.~$\t''^\perp\times\t''^\perp$ and $\t^\perp\times\t^\perp$). 
Since $\t' \subset \g_{\pi'}\oplus \m_{\pi}^+$, we get $w'' \in \t'^{\perp}$ and $w_0 \in \t^{\perp}=\t'^\perp\cap\t''^\perp$.
Thus, $v_{\t'^{\perp}}=w_0+w'_{\t'^{\perp}} + w'' + (u_{\pi}^-)_{\t'^\perp}$ 
and $v_{\t^{\perp}}=w_0+w'_{\t'^{\perp}}+w''_{\t''^{\perp}}+ (u_{\pi}^-)_{\t'^\perp \cap \t''^\perp}$.
So, $v_{\t^{\perp}} - v_{\t'^{\perp}} = w''_{\t''}  + ((u_{\pi}^-)_{\t'^\perp})_{\t''}$.
In particular, since $\t''$ is contained in $\g_{\pi''} \oplus \m_{\pi}^+$ by Lemma~\ref{l:add1}, 
we see that $v_{\t^{\perp}}- v_{\t'^\perp}$ lies in $\g_{\pi''} \oplus \m_{\pi}^+$.
Hence, from the inclusion $\r' \subset \g_{\pi'}\oplus \m_{\pi}^+$ (Lemma~\ref{l:add2}{\sf (i)}), 
we deduce that
$[v_{\t^{\perp}}- v_{\t'^{\perp}},\r']$ is contained in $\m_{\pi}^+$.
We further deduce:
\begin{eqnarray} \label{e2:add}
\kappa(v_{\t^{\perp}}- v_{\t'^{\perp}},[\r',\p_{\pi}^+]) &=& \{0\}.
\end{eqnarray}
Finally, it follows from~(\ref{e:add}) and~(\ref{e2:add}) that $\kappa(v_{\t^{\perp}},[\r',\p_{\pi}^+])=\{0\}$,
whence $\r' \subset \r$.

\medskip

Turn now to the inclusion $\r \subseteq \r' + \r''$.

\medskip

\noindent
Step 1: By what foregoes, $v_{\t'^{\perp}}=v_{\t^{\perp}} + y$, where $y$ is an element of $\t''\subset \g_{\pi''} \oplus \m_{\pi}^+$. 
On one hand, we have $\kappa(y,[\g_{\pi'}\oplus \m_{\pi}^+,\p_{\pi'}^+])=\{0\}$. 
On the other hand, $\kappa(v_{\t^{\perp}},[\r,\p_{\pi'}^+])=\{0\}$, since $\p_{\pi'}^+ \subset \p_\pi^+$. 
From this, we deduce that $\r \cap (\g_{\pi'}\oplus \m_{\pi}^+)$ is contained in $\r'$.
Similarly, we show that $\r \cap (\g_{\pi''}\oplus \m_{\pi}^+)$ is contained in $\r''$.

\medskip

\noindent
Step 2: Denote by $\a'$ (resp.~$\a''$) the sum of all simple factors of $\r$ that are contained in $\g_{\pi'}\oplus \m_{\pi}^+$
(resp.~$\g_{\pi''}\oplus \m_{\pi}^+$), and denote by $\a_0$ the sum of all simple factors
that are neither contained in $\g_{\pi'}\oplus \m_{\pi}^+$ nor in $\g_{\pi''}\oplus \m_{\pi}^+$.
Thus, $\r=\a'\oplus \a''\oplus \a_0$.
By Step 1, $\a' \subset \r'$ and $\a''\subset \r''$.
Furthermore, the inclusions $\r' \subset \r$ and $\r''\subset \r$, shown just above, imply $\r' \subset \a'$ 
and $\r''\subset \a''$ by Lemma~\ref{l:add2}{\sf (ii)}. 
So, $\a' = \r'$ and $\a''= \r''$. 
Hence $\r'$ and $\r''$ are in direct sum and $\t'\oplus \t''$ is a maximal torus of $\a'\oplus \a''$.
But, $\t=\t' \oplus \t''$ is a maximal torus of $\r$.
This forces $\a_0=\{0\}$, whence $\r =\r' \oplus \r''$.
\end{proof}

\subsection{} \label{ss:cas}
We conclude this section with one remark about the use of Kostant's cascade in the description of MRS.

Let $\q:=\p_{\pi_1}^+ \cap \p_{\pi_2}^-$ be a quasi-reductive biparabolic subalgebra of $\g$.
In most cases, there is a tuple $(a_K)_{K \in \mathcal{E}_{\pi_1}} \cup (b_L)_{L \in \mathcal{E}_{\pi_2}}\in(\C^*)^{\kg{\pi_2}+\kg{\pi_2}}$
such that the element
$x=\sum\limits_{\theta \in \mathcal{E}_{\pi_1}} a_K e_{\theta} +
 \sum\limits_{\nu \in \mathcal{E}_{\pi_2}} b_L e_{-\nu}$ is in $\q_{\reg}^- \cap \q_{\red}^-$ (see \cite{BM}).
In~\cite[Appendix~A]{BM}, it is explained how to use computer programs like \texttt{GAP} to compute $\t$,
the stabiliser in $\q$ of $\var{x_\t}$, in such a case.
Then \texttt{GAP} gives the orthogonal complement to $\t$ in $\g$ and can compute the stabiliser of $\var{x_{\t^{\perp}}}$ in $\q$.
In this paper, we do not use \texttt{GAP} to describe the maximal reductive stabilisers.
However, some of the results in the exceptional case may be verified thank to it.

\begin{rmk}
As the computations of~\cite{BM} show, there is always an element $x$ of the form $x=\sum\limits_{\theta \in \mathcal{E}_{\pi_1}} a_K e_{\theta} +
 \sum\limits_{\nu \in \mathcal{E}_{\pi_2}} b_L e_{-\nu}$ which belongs to $\q_{\reg}^- \cap \q_{\red}^-$ in all cases of quasi-reductive parabolic subalgebras (with $\pi_2=\Pi$) in the exceptional cases.
\end{rmk}

\section{Reductions in graded Lie algebras and some parabolic subalgebras} \label{S:red}

In order to compute $M_*(\gt q)$, one may use approach of
Lemma~\ref{unipotent-exists} and ``cut'' the nilpotent radical of $\gt q$.
In practise this is a rather complicated task. Here we present some algorithms
for cutting small pieces of the nilpotent radical.

\subsection{}
In this section $Q$ is a linear algebraic group and $\gt q=\Lie Q$.

For a linear function $\alpha$ on a Lie algebra $\gt q$,
let $\hat\alpha$ denote the skew-symmetric form on $\gt q$
(or any of its subspaces) given by $\hat\alpha(\xi,\eta)=\alpha([\xi,\eta])$.
Let $\gt a\lhd \gt q$ be an ideal.
Then the Lie algebra $\gt q$ acts
on $\gt a$ and also on $\gt a^*$.
For $\gamma\in\gt a^*$,
$\gt q_\gamma$ will always refer to the stabiliser in $\q$ for this action. If $\gamma$ is extended
to a linear function on $\gt q$, this is explicitly stated
and the extension is denoted by some other symbol. 
Suppose that $\gt a\lhd\gt q$ is an Abelian ideal and
$\gamma\in\gt a^*$.  Let $\gt q_{\cp\gamma}$ be
the normaliser of the line $\cp\gamma$  in $\gt q$.
Then $\gt a\subset\gt q_{\cp\gamma}$. Let
$\q(\gamma)$ denote the quotient Lie algebra
$\gt q_{\cp\gamma}/\ker\gamma$, where
$\ker\gamma\subset\gt a$ is an ideal of $\gt q_{\cp\gamma}$.
(It is not always an ideal in $\gt q$.)

\begin{lm}\label{restriction}
Let $\gt a\lhd\gt q$ be an Abelian ideal,
$\beta\in\gt q^*$, $\gamma=\beta|_{\gt a}$
and $\check\beta=\beta|_{\gt q_\gamma}$.
Then
$(\gt q_\gamma)_{\check\beta}=\gt q_\beta +\gt a$.
\end{lm}
\begin{proof}
Since $[\gt a,\gt a]=0$, we have $\gt a\subset\gt q_\gamma$.
Besides
$\check\beta([\gt a,\gt q_\gamma])=\gamma([\gt a,\gt q_y])=0$.
Hence $\gt a\subset(\gt q_\gamma)_{\check\beta}$.
Next, $\gt q_\beta\subset\gt q_\gamma$, because $\gt a$ is an ideal,
and obviously  $\gt q_\beta\subset(\gt q_\gamma)_{\check\beta}$.
Now take $\xi\in(\gt q_\gamma)_{\check\beta}$. Then
$\ad^*(\xi)\beta\in(\gt q/\gt q_\gamma)^*$. On the other
hand, $\hat\beta$ defines a pairing between $\gt a$ and $\gt q$,
which is non-degenerate on $\gt q/\gt q_\gamma$ and $\gt a/\gt a_\beta$.
Hence
$$
\dim(\ad^*(\gt a)\beta)=\dim\gt a-\dim\gt a_\beta=
\dim\gt q-\dim\gt q_\gamma.
$$
Since all elements of $\ad^*(\gt a)\beta$ are zero on $\gt q_\gamma$,
we obtain
$\ad^*(\gt a)\beta=(\gt q/\gt q_\gamma)^*$
and
there is an element $\eta\in\gt a$ such that
$\ad^*(\xi)\beta=\ad^*(\eta)\beta$.
We have $\xi-\eta\in\gt q_\beta$ and hence
$\xi\in\gt q_\beta+\gt a$.
\end{proof}


\begin{lm}\label{ideal}
Assume that the centre of $\gt q$ consists of semisimple
elements.
Let $A\lhd Q$ be a non-trivial normal Abelian unipotent subgroup and
$\gt a=\Lie A$. Then
$\gt q$ is quasi reductive if and only if
the action of
$Q/A$ on $\gt a^*$ has an open orbit and, for generic $\gamma\in\gt a^*$,
$\gt q(\gamma)$ is strongly quasi-reductive.

\noindent
Moreover, if $\q$ is quasi-reductive, then
$M_*(\gt q)$  coincides with $M_*(\gt q(\gamma))$.
\end{lm}
\begin{proof}
Suppose  first that $\gt q$ is quasi-reductive.
Then $\gt q_\beta$ is reductive for some
$\beta\in\gt q^*$. By assumptions, the ideal $\gt a$ is contained in
the nilpotent radical of $\gt q$. Hence $\gt a\cap\gt q_\beta=0$.
Set $\gamma:=\beta|_{\gt a}$. Since $[\gt q,\gt a]\subset\gt a$,
the functions $\gamma$ and $\beta$ define the same pairing
$\hat\beta|_{\gt q{\times}\gt a}$
between $\gt q$ and $\gt a$, which is non-degenerate on $\gt a$,
since $\ker\hat\beta\cap\gt a=\gt a\cap\gt q_\beta=0$.
On the other hand, $\gt q_\gamma=\{y\in\gt q\mid \hat\beta(y,\gt a)=0\}$.
Hence $\dim Q\gamma=\dim\gt a$ and $Q\gamma$ is a required open orbit.
Since $A$ acts trivially on $\gt a^*$, $Q\gamma$ coincides with the
$Q/A$-orbit of $\gamma$.
Let $\check\beta$ be the restriction of $\beta$ to $\gt q_\gamma$.
According to Lemma~\ref{restriction}, we have
$(\gt q_\gamma)_{\check\beta}=\gt q_\beta\oplus\gt a$,
in our case the sum is direct.
In  $\gt q_{\cp\gamma}$
the stabiliser of (the restriction of) $\beta$ is almost the same, more precisely,
it is equal to $\gt q_\beta\oplus\ker\gamma$, since there is an
element multiplying $\gamma$ by a non-zero number and
$\xi\in\gt a$ such that $\gamma(\xi)\ne 0$ does not stabilise $\beta$.
Hence, if we consider $\beta$ is a function on
$\gt q(\gamma)$, the stabiliser $\gt q(\gamma)_\beta$
is equal to $\gt q_\beta$ and therefore is reductive.
Thus $\gt q(\gamma)$ is quasi reductive and its
centre necessary consists of semisimple elements.
We have proved that each reductive stabiliser in $\gt q$
(of a linear function on $\gt q$) is also
a stabiliser in $\gt q(\gamma)$.

Now suppose that $\gamma\in\gt a^*$ belongs to an open
$Q/A$-orbit. Then $\gamma$ is generic. Suppose also
that $\gt q(\gamma)$ is quasi reductive and its centre
consists of semisimple elements. Then there is
$\beta\in\gt q(\gamma)^*$ such that
$\gt q(\gamma)_\beta$ is reductive.
If $\beta$ were zero on $\gt a/\ker\gamma$, its stabiliser
would have contained this quotient, and hence a
non-zero nilpotent element. Thereby the restriction
of $\beta$ to $\gt a/\ker\gamma$ is non-zero and
rescaling $\beta$ if necessary we may assume that it
coincides with $\gamma$. Let $\tilde\beta\in\gt q^*$
be a lifting of $\beta$ such that $\tilde\beta|_{\gt a}=\gamma$.
Let also $\check\beta$ be the restriction of $\tilde\beta$ to
$\gt q_\gamma$.
By Lemma~\ref{restriction},
$(\gt q_\gamma)_{\check\beta}=\gt q_{\tilde\beta}+\gt a$.
Since $\gamma$ belongs to an open $Q$-orbit,
$\gt a\cap\gt q_{\tilde\beta}=0$ and, hence,
$\gt q_{\tilde\beta}$ is isomorphic to $\gt q(\gamma)_\beta$.
In particular, $\gt q$ is quasi-reductive. Also
each reductive stabiliser in $\gt q(\gamma)$
(of a linear function on $\gt q(\gamma)$) is also a stabiliser in $\gt q$.
This proves that maximal reductive stabilisers of
two algebras coincide.
\end{proof}

\begin{lm}\label{2-grad}
Suppose that $\gt q=\gt q(0)\ltimes\gt a$ is a semi-direct
product of an Abelian ideal $\gt a$ 
consisted of nilpotent elements
and a Lie subalgebra $\gt q(0)$.
Assume further that the
centre of $\gt q$ consists of semisimple
elements.
Denote by $Q(0)$ the connected subgroup of $Q$ of Lie algebra $\q(0)$.
Then $\gt q$ is quasi reductive if and only if the action of
$Q(0)$ on $\gt a^*$ has an open orbit and
for generic $\gamma\in\gt a^*$, $\gt q(0)_\gamma$ is strongly quasi-reductive.

\noindent
Moreover, if $\q$ is quasi reductive, then
$M_*(\gt q)$  coincides with $M_*(\gt q(0)_\gamma)$.
\end{lm}
\begin{proof}
According to Lemma~\ref{ideal}, all the statements become true
if we replace $\gt q(0)_\gamma$ by $\gt q(\gamma)$.
In our case, there is a complementary to $\gt a$ subalgebra,
namely $\gt q(0)$. This implies that
$\gt q(\gamma)=\gt q(0)_\gamma\oplus\cp h\oplus\cp\xi$,
where $\gt q(0)_\gamma$ is a subalgebra commuting with
$\xi$, $\cp\xi$ is a commutative ideal contained in the nilpotent radical,
and $h$ acts on
$\cp\xi$ via a non-trivial character. Whenever a stabiliser
$\gt q(\gamma)_\beta$ with $\beta\in\gt q(\gamma)^*$ is reductive,
$\beta$ is necessary non-zero on $\xi$  and
$\gt q(\gamma)_\beta=(\gt q(0)_\gamma)_{\check\beta}$ for
$\check\beta$ being a restriction of $\beta$ to $\gt q(0)_\gamma$.
This completes the proof.
\end{proof}

\begin{df}\label{m-grading}
We will say that a Lie algebra $\gt q$ is $m$-graded, if
it has a $\mathbb Z$-grading with only $m$ non-trivial components:
$\gt q=\bigoplus\limits_{i=0}^{m{-}1} \gt q(i)$.
\end{df}

\begin{rmk}
A semi-direct product structure $\gt q(0)\ltimes\gt a$ can be also
considered as a $2$-grading
$\gt q=\gt q(0)\oplus\gt q(1)$, where $\gt q(1)=\gt a$.
In the same spirit, if $\gt q$ is
$m$-graded, then $\gt q(m{-}1)$ is an Abelian ideal.
\end{rmk}

\begin{lm}\label{3-grad}
Assume that the centre of $\gt q$ consists of semisimple elements.
Suppose that $\gt q$ is $3$-graded
$\gt q=\gt q(0)\oplus\gt q(1)\oplus\gt q(2)$ and
for generic $\alpha\in\gt q(2)^*\subset\gt q^*$
the skew-symmetric form
$\hat\alpha$ is non-degenerate on $\gt q(1)$.
Set $\gt a=\gt q(2)$.
Then $\gt q$ is quasi reductive if and only if the action of
$Q(0)$ on $\gt a^*$ has an open orbit and
for generic $\gamma\in\gt a^*$, $\gt q(0)_\gamma$ is strongly quasi reductive.

\noindent
Moreover, if $\q$ is quasi reductive, then
$M_*(\gt q)$  coincides with $M_*(\gt q(0)_\gamma)$.
\end{lm}
\begin{proof} First of all, if $\gt q$ is
quasi reductive, then by Lemma~\ref{ideal},
$Q$ has an open orbit in $\gt a^*$.
Since the normal subalgebra $\gt q(1)\oplus\gt q(2)$
acts on $\gt a^*$ trivially, that open orbit is
also an open orbit of $Q(0)$. For the rest of the proof
assume that $\overline{Q(0)\gamma}=\gt a^*$.
If $\gt q_\beta$ is reductive for some $\beta\in\gt q^*$, then
$\beta|_{\gt a}$ lies in $Q(0)\gamma$. Replacing $\beta$ by a conjugate
one, we may (and will) assume that the restriction of $\beta$ to
$\gt a$ equals $\gamma$. By our assumptions,
$\hat\beta$ is non-degenerate on $\gt q(1)$.
Again replacing $\beta$ by a $Q$-conjugate
linear function, one may assume that $\beta$ is zero on $\gt q(1)$.
Summing up, if we are interested in the existence of reductive
stabilisers and their possible types it suffices to consider
$\beta\in\gt q^*$ such that $\beta|_{\gt q(1)}=0$ and
$\beta|_{\gt a}=\gamma$. In this case
$\beta$ can be considered as a function on
$\hat{\gt q}=\gt q(0)\ltimes\gt q(2)$ and
$\gt q_\beta=\hat{\gt q}_\beta$.
Now all the claims follow from Lemma~\ref{2-grad}
applied to $\hat{\gt q}$.
\end{proof}

%
Let $\gt g$ be a simple Lie
algebra and $e\in\gt g$ a minimal nilpotent element.
It can be included into an $\gt{sl}_2$-triple
$\left( e,h,f\right)$ in $\gt g$. Then the $h$-grading
of $\gt g$ looks as follows
$$
\gt g=\gt g(-2)\oplus\gt g(-1)\oplus\gt g(0)\oplus\gt g(1)\oplus\gt g(2),
$$
where $\gt g(2)=\cp e$ and $\gt g(-2)=\cp f$.
Set $\gt p(e):=\gt g(0)\oplus\gt g(1)\oplus\gt g(2)$.
It is a parabolic subalgebra in $\gt g$ and
it is $3$-graded. The Levi part,
$\gt g(0)$, is generated by all simple roots orthogonal
to the highest root, the semisimple element $h$, and, in type A only, another
semisimple element.
Note that $\gt g(0)_e=\gt g(0)\cap\gt g_e$ is a reductive Lie subalgebra of $\gt g$.

\begin{prop}\label{min-n-red}
Let $\gt p\subset\gt p(e)$ be a parabolic subalgebra of $\g$.
Then $M_*(\gt p)=M_*(\gt g(0)_e\cap\gt p)$
and the conjugacy class of $M_*(\gt p)$ has a representative embedded into the intersection on the right hand side.
\end{prop}
\begin{proof} Set $\hat{\gt p}=\gt p\cap\gt g(0)_e$.
Then $\gt p=\hat{\gt p}\oplus\gt g(1)\oplus\gt g(2)$
and this decomposition is a $3$-grading.

Using the Killing form $\left<\,\,,\,\right>$
of $\gt g$, one can consider
$f$ as a linear function on $\gt g$ and $\gt p(e)$.
Since the pairing $\left<\gt g(-1),\gt g(1)\right>$
is non-degenerate, the skew-symmetric form
$\hat f$ is non-degenerate on $\gt g(1)$.
Thus conditions of Lemma~\ref{3-grad} are satisfied.
To conclude, note that $\gt g(0)_e=\gt g(0)_{f}$.
\end{proof}

\begin{rmk} The intersection $\gt p\cap\gt g(0)_e$ is a parabolic
subalgebra of $\gt g(0)_e$. If $\gt g$ is an exceptional Lie algebra, then
$\gt g(0)_e$ is simple and is generated by all simple roots orthogonal to the
highest one. For extended Dynkin diagrams see Table~\ref{extended} in
Section~\ref{S:exc}.
\end{rmk}

\subsection{}
In our reductions
the following situation will appear
quite often. Let $P$ be a standard proper parabolic subgroup
of $Spin_n$ (or $SO_n$) such that the Levi part of
$\gt p$ contains the first $k$ simple roots
$\alpha_1,\ldots,\ldots,\alpha_k$ of
$\gt{so}_n$ and does not contain $\alpha_{k+1}$.
Let $Pv$ be an orbit of the maximal dimension  in
the defining representation $\cp^n$.
It is assumed that $v$ is not an isotropic vector,
i.e., $(SO_n)_v=SO_{n-1}$. The orbit
$Pv$ is an open subset of a complex sphere.
In order to understand the stabiliser
$\gt p_v$, we choose a complementary to $\gt p$
subalgebra $\gt{so}_{n-1}\subset\gt{so}_n$
embedded as shown in Picture~\ref{so-1}.
As a vector space, we have  $\gt{so}_{n-1}=\gt{so}_{n-2}\oplus\cp^{n-2}$.
The big matrix is skew-symmetric with respect to
the anti-diagonal, $v$ stands for a vector in
$\cp^{n-2}$, and $v^t$ is $v$ transposed with respect to the
anti-diagonal.
The parabolic $\gt p$ is schematically shown by a dotted line,
that is, it lies above the dotted line.

\begin{figure}[htb]
{\setlength{\unitlength}{0.024in}
\begin{center}
\begin{picture}(70,70)(0,0)

\put(0,0){\line(1,0){70}}
\put(0,0){\line(0,1){70}}
\put(70,0){\line(0,1){70}}
\put(0,70){\line(1,0){70}}
\put(0,60){\line(1,0){70}}
\put(0,10){\line(1,0){70}}
\put(10,0){\line(0,1){70}}
\put(60,0){\line(0,1){70}}

\put(28,32){$\gt{so}_{n-2}$}
\put(64,32){$v$}
\put(0,32){${-}v$}
\put(4,63){$0$}
\put(30,62.5){${-}v^t$}
\put(64,63){$0$}
\put(4,2.5){$0$}
\put(33,2.5){$v^t$}
\put(64,2.5){$0$}

\qbezier[30](0,46),(12,46),(24,46)
\qbezier[30](46,0),(46,12),(46,24)

\qbezier[27](24,46),(24,35),(24,24)
\qbezier[27](24,24),(35,24),(46,24)

\put(0,47){\line(1,0){10}}
\put(0,48){\line(1,0){10}}
\put(0,49){\line(1,0){10}}
\put(0,50){\line(1,0){10}}
\put(0,51){\line(1,0){10}}
\put(0,52){\line(1,0){10}}
\put(0,53){\line(1,0){10}}
\put(0,54){\line(1,0){10}}
\put(0,55){\line(1,0){10}}
\put(0,56){\line(1,0){10}}
\put(0,57){\line(1,0){10}}
\put(0,58){\line(1,0){10}}
\put(0,59){\line(1,0){10}}

\put(60,46){\line(1,0){10}}
\put(60,47){\line(1,0){10}}
\put(60,48){\line(1,0){10}}
\put(60,49){\line(1,0){10}}
\put(60,50){\line(1,0){10}}
\put(60,51){\line(1,0){10}}
\put(60,52){\line(1,0){10}}
\put(60,53){\line(1,0){10}}
\put(60,54){\line(1,0){10}}
\put(60,55){\line(1,0){10}}
\put(60,56){\line(1,0){10}}
\put(60,57){\line(1,0){10}}
\put(60,58){\line(1,0){10}}
\put(60,59){\line(1,0){10}}

\put(47,0){\line(0,1){10}}
\put(48,0){\line(0,1){10}}
\put(49,0){\line(0,1){10}}
\put(50,0){\line(0,1){10}}
\put(51,0){\line(0,1){10}}
\put(52,0){\line(0,1){10}}
\put(53,0){\line(0,1){10}}
\put(54,0){\line(0,1){10}}
\put(55,0){\line(0,1){10}}
\put(56,0){\line(0,1){10}}
\put(57,0){\line(0,1){10}}
\put(58,0){\line(0,1){10}}
\put(59,0){\line(0,1){10}}

\put(46,60){\line(0,1){10}}
\put(47,60){\line(0,1){10}}
\put(48,60){\line(0,1){10}}
\put(49,60){\line(0,1){10}}
\put(50,60){\line(0,1){10}}
\put(51,60){\line(0,1){10}}
\put(52,60){\line(0,1){10}}
\put(53,60){\line(0,1){10}}
\put(54,60){\line(0,1){10}}
\put(55,60){\line(0,1){10}}
\put(56,60){\line(0,1){10}}
\put(57,60){\line(0,1){10}}
\put(58,60){\line(0,1){10}}
\put(59,60){\line(0,1){10}}

\put(-6,56.3){$\left\{ {\parbox{1pt}{\vspace{24\unitlength}}}  \right.$}
\put(-18,56){$k{+}1$}
\end{picture}
\end{center}}
\caption{}\label{so-1}
\end{figure}

\begin{lm}\label{wings} Let $\gt p\subset\gt{so}_n$ be
a parabolic subalgebra as above and $\cp^n$ a defining
representation of $\gt{so}_n$ with an invariant
bilinear form chosen as $(x_j,x_{n-j})=1$,
$(x_j,x_t)=0$ for $t\ne n-j$. Take $v=x_1+x_n$.
Assume further that $k\ge 2$ is even.
Then $M_*(\gt p_v)=M_*(\gt p_0)$, where
$\gt p_0\subset\gt p_v$ is a point wise stabiliser
of the plane $\cp x_1{\oplus}\cp x_n$ and a parabolic
subalgebra of $\gt{so}_{n-2}$.
\end{lm}
\begin{proof}
According to our choice, $v$ is not isotropic and
$Pv$ is an open subset of the sphere containing $v$.
Hence for $(\gt{so}_n)_v=\gt{so}_{n-1}$ we have
$(\gt{so}_n)_v+\gt p_v=\gt{so}_{n}$.
The stabiliser $\gt p_v$ decomposes
as $\gt p_v=\gt p_0\oplus V$, where
$\gt p_0=\gt p\cap\gt{so}_{n-2}$ is a parabolic in
$\gt{so}_{n-2}$, and $V$ comes, partly, from the intersection of
$\gt{gl}_{k{+}1}$-part of the Levi with "$v$-part" of the $\gt{so}_{n-1}$.
In Picture~\ref{so-1}, the embedding $V\subset\gt{so}_n$
is shown by four coated segments.
Here $\dim V=k$ and  $[\gt p_0,V]\subset V$.

Let $\gt n_0$ be the nilpotent radical
of $\gt p_0$
and $\gt p=\gt l\oplus\gt n$ a Levi decomposition.
Then $\gt n_0=\gt n\cap\gt{so}_{n-2}=\gt n\cap\gt p_v$.
On one hand, $[\gt n_0,V]\subset\gt n$,
on the other hand,
$[\gt n_0,V]\subset V$, since it lies in
the complementary to $\gt{so}_{n-2}$ subspace.
Hence $[\gt n_0, V]=0$ for the simple reason that
$V\cap\gt n=0$.
Observe that
$[V,V]\cong \Lambda^2 V$ lies in the centre of $\gt n_0$.
For $k\ge 2$, the subspace $[V,V]=\Lambda^2\cp^{k}$ is
non-zero, and in that case it coincides with the
centre of $\gt n_0$.

In the Levi part of $P$ the passage from $P$ to $P_v$ results in
replacing $GL_{k{+}1}$  by $GL_{k}\ltimes V$.
The subspace $V$ is acted upon
only by the reductive part of $P_0$, more precisely,
by $GL_{k}$.

Set $\gt a=[V,V]$. It is an Abelian ideal in $\gt p_v$.
Note that $GL_{k}$, and hence $P_0$ and $P_v$,
acts on $\gt a^*$ with an open orbit.
Suppose that a stabiliser $(\gt p_v)_\beta$ (with $\beta\in\gt p_v^*$)
is reductive. Then
the restriction $\gamma=\beta|_{\gt a}$ lies
in the open $P_v$-orbit and,
since $k$ is even, the form $\hat\beta$ is non-degenerate on $V$.
Replacing $\beta$ by an element of $\exp(V)\beta$ we may assume that
$\beta$ is zero on $V$.
Then it can be also considered as a function on
$\gt p_0$. Since $[\gt p_0,V]\subset V$,
all elements in the orbit $P_0\beta$ are zero
on $V$ and hence
$(\gt p_v)_\beta=(\gt p_0)_\beta$.

Now suppose that $\beta\in\gt p_0^*$ is of reductive type.
Then again $\gamma=\beta|_{\gt a}$ lies in the open $P_0$-orbit.
We also consider $\beta$ as a linear function on $\gt p_v$ such that
$\beta(V)=0$.
Since the form $\hat\beta$ is non-degenerate on $V$,
one gets
$(\gt p_0)_\beta=(\gt p_v)_\beta$.
Thus
the maximal reductive stabiliser of $P_v$ is the same as of $P_0$.
\end{proof}


\subsection{}\label{sb-EF}

Suppose that $\gt g$ is a Lie algebra of type F$_4$ or E$_8$.
We use Vinberg-Onishchik numbering of simple root
(see Table~\ref{VO-numb} in Section~\ref{S:exc}).
Let $H\subset G$ be a maximal parabolic subgroup
with the Lie algebra $\gt h=\gt p_{\pi}^+$
such that $\pi=\{\alpha_2,\alpha_3,\alpha_4\}$ in the  F$_4$
case and
in the E$_8$ case $\pi$
contains  to all simple roots
except $\alpha_7$.
Both these parabolics are $3$-graded, more precisely,
they are
$$(\cp^{^\times}{\times}Spin_7)\ltimes\exp(\cp^8\oplus\cp^7)
\ \text{ and } \
(\cp^{^\times}{\times}Spin_{14})\ltimes\exp(\cp^{64}\oplus\cp^{14}).
$$
In both cases $H(0)$ acts on $\gt h(1)$ via a half-spin representation and
on $\gt h(2)$ via the defining representation.
The reductive part, $H(0)$, has an open orbit in $\gt h(2)^*$ and
$\gt h(1)$ remains irreducible after the restriction to a
generic stabiliser $H(0)_*(\gt h(2)^*)$.
Hence for a generic point
$\alpha\in\gt h(2)^*$ the skew-symmetric form $\hat\alpha$
is non-degenerate on
$\gt h(1)$. Here generic means that $\alpha$ lies in
the open $H(0)$-orbit.
(Another way to see that $\hat\alpha$ is non-degenerate,
is to notice that
the above grading is related to a nilpotent element of height $2$.)
By Lemma~\ref{3-grad}, $M_*(\gt h)=\gt h(0)_\alpha$ is either $\gt{so}_6$ or $\gt{so}_{13}$,
depending on $\gt g$.

\begin{prop}\label{red-E8,F4}
Let $H\subset G$  be as above,
$P\subset G$ a parabolic subgroup, which is contained in $H$,
and $\alpha\in\gt h(2)^*$ generic.
Then $M_*(\gt p)$ is equal to the maximal reductive stabiliser of
$\gt p\cap\gt h(0)_\alpha=(\gt p\cap\gt h(0))_\alpha$.
\end{prop}
\begin{proof}
Clearly $P\cap H(0)$ is a parabolic in $H(0)$.
Let $B\subset (P\cap H(0))$ be a Borel subgroup.
As is well known, $B$ acts on $\gt h(2)^*$ with an open orbit.
Let us choose $\alpha$ such that
$B\alpha$ is that open orbit.
The nilpotent radical $\gt h(1){\oplus}\gt h(2)$ of $\gt h$ is
contained in $\gt p$.
Thus $\gt p$ is $3$-graded:
$\gt p=(\gt p\cap\gt h(0))\oplus\gt h(1)\oplus\gt h(2)$.
Since $\hat\alpha$ is still non-degenerate on $\gt h(1)$,
the claim follows from
Lemma~\ref{3-grad}.
\end{proof}

Combining Proposition~\ref{red-E8,F4} and Lemma~\ref{wings}
we get the following.

\begin{cl}\label{red-other-end}
Let $P$ be a parabolic subgroup of $G$ such that
$P\subset H$ and
\begin{itemize}
\item\ if $\gt g$ is of type F$_4$, then $\gt p=\gt p_{\pi}^+$
with  $\pi=\{\alpha_3,\alpha_4\}$;
\item\ if $\gt g$ is of type E$_8$, then the Levi part of $\gt p$
contains $\alpha_1$ and
is of type A$_k$ with even $k$.
\end{itemize}
Let $\gt{so}_m\subset\gt h(0)$ be the standard Levi subalgebra of
$\gt g$ corresponding to $\alpha_2,\alpha_3$ in
F$_4$ and $\alpha_2,\alpha_3,\ldots,\alpha_6,\alpha_8$ ($m=12$) in
E$_8$. Then
$M_*(\gt p)$ is equal to the maximal reductive stabiliser
of $\gt p_0=\gt p\cap\gt{so}_m$, where
$\gt p_0$ is a parabolic subalgebra of $\gt{so}_m$.
\end{cl}

\section{Classical Lie algebras} \label{S:cla}

\subsection{The $\gt{gl}_n$ case}\label{A}

Recall that a biparabolic (or seaweed) subalgebra $\gt q$ of a  reductive
Lie algebra $\gt g$ is an intersection
$\gt p_1\cap\gt p_2$ of two parabolic subalgebras such that
$\gt p_1+\gt p_2=\gt g$.
In case $\gt g=\gt{gl}_n$ a parabolic subalgebra
is defined up to conjugation by a flag of $\cp^n$ or
by a composition of $n$. Fixing a maximal torus in
$\gt{gl}_n$ and a root system, one may say that
a seaweed is given by two compositions of $n$.
Our goal is to describe $M_*(\gt q)$ in terms of these compositions.

\begin{rmk} If $\gt q\subset\gt{gl}_n$ is a seaweed, then
$\gt q_0:=\gt q\cap\gt{sl}_n$ is a seaweed in $\gt{sl}_n$
and $[\gt q,\gt q]\subset\gt q_0$. Since also the centre of
$GL_n$ acts on $\gt q^*$ trivially, we conclude that
$M_*(\gt q_0)=M_*(\gt q)\cap\gt{sl}_n$.
\end{rmk}

Let $(\bar a|\bar b)=(a_1,a_2,\ldots,a_m|b_1,b_2,\ldots b_l)$
be two compositions of $n$ and $\gt q=\gt q(\bar a|\bar b)$ a
corresponding seaweed in $\gt{gl}_n$.
Following \cite{DK}, we associate
to this object a graph with $n$ vertices and
several edges constructed by the following principle:
take first $a_1$ vertices and connect vertex $1$ with $a_1$,
$2$ with $a_1{-}1$ and so on; repeat it for
vertices $a_1{+}1,a_1{+}2,\ldots,a_1{+}a_2$, namely
connecting $a_1{+}1$ with $a_1{+}a_2$; do the same for all
intervals $(a_1+\ldots+a_k+1,a_1+\ldots+a_k+a_{k+1})$;
finally, repeat the procedure using the composition $\bar b$ (and the same set
of vertices).  Let $\Gamma(\bar a|\bar b)=\Gamma(\gt q)$ denote the
obtained graph. Each vertex has valency $1$ or $2$, hence
connected components of $\Gamma(\bar a|\bar b)$ are
simple cycles or segments. Below is an example of
such a graph of a seaweed in $\gt{gl}_9$.

$\Gamma(5,2,2|2,4,3)$=
{\setlength{\unitlength}{0.021in}
\raisebox{-12\unitlength}{%
\begin{picture}(110,29)(-5,-11)
\multiput(0,3)(10,0){9}{\circle*{2}}
\put(20,5){\oval(40,20)[t]}
\put(20,5){\oval(20,10)[t]}
\put(55,5){\oval(10,5)[t]}
\put(75,5){\oval(10,5)[t]}
\put(5,1){\oval(10,5)[b]}
\put(35,1){\oval(30,15)[b]}
\put(35,1){\oval(10,5)[b]}
\put(70,1){\oval(20,10)[b]}
\end{picture}
} }

\begin{df}\label{max-cycles}
Let $Y$ be a cycle, $X$ either a segment or a cycle in $\Gamma(\bar a|\bar b)$ and
$x_1>\ldots>x_r$, $y_1>\ldots>y_t$ the vertices of $X,\ Y$, respectively.
We say that $X$ {\it lies inside} $Y$ if
$y_{2i}<x_1<y_{2i-1}$ for some $i$.
(This means that in the  $2$-dimensional picture of $\Gamma(\bar a|\bar b)$
$X$ lies inside $Y$.)
We say that $X$ is {\it maximal} if it does not lie inside  any 
cycle.
\end{df}

Let us consider the simplest example when the seaweed is just
$\gt{gl}_n$. Then the corresponding graph $\Gamma(n|n)$ consist of
$[n/2]$ cycles and for $n$ odd there is also a single vertex
in the middle. Only one cycle is maximal.

$\Gamma(2|2)$=
{\setlength{\unitlength}{0.021in}
\raisebox{-12\unitlength}
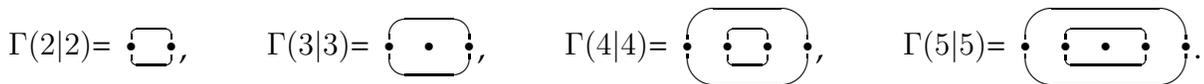
{%
\begin{picture}(15,15)(-1,-11)
\multiput(0,3)(10,0){2}{\circle*{2}}
\put(5,5){\oval(10,5)[t]}
\put(5,1){\oval(10,5)[b]}
\put(12,0){{,}}
\end{picture}
}}
\qquad
$\Gamma(3|3)$=
{\setlength{\unitlength}{0.021in}
\raisebox{-12\unitlength}{%
\begin{picture}(25,20)(-1,-11)
\multiput(0,3)(10,0){3}{\circle*{2}}
\put(10,5){\oval(20,10)[t]}
\put(10,1){\oval(20,10)[b]}
\put(22,0){{,}}
\end{picture}
}}
\qquad
$\Gamma(4|4)$=
{\setlength{\unitlength}{0.021in}
\raisebox{-12\unitlength}{%
\begin{picture}(35,29)(-1,-11)
\multiput(0,3)(10,0){4}{\circle*{2}}
\put(15,5){\oval(30,15)[t]}
\put(15,1){\oval(30,15)[b]}
\put(15,5){\oval(10,5)[t]}
\put(15,1){\oval(10,5)[b]}
\put(32,0){{,}}
\end{picture}
}}
\qquad
$\Gamma(5|5)$=
{\setlength{\unitlength}{0.021in}
\raisebox{-12\unitlength}{%
\begin{picture}(45,29)(-1,-11)
\multiput(0,3)(10,0){5}{\circle*{2}}
\put(20,5){\oval(40,15)[t]}
\put(20,1){\oval(40,15)[b]}
\put(20,5){\oval(20,5)[t]}
\put(20,1){\oval(20,5)[b]}
\put(42,0){{.}}
\end{picture}
}}

\noindent
To each cycle we attach a number, its {\it dimension},
which is equal to the sum
2\#(cycles lying inside)+\,\#(segments lying inside)+2.
We will see later that
the second summand is either $1$ or $0$. According to this formula,
the maximal cycle arising in the $\gt{gl}_n$ example has dimension $n$.
By convention,  segments are of dimension one and
if a segment is maximal,
it is considered as a maximal cycle of dimension $1$.

To each maximal cycle $X\subset\Gamma(\bar a|\bar b)$
of dimension $r$ we associate
a subgroup $GL_r\subset GL_n$,  embedded
in the following way.
Let $x_1>\ldots>x_t$ be the vertices of $X$.
If $X$ is a segment,
the corresponding $GL_1$ is
a diagonal torus with the same $c\in\mathbb C^{^\times}$ on places
$x_i$ and $1$'s on all other places.
If $X$ is not a segment, then necessary  $t$ is even
and $x_{2i-1}-x_{2i}=r{-}1$ for all $i$ (see Lemma~\ref{gl-reduction}({\sf iii}) below).
Our $GL_{r}$ is the diagonal in the product of
$t/2$ copies of $GL_{r}$ corresponding to
columns and rows
intervals $[x_{2i},x_{2i-1}]$.

\begin{thm}\label{MRS-gl_n}
The product of all $GL_r$ over all maximal cycles
in $\Gamma(\bar a,\bar b)$ is a maximal reductive stabiliser
for the corresponding seaweed.
\end{thm}

The proof is based on a reduction procedure introduced
by Panyushev \cite[Prop.~4.1]{P1}.
First, notice that if $a_1=b_1$, then
$\gt q$ is a direct sum of $\gt{gl}_{a_1}$ and
a seaweed $\gt q''\subset\gt{gl}_{n-a_1}$;
the graph $\Gamma(\bar a,\bar b)$ is a disjoint union
of $\Gamma(a_1|a_1)$ and
$\Gamma(\gt q'')$.

Suppose that $a_1\ne b_1$.
Interchanging $\bar a$ and $\bar b$, one may assume that
$a_1<b_1$.
We define a new seaweed $\gt q''\subset\gt{gl}_{n''}$
and a new tuple of numbers according to the following rule:
$$
 \begin{array}{l}
\text{ if $2a_1\le b_1$, then $n''=n-a_1$ and } \ \gt q'' \ \text{ corresponds to } \
(a_2,\ldots,a_m|b_1-2a_1,a_1,b_2,\ldots,b_l); \\
\text{ if $2a_1 > b_1$, then $n''=n-b_1+a_1$ and } \ \gt q'' \ \text{ corresponds to } \
(2a_1-b_1,a_2,\ldots,a_m|a_1,b_2,\ldots,b_l). \\
\end{array}
$$

\begin{lm}\label{gl-reduction}
\begin{itemize}
\item[({\sf i})] \ In case $a_1\ne b_1$, the passage from  $\Gamma(\gt q)$
to $\Gamma(\gt q'')$
preserves segments and cycles as well as inclusion relation
among then and therefore dimensions.
\item[({\sf ii})] \ If two connected components $X_1$, $X_2$ of
$\Gamma(\gt q)$ lie inside a cycle $Y$, then one
of them lies inside the other.
\item[({\sf iii})] \ For each cycle $X\subset\Gamma(\gt q)$ of dimension $r$ ($r>1$) with vertices
$x_1>\ldots >x_t$ all differences  $x_{2i-1}-x_{2i}$
are equal to $r{-}1$.
\end{itemize}
\end{lm}
\begin{proof}
For convenience, we draw $\Gamma(\bar a,\bar b)$ in
a 3-dimensional space, putting all vertices on a line,
preserving the order and choosing a separate plane for each edge.

For the proof of part ({\sf i}) we may assume that $a_1<b_1$.
Consider first the case where $2a_1\le b_1$. Here we contract
$\bar b$-edges connecting $1$ with $b_1$,
$2$ with $b_1{-}1$, and so on finishing with
the $\bar b$-edge $(a_1,b_1-a_1+1)$.
At the same time also  the vertices are identified in each
pair. More precisely, the pair $(i,b_1-i+1)$ is now
a single vertex $b_1-a_1+i$ for $1\ge i\ge a_1$.
Other vertices are renumbered $j\to j{-}a_1$.
Topologically speaking, the transformation
was just a contraction. Thus all cycles and segments remain cycles and segments,
no new connected components appear.
We still have a graph with no self-intersections and this new
graph corresponds now to the seaweed $\gt q''$.

Now suppose a connected component $X\subset \Gamma(\gt q)$
is lying inside a cycle $Y$ with vertices $y_1>\ldots >y_t$.
Let $x_1>\ldots>  x_r$ be the vertices of $X$.
By definition, there is $i$ such that
$y_{2i-1}>x_1>y_{2i}$. We need to show that this kind of an
inequality remains after the modification.
Note that because of the $\bar b$-edges, if $X$ has a vertex in
the interval $[b_1{-}a_1{+}1,b_1]$, then it also have a vertex
in $[1,a_1]$. (The same holds for $Y$.)
Therefore necessary  $x_1 > a_1$ and there is nothing to prove.

Consider now the second case, where $2a_1>b_1$.
Here $\Gamma(\gt q'')$ is obtained from $\Gamma(\gt q)$
in three steps. First we contract
$\bar b$-edges connecting $1$ with $b_1$,
$2$ with $b_1{-}1$, and so on finishing with
the $\bar b$-edge $(b_1{-}a_1,a_1+1)$.
At the same time also  the vertices are identified in each
pair. More precisely, the pair $(i,b_1-i+1)$ is now
a single vertex $i$ for $1\ge i\ge b_1{-}a_1$.
Next we apply the central symmetry to $[1,a_1]$
and renumber the vertices $i\to a_1{-}i{+}1$ for
$1\ge i\ge a_1$; $j\to j{-}b_1{+}a_1$ for
 $j>b_1$. Finally  the first $a_1$ $\bar a$-edges
are turned into $\bar b$-edges and
the first $(2a_1{-}b_1)$ of the $\bar b$-edges into $\bar a$-edges.
One can easily see that the passage preserves connected components.

It remains to treat connected components $X$ and $Y$.
Let $X'',\,Y''\subset\Gamma(\gt q'')$ be the modified
$X$ and $Y$.
If $x_1>a_1$, then $y_{2i-1}''>x_1''>y_{2i}''$
for the same $i$ as before the modification.
The other possibility is that $x_1\in[b_1{-}a_1{+}1,a_1]$.
Nevertheless, in that case $x_1''=x_1{-}b_1{+}a_1$ and again
$y_{2i-1}''>x_1''>y_{2i}''$
for the same $i$ as before the modification.

\vskip1ex

In order to prove parts ({\sf ii}) and ({\sf iii}), we argue by induction on
$n$. For $n=1$ there is nothing to prove.
If $a_1\ne b_1$, we can pass to $\Gamma(\gt q'')$.
By part ({\sf i}) the passage preserves connected components and inclusions among them. Moreover, $\gt q''\subset\gt{gl}_{n''}$ with $n''<n$.
Thus it remains to show that this passage preserves the
differences between vertices in a cycle $Y$.
In case $2a_1\le b_1$, we only have to notice that the
number $\#(Y\cap[1,a_1])=\#(Y\cap[b_1{-}a_1{+}1,b_1])$ is even.
In case $2a_1>b_1$, an observation that the set
$Y\cap [b_1{-}a_1{+}1,a_1]$ is invariant under the central symmetry
does the job.

If $a_1=b_1$, then
$\Gamma(\gt q)$ is a disjoint union of $\Gamma(a_1|a_1)$
and $\Gamma(\gt q'')$ with $\gt q''\subset\gt{gl}_{n-a_1}$.
Clearly both statements hold for $\Gamma(a_1|a_1)$ and by induction they
hold for $\Gamma(\gt q'')$.
This finishes the proof of the lemma.
\end{proof}

\begin{ex} We illustrate reductions of Lemma~\ref{gl-reduction}
by a seaweed $\gt q(9,3,4|4,1,11)$ in $\gt{gl}_{16}$.

$\Gamma(9,3,4|4,1,11)$=
{\setlength{\unitlength}{0.021in}
\raisebox{-12\unitlength}
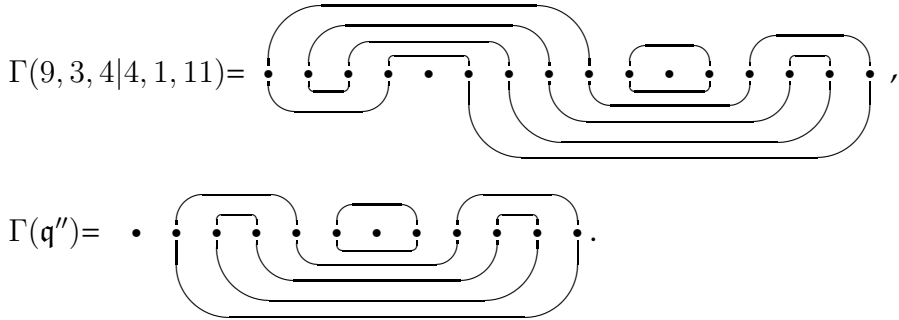
{%
\begin{picture}(160,39)(-2,-11)
\multiput(0,3)(10,0){16}{\circle*{2}}
\put(40,5){\oval(80,30)[t]}
\put(40,5){\oval(60,20)[t]}
\put(40,5){\oval(40,12)[t]}
\put(40,5){\oval(20,5)[t]}
\put(15,1){\oval(10,5)[b]}
\put(15,1){\oval(30,15)[b]}
\put(100,1){\oval(100,38)[b]}
\put(100,1){\oval(80,30)[b]}
\put(100,1){\oval(60,20)[b]}
\put(100,1){\oval(40,12)[b]}
\put(100,1){\oval(20,5)[b]}
\put(100,5){\oval(20,10)[t]}
\put(135,5){\oval(30,15)[t]}
\put(135,5){\oval(10,5)[t]}
\put(155,2){{,}}
\end{picture}
}}

$\Gamma(\gt q'')$=
{\setlength{\unitlength}{0.021in}
\raisebox{-12\unitlength}{%
\begin{picture}(120,39)(-5,-11)
\multiput(0,3)(10,0){12}{\circle*{2}}
\put(25,5){\oval(30,15)[t]}
\put(25,5){\oval(10,5)[t]}
\put(60,1){\oval(100,38)[b]}
\put(60,1){\oval(80,30)[b]}
\put(60,1){\oval(60,20)[b]}
\put(60,1){\oval(40,12)[b]}
\put(60,1){\oval(20,5)[b]}
\put(60,5){\oval(20,10)[t]}
\put(95,5){\oval(30,15)[t]}
\put(95,5){\oval(10,5)[t]}
\put(113,1){{.}}
\end{picture}
} }

\vskip3.5ex

\noindent
We are reduced to a seaweed (parabolic) $\gt p=\gt q(4,3,4|11)\subset\gt{gl}_{11}$ and

$\Gamma(\gt p'')$=
{\setlength{\unitlength}{0.021in}
\raisebox{-12\unitlength}
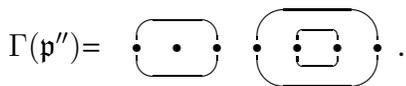
{%
\begin{picture}(70,30)(-5,-11)
\multiput(0,3)(10,0){7}{\circle*{2}}
\put(10,5){\oval(20,10)[t]}
\put(10,1){\oval(20,10)[b]}
\put(45,5){\oval(30,15)[t]}
\put(45,5){\oval(10,5)[t]}
\put(45,1){\oval(30,15)[b]}
\put(45,1){\oval(10,5)[b]}
\put(65,1){{.}}
\end{picture}
}}

\noindent
One can conclude that $M_*(\gt q)=\cp\oplus M_*(\gt p)=\cp\oplus\gt{gl}_3\oplus\gt{gl}_4$.
\end{ex}

\vskip1ex

\noindent
\begin{proof}[Proof of Theorem~\ref{MRS-gl_n}] (cf.
~\cite[proof of Theorem~4.3]{P3}.)
We argue by induction on $n$.
If $n=1$, there is nothing to prove.
Assume that $n>1$ and suppose first that $a_1=b_1$.
Then ${\rm MRS}$ of $\gt q$ is equal to
the direct product of $GL_{a_1}$ and ${\rm MRS}$ of
$\gt q''$. By induction, the theorem holds for $\gt q''\subset\gt{gl}_{n{-}a_1}$.
Suppose now that $a_1\ne b_1$.
Then
$$
\gt q=\widehat{\gt q}\ltimes V,
$$
where $\widehat{\gt q}$ is a seaweed subalgebra of $\gt{gl}_n$ defined by
$(\bar a|a_1,b_1{-}a_1,b_2,\ldots,b_l)$ and $V$ is an Abelian ideal.
Moreover $\widehat Q$ acts on $V^*$ with an open orbit.
More precisely, $V=V_1{\oplus}V_2$ with
$V_1\cong (\cp^{a_1})^*{\otimes}\cp^{c}$,
$V_2$ being isomorphic to a space of $a_1{\times}a_1$-matrices
in case $2a_1\le b_1$; and
$V_1\cong (\cp^{c})^*{\otimes}\cp^{d}$,
$V_2$ being isomorphic to a space of $d{\times}d$-matrices
for $d=b_1{-}a_1$
in case $2a_1 >b_1$;  in both cases $c=|2a_1{-}b_1|$.
We illustrate this decomposition in Picture~\ref{red-II} for
$2a_1>b_1$.
\begin{figure}[htb]
{\setlength{\unitlength}{0.1in}
\begin{center}
\begin{picture}(14,16)(0,-2)
\put(-3.1,8.7){$a_1$}
\put(-1.5,8.6){$\left\{ {\parbox{1pt}{\vspace{10\unitlength}}}  \right.$}

\put(0,4){\line(0,1){10}}
\put(14,0){\line(0,1){14}}
\put(0,14){\line(1,0){14}}

\qbezier[25](0,14)(8,6)(16,-2)

\put(0,4){\line(1,0){10}}
\put(10,-2){\line(0,1){6}}
\put(14,0){\line(1,0){4}}

\qbezier[70](0,10)(7,10)(14,10)
\qbezier[70](10,0)(10,7)(10,14)
\qbezier[50](4,4)(4,9),(4,14)
\qbezier[20](10,4)(12,4)(14,4)

\qbezier[80](0,4)(0.5,9)(1,14)
\qbezier[60](1,4)(1.5,9)(2,14)
\qbezier[60](2,4)(2.5,9)(3,14)
\qbezier[60](3,4)(3.5,9)(4,14)
\qbezier[60](4,4)(4.5,9)(5,14)
\qbezier[60](5,4)(5.5,9)(6,14)
\qbezier[60](6,4)(6.5,9)(7,14)
\qbezier[60](7,4)(7.5,9)(8,14)
\qbezier[60](8,4)(8.5,9)(9,14)
\qbezier[60](9,4)(9.5,9)(10,14)

\qbezier[40](10,-2)(10.5,1)(11,4)
\qbezier[40](11,-2)(11.5,1)(12,4)
\qbezier[40](12,-2)(12.5,1)(13,4)
\qbezier[40](13,-2)(13.5,1)(14,4)
\qbezier[10](14,-2)(14.25,-1)(14.5,0)
\qbezier[10](15,-2)(15.25,-1)(15.5,0)
\qbezier[10](16,-2)(16.25,-1)(16.5,0)
\qbezier[10](17,-2)(17.25,-1)(17.5,0)

\put(11,11.5){$V_2$}
\put(11,6.5){$V_1$}

\put(16.3,7){$b_1$}
\end{picture}
\lefteqn{\raisebox{8.73\unitlength}%
{$\left. {\parbox{1pt}{\vspace{13.65\unitlength}}}
\right\}$}}

\mbox{The subalgebra $\widehat{\gt q}$ is shaded.}
\end{center}}
\caption{}\label{red-II}
\end{figure}
In both cases there is a point $\gamma\in V^*$ in the $\widehat Q$-open orbit
such that $\gamma(V_1)=0$.
One can easily see that  $\widehat Q_\gamma\cong Q''$ for this $\gamma$.
By Lemma~\ref{2-grad}, MRS of $Q$ is equal to MRS of $Q''$
and its embedding into $Q$ can be read from the embedding of
$\widehat Q_\gamma$.
By Lemma~\ref{gl-reduction}({\sf i}),
the graphs $\Gamma(\gt q)$ and $\Gamma(\gt q'')$ have the same maximal cycles.
Therefore our descriptions of {\rm MRS} as an abstract group is justified and
it remains only to specify the embedding into $GL_n$.

Assume that $\gamma$ is given as the identity matrix in $V_2^*$.
We detail the first case, $2a_1\le b_1$.
Let $H\cong GL_{a_1}$ be a subgroup of $GL_{n}$
corresponding to a diagonal square at rows
$[(b_1{-}a_1{+}1),b_1]$. Then $\gt q''$ has a subalgebra
isomorphic to  $\gt h\cap\gt q=\gt h\cap\widehat{\gt q}$; and
$\widehat Q_\gamma$ has a subgroup isomorphic to
$H\cap Q$ embedded diagonally into
$GL_{a_1}{\times}H$.
Let $X\subset\Gamma(\gt q)$ be a maximal cycle of dimension $1$
and $GL_1$ the corresponding subgroup of the MRS of $\gt q$.
As a subgroup of $Q''$ it is defined by the properties that diagonal entries
on places $x_i''$ are all equal and entries on other diagonal places are
$1$'s.  By means of $\widehat Q_\gamma$ this embedding is extended to
$GL_{a_1}$.  If $x_i\in[1,a_1]$ is a vertex of $X$,
then also $a_1{-}x_1{+}1$ and $b_1{-}a_1{+}x_i$ are.
From the description of $\widehat Q_\gamma$, we get that diagonal entries on
places $x_i$ and $b_1{-}a_1{+}x_i$ are equal.
If $i\in[1,a_1]$ is not a vertex of $X$, then neither is
$a_1{-}i+1$ or $b_1{-}a_1{+}i$, and the diagonal entry on place $i$ is equal
to $1$.

Now let $Y\subset\Gamma(\gt q)$ be a maximal cycle of dimension $r>1$ and
$H_r\cong GL_r$ the corresponding subgroup of the MRS of $Q''$.
Recall that because of the $\bar b$-edges the number
$\#(Y\cap[a_1{+}1,b_1{-}a_1])$ is even. Hence each
interval
$[y_{2i},y_{2i-1}]$ with $y_{2i-1}\le a_1$
after going over $\bar a$-edges and $\bar b$-edges
is shifted by $b_1{-}a_1$ to an interval
$[y_{2j},y_{2j-1}]$.  Hence $H_r$ projects isomorphically on
a $GL_r$ subgroup of $GL_{a_1}$ corresponding to
columns and rows interval $[y_{2i},y_{2i-1}]$.

In case $2a_1>b_1$
the  proof goes with evident changes.
For example, here $H\cong GL_d$ corresponds to rows and columns
$[a_1{+}1,b_1]$.
One also has to notice that for a maximal cycle $Y$ of dimension $r>1$, the
number $\#(Y\cap[b_1{-}a_1{+}1,a_1])$ is even.
\end{proof}

\subsection{The $\gt{sp}_{n}$ and $\gt{so}_{n}$ cases}

In this subsection, $E$ is the vector space $\C^n$ endowed with
a non-degenerate bilinear form $B$ which is either symmetric or alternating.
Set $\ell=\left[\frac{n}{2}\right]$ and assume that $\ell\ge 1$.
We have $B(v,w)=\ep B(w,v)$ for all $v,w \in E$ where $\ep \in \{1,-1\}$.
The Lie subalgebra of $\gt{gl}_n(E)$ preserving $B$ is denoted by $\g^{\ep}$.
Thus $\g^{+1}$ is $\gt{so}(E) \simeq \gt{so}_{n}$ and $\g^{-1}$ is $\gt{sp}(E) \simeq \gt{sp}_{2\ell}$.
Let $\Pi=\{\alpha_1,\ldots,\alpha_\ell\}$ be a set of simple roots of $\g^{\ep}$.
We use
Vinberg-Onishchik numbering of simple roots, which in the classical case
coincides with the Bourbaki numbering.

The stabiliser of a flag of isotropic subspaces of $E$ in $\gt{g}^{\ep}$ is a parabolic subalgebra of $\gt{g}^{\ep}$ and
any parabolic subalgebra of $\gt{g}^{\ep}$ is obtained in this way.
A composition ${\bar a}=(a_1,\ldots,a_t)$ of an integer $r\in\{1,\ldots,\ell\}$ determines a standard (with respect to $\Pi$)
flag
$\cal{V}(\bar a)=\{\{0\}=V_0(\bar a)\varsubsetneq\cdots\varsubsetneq V_t(\bar a)\}$ of isotropic subspaces of $E$.
Note that $\dim V_t(\bar a)=r$ and $\dim V_i(\bar a)-\dim V_{i-1}(\bar a)=a_i$ for $i=1,\ldots,t$.
Let  $\p^{\ep}_{n}(\bar a)$ denote the stabiliser of $\cal{V}(\bar a)$ in $\gt{g}^{\ep}$,
and $\p_{r}(\bar a)$ will stand for the stabiliser in $\gt{gl}(V_t(\bar a))$ of the flag $\cal{V}(\bar a)$.
The Levi part of $\p^{\ep}_{n}(\bar a)$ is
$\gt{gl}_{a_1}{\oplus}\ldots{\oplus}\gt{gl}_{a_t}\oplus\gt g^{\ep}_{n-2r}$,
where $\gt g^{\ep}_{n-2r}$ is either $\gt{sp}_{n-2r}$ or
$\gt{so}_{2n-r}$, depending on $\ep$.

As has been noticed, any parabolic subalgebras of $\g^{-1} \simeq \gt{sp}_{2\ell}$ is (strongly) quasi-reductive, see~\cite{P3}.
Using a reduction of Panyushev \cite[Proof of Theorem 5.2]{P2},
we describe in the following theorem the maximal reductive stabilisers
for parabolics $\gt p\subset\gt{sp}_{2l}$.

\begin{thm} \label{t:sp}
Let $\bar a=(a_1,\ldots,a_t)$ be a composition of $r$, with $1\le r\le \ell$.
Then $$
M_*(\p_n^{-1}(\bar a)) = \gt{so}_{a_1} \oplus \cdots \oplus \gt{so}_{a_t} \oplus \gt{sp}_{n-2r}\ .
$$
The embedding of $M_*(\gt p_n^{-1}(\bar a))$ into a Levi of
$\p_n^{-1}(\bar a)$ is obvious.
\end{thm}
\begin{proof}
We argue by induction on $n=2\ell$.
Clearly the theorem is true for $n=2$.
Let $n \ge 4$ and assume that the theorem is true for any standard parabolic subalgebra of
$\gt {sp}_{n'}$ with $n'<n$.

For $\gamma \in \Delta$ and $\alpha\in\Pi$,
$[\gamma:\alpha]$ will stand for the component of $\gamma$ in $\alpha$ written in the basis $\Pi$.
We define a $\mathbb{Z}$-grading on $\g^{-1}$ by letting $\g^{-1}(i)$, for $i \not=0$, be the
sum of all root spaces $\g_{\gamma}^{-1}$ with $[\gamma:\alpha_{a_1}]=i$, and
$\g^{-1}(0)$ be the sum of $\h$ and all root spaces $\g^{-1}_{\gamma}$ with $[\gamma:\alpha_{a_1}]=0$.
Restricting this grading to $\p=\p_{n}^{-1}({\bar a})$, we obtain a 3-term grading on $\p$,
$\p=\p(0)\oplus\p(1)\oplus\p(2)$, where $\p(0)=\p\cap \g^{-1}(0)$ and $\p(i)=\g^{-1}(i)$ for $i=1,2$, see Picture~\ref{pic-red-sp}.
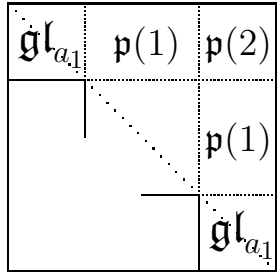
\begin{figure}[htb]
{\setlength{\unitlength}{0.1in}
\begin{center}
\begin{picture}(14,14)(0,0)

\put(0,0){\line(1,0){14}}\put(0,0){\line(0,1){14}}
\put(14,0){\line(0,1){14}}\put(0,14){\line(1,0){14}}

\qbezier[35](0,14)(7,7)(14,0)

\put(0,10){\line(1,0){4}}
\put(0.5,11.5){{\Large$\gt{gl}_{a_1}$}}
\put(4,10){\line(0,-1){3}}
\qbezier[20](4,14)(4,12)(4,10)
\qbezier[50](4,10)(9,10)(14,10)
\put(10,0){\line(0,1){4}}
\put(10,4){\line(-1,0){3}}
\put(10.5,1.5){{\Large$\gt{gl}_{a_1}$}}
\qbezier[20](10,4)(12,4)(14,4)
\qbezier[50](10,4)(10,9)(10,14)
\put(10.3,6.3){{\large$\gt p(1)$}}
\put(10.3,11.3){{\large$\gt p(2)$}}
\put(5.5,11.3){{\large$\gt p(1)$}}

\end{picture}
\end{center}}
\caption{Reduction for parabolics in $\gt{sp}_{2l}$.}\label{pic-red-sp}
\end{figure}
It follows from the construction that $\p(0)$ is isomorphic to
$\mathfrak{gl}_{a_1} \oplus \p_{n-2a_1}^{-1}(\bar b)$, where ${\bar b}=(a_2,\ldots,a_t)$,
and $\gt p(2)\cong S^2\cp^{a_1}$. For a non-degenerate (as a matrix)
$\xi \in \p(2)^*$, we have $\ad^*(\p(0))\xi=\gt p(2)$
and $\hat\xi$ is non-degenerate on $\gt p(1)$.
Therefore
Lemma~\ref{3-grad} applies and
$M_*(\p)=M_*(\p(0)_{\xi})=M_*(\gt p_{n-2a_1}^{-1}(\bar b))\oplus\gt{so}_{a_1}$.
%
%
%
By our induction applied to the parabolic subalgebra $\p_{n-2a_1}^{-1}(\bar b)$ of
$\mathfrak{sp}_{n - 2a_1}$, we obtain the expected result.
\end{proof}

We now turn to the $\mathfrak{so}_n$ case.
From now on, $B$ is assumed symmetric.
\begin{df} \label{d:star}
We will say that a composition ${\bar a}=(a_1,\ldots,a_t)$ of $r$, with $1\le r\le \ell$, satisfies the property $(\ast)$ if it does not contain
pairs $(a_i,a_{i+1})$ with $a_i$ odd and $a_{i+1}$ even.
\end{df}

If ${\bar a}=(a_1,\ldots,a_t)$ is a composition of $r$, with $1\le r\le \ell$, we set
${\bar a'}:=(a_1,\ldots,a_{t-1})$ if $r$ is odd and equal to $n/2$,
and ${\bar a'}:={\bar a}$ otherwise.

The characterisation of quasi-reductive parabolic subalgebras in
$\gt{so}_n \simeq \g^{+1}$
in term of flags has been established in~\cite[Th\'eor\`eme 5.15.1]{DKT}. 
The autors give also a characterisation in term of root systems and Dynkin diagrams~\cite[\S 5.22]{DKT}.
Recall here this result, which is also stated in~\cite[Theorem 1.7]{BM}:

\begin{thm}[\cite{DKT}] \label{t:dkt}
Let $\bar a=(a_1,\ldots,a_t)$ be a composition of $r$, with $1\le r\le \ell$.
Then, $\p^{+1}_n(\bar a)$ is quasi-reductive if and only if ${\bar a'}$ has property~$(\ast)$.
\end{thm}

\begin{rmk}
Explicit description of $M_*(\gt p)$ also proves that $\gt p$ is quasi-reductive if
it satisfies condition~$(\ast)$.
\end{rmk}

To each composition $\bar a=(a_1,\ldots,a_t)$ of $r$, with $1\le r\le \ell$,
such that ${\bar a'}$ has property $(\ast)$, and each $s \in\{1,\ldots,t\}$, we assign a subalgebra $\gt{r}_s(\bar a)$:
$$
\gt{r}_s(\bar a)=(\bigoplus\limits_{i\in\{1,\ldots,s\}, \atop a_{i} \textrm{ even} } \gt{sp}_{a_i})
\oplus (\bigoplus\limits_{i\in\{1,\ldots,s\}, \atop a_{i-1},\ a_{i} \textrm{ odd} }
\gt{sp}_{a_{i-1}-1} \oplus \gt{sp}_{a_i-1} ).
$$
By convention, $a_0:=0$ and $a_0$ is even.
Moreover, $r_0(\bar a):=0$ and $\gt{sp}_0:=0$, $\gt{so}_0:=0$.

\begin{thm} \label{t:so}
Let $\bar a=(a_1,\ldots,a_t)$ be a composition of $r$, with $1\le r\le \ell$,
such that ${\bar a}'$ has property $(\ast)$.
Then, setting $\p=\p_{n}^{+1}(\bar a)$, $M_*(\p)$ is given by the following formulas,
depending on the different cases:

\begin{itemize}
\item[{\sf (1)}]  $r$ is even:
$M_*(\p)=\r_t(\bar a)\oplus\gt{so}_{n-2r}$.
\item[{\sf (2)}] $r < \ell$ is odd: 
$M_*(\p)=\r_{t-1}(\bar a)\oplus \gt{sp}_{a_t-1} \oplus \gt{so}_{n-2r-1}$.
\item[{\sf (3)}] $r = \ell$ is odd and $a_t=1$:
$M_*(\p)=\r_{t-1}(\bar a)\oplus \C$.
\item[{\sf (4)}] $r = \ell$ is odd and $a_t > 1$ is odd:
$M_*(\p)=\r_{t-1}(\bar a)\oplus \gt{sp}_{a_t-3}$.
\item[{\sf (5)}] $r=\frac{n}{2}$ is odd and $a_t$ is even:
$M_*(\p)=\r_{t-2}(\bar a)\oplus \gt{sp}_{a_{t-1}-1} \oplus \gt{sp}_{a_t-2}$.
\end{itemize}
\end{thm}
The above cases are the only possibilities since ${\bar a}'$ has property $(\ast)$.
Note that the index of $\p$ is described in~\cite[Th\'eor\`eme 5.15.1]{DKT} for each of these cases.
More generally, a formula for the index of any biparabolic
subalgebra has been obtained in~\cite{Jo1}.

\begin{proof}
We argue by induction on $n$.
By small rank isomorphisms, e.g. $\gt{so}_5\cong\gt{sp}_4$, the statement
is known
for $n\le 6$. 
Let $n > 6$ and assume that the theorem is true for any standard
parabolic subalgebra of
$\mathfrak{so}_{n'}$ with $n'<n$.
Set $\p=\p_n^{+1}(\bar a)$ as in the theorem and set
$r_i:=a_1+\cdots +a_i$ for all $i\in\{1,\ldots,t\}$; thus $r_t=r$.

\medskip

\noindent
Step 1:
Assume first that there is $k \in\{1,\ldots,t\}$ such that $r_k$ is even.
Define a $\mathbb{Z}$-grading on $\g$ as in the proof of Theorem~\ref{t:sp} with respect
to $\alpha_{r_k}$.
So, $\p=\p(0)\oplus\p(1)\oplus\p(2)$ where $\p(i)=\g^\ep(i)$ for $i=1,2$.
We have $\p(0) = \p' \oplus \p''$ where $\p'$ and $\p''$ are parabolic subalgebras of
$\gt{gl}_{r_k}$ and $\gt{so}_{n-2r_k}$, respectively.
Namely, $\p'=\p_{r_k}({\bar a}^{(k)})$ and $\p''=\p^{+1}_{n-2r_k}(\bar b)$
with ${\bar a}^{(k)}:=(a_1,\ldots,a_k)$ and ${\bar b}:=(a_{k+1},\ldots,a_t)$.
Moreover,
$\gt p(2)\cong \Lambda^2\cp^{r_k}$ and
$\p' \oplus \p(2) \cong \p_{2r_k}^{+1}({\bar a}^{(k)})$
(cf.~\cite[proof of Theorem 6.1]{P1}),
finally $[\p'',\gt p(2)]=0$.
Let $\xi$ be a generic element of $\p(2)^\ast$.
Then
$\ad^*(\gt p(0))\xi=\gt p(2)$ and the form
$\hat{\xi}$ is non-degenerate on $\p(1)$.
By Lemma~\ref{3-grad}, $M_*(\p)$ is equal to the $M_*$ of
the stabiliser $\gt p(0)_\xi$.
In view of Lemma~\ref{2-grad}, this also can be expressed as
$M_*(\gt p)=M_*(\p(0)\oplus\p(2))$.
To summarise,
$$
M_*(\p)=M_*(\p^{+1}_{2r_k}({\bar a}^{(k)})) \oplus M_*(\p^{+1}_{n-2r_k}(\bar b))\ .
$$
Since $\bar a'$ satisfies  the property~$(\ast)$, the  inductive step does not work only
in the following three cases:

(a) $t=1$ and $a_1=n/2$;

(b) $t=1$ and $a_1$ is odd;

(c) $t=2$,  $a_1$ is odd,  and $r=n/2$.

\noindent
These different cases will be discussed in Step 2.

\medskip

\noindent
Step 2:
Define a 3-term $\mathbb{Z}$-grading on $\p$ as above with respect to $\alpha_{r}$.
We have $\p(0) \cong \gt{p}_r(\bar a) \oplus \gt{so}_{n-2r}$.
Whenever $r > 1$, $\p(1)$ is isomorphic to $\C^r \otimes \C^{n-2r}$ as
a $\p(0)$-module, and
$\p(2)\cong \Lambda^2 \C^r$.
Otherwise, $\p(1) \cong \C^{n-2}$ is Abelian, and $\p(2)=0$.
%
%
Note that there is an open $P(0)$-orbit in $\gt p(2)^*$.
From now on, we fix an element $\xi\in\p(2)^*$ in this open orbit.

\medskip

\noindent
Case (a): Here $\gt p(1)=0$ and
by Lemma~\ref{2-grad},
$M_*(\gt p)=M_*(\gt p(0)_\xi)$. If
$a_1$ is even, this is just $\gt{sp}_{a_1}$.
If $a_1=2s+1$ is odd, then
$\gt p(0)_\xi=(\cp{\oplus}\gt{sp}_{2s})\ltimes\cp^{2s}$.
This algebra is $2$-graded and, by Lemma~\ref{2-grad},
we are reduced to computing MRS of
$\gt q=\gt{sp}_{2s-2}\oplus \cp\ltimes (\cp^{2s-2}{\oplus}\cp)$,
where the second summand is a semidirect product of a
one-dimensional reductive torus and a Heisenberg Lie algebra.
Since the second summand is a Lie algebra of index zero,
we conclude that  $M_*(\gt q)=\gt{sp}_{2s-2}$ and
$M_*(\gt p)=\gt{sp}_{a_1-3}$.

\medskip

\noindent
Case (b): We may safely assume that $a_1<n/2$, i.e.,
either $1 \le r < \ell$ or $r=\ell$ and $n=2\ell+1$.
If $a_1=1$, then
$\p=\p(0)\oplus \p(1)$ and $\p(1)\cong \C^{n-2}$ is an Abelian ideal of $\p$.
Let $v$ be a non-isotropic vector of $(\C^{n-2})^*$.
Then it belongs to an open $P(0)$-orbit and
$\p(0)_v \cong (\gt{so}_{n-2})_v \cong \gt{so}_{n-3}$.
Hence, Lemma~\ref{2-grad} yields $M_*(\p) \cong \gt{so}_{n-3}$.

 Assume now that $r>1$.
 Appling Lemma~\ref{ideal} to the Abelian nilpotent ideal
 $\p(2)$, 
we get
$$
\p_{\C \xi}=(\gt{gl}_r(\xi)\oplus\gt{so}_{n-2r}) \ltimes ((\C^{r}{\otimes}\C^{n-2r} )
 \oplus \Lambda^2 \C^r),
$$
where $\gt{gl}_r(\xi)\subset \gt{gl}_r$ is the nomaliser
of $\cp\xi$.
Thus,
$\p(\xi)\cong (\gt{gl}_r(\xi)\oplus\gt{so}_{n-2r}) \ltimes (( \C^r{\otimes}\C^{n-2r} ) \oplus \C)$.
There exists a one-dimensional $\gt{gl}_r(\xi)$-invariant subspace
$\cp w\subset \cp^r$ such that $\cp w{\otimes}\cp^{n{-}2r}\subset\gt p(1)$
is the kernel of $\hat\xi$.
Let $W\subset \cp^r$ be a complement of $\cp w$ and
$\gt a=\cp w{\otimes}\cp^{n{-}2r}\oplus\cp$ an Abelian ideal
of $\gt p(\xi)$. Note that $P(\xi)$ acts on $\gt a^*$ with an open orbit.
Assertion of Lemma~\ref{ideal}, reduces computation of the MRS to
a $3$-graded Lie algebra
$$
\gt q=(\gt{so}_{n{-}2r{-}1}{\oplus}((\gt{sp}_{r-1}{\oplus}\cp)\ltimes\cp^{r-1}))
 \oplus W{\otimes}\cp^{2n-r}\oplus\cp,
$$
with the last $\cp$ being $\gt q(2)$.
We can change the first and the second grading components,
making $\gt q(0)=\gt{sp}_{r-1}{\oplus}\gt{so}_{n{-}2r{-}1}$
reductive and putting $\gt q(1)=\cp^{r-1}\oplus W{\otimes}\cp^{2n-r}$.
Here $W{\otimes}\cp^{2n-r}$ decomposes as a sum of two
$\gt q(0)$-stable subspaces $W_1\oplus W_2$ with
$\dim W_1=r{-}1$ and 
$\dim W_2=(r{-}1)(2n{-}r{-}1)$.
Moreover $[x,W_1]=\cp$ for each non-zero $x\in\cp^{r{-}1}$ and $[\cp^{r{-}1},W_2]=0$.
Since our old form $\hat\xi$ is non-degenerate on $W_2$,
we conclude that $\hat\eta$ is non-degenerate on $\gt q(1)$
for a non-zero $\eta\in\gt q(2)^*$.
By Lemma~\ref{3-grad},
$M_*(\gt q)=M_*(\gt q(0)_\eta)$.
The finall result is that
$M_*(\p)=\gt{sp}_{a_1-1}\oplus \gt{so}_{n-2r-1}$.

\medskip

\noindent
Case (c): Here $\gt p(1)=0$ and $\gt p$ has only two graded components,
with $\gt p(0)=\gt p_r(\bar a)$ being a parabolic in $\gt{gl}_r$ and
$\gt p(2)=\Lambda^2\cp^{r}$.
As above, $P(0)\xi$ is an open orbit in $\gt p(2)^*$.
By  virtue of Lemma~\ref{2-grad},
$M_*(\gt p)=M_*(\gt q)$ for $\gt q=\gt p_r(\bar a)_\xi$.
Description of $\gt q$ depends on the parity of $a_2$.
In order to calculate this stabiliser one may consider the
intersection of $\gt p(0)$ with a complementary
(in $\gt{gl}_r$) subalgebra $(\gt{gl}_r)_\xi$.

If $a_2$ is odd, then
$\gt q$ is isomorphic to
$(\gt{sp}_{a_1{-}1}{\oplus}\gt{sp}_{a_2{-}1}{\oplus}\C)\ltimes
   ((\C^{a_1-1}{\oplus} \C^{a_2-1})\oplus \C)$.
Here $(\C^{a_1-1}{\oplus} \C^{a_2-1})\oplus \C$ is the nilpotent radical,
which is also a Heisenberg Lie  algebra
with the centre $\C$.
According to Lemma~\ref{3-grad},
$M_*(q)=\gt{sp}_{a_1{-}1}{\oplus}\gt{sp}_{a_2{-}1}$.
Therefore,
we have obtained that
$M_*(\p)=\gt{sp}_{a_1{-}1}{\oplus}\gt{sp}_{a_2{-}1}$.

If $a_2$ is even, then
$\gt q=\gt q(0)\oplus\gt q(1)\oplus\gt q(2)$
with $\gt q(0)=\cp{\oplus}\cp{\oplus}\gt{sp}_{a_1{-}1}{\oplus}\gt{sp}_{a_2{-}2}$
being reductive,
$\gt q(1)=\cp^{a_1{-}1}{\oplus}V_1{\oplus}V_2$,
and $\gt q(2)=\cp_{I}{\oplus}\cp_{II}$.
The non-evident commutator relations are
$[\cp^{a_1{-}1},\cp^{a_1{-}1}]=\cp_{I}$,
$[V_1,V_1]=\cp_{I}$, $[V_1,V_2]=\cp_{II}$, and
$[V_2,V_2]=0$. Note also that
$V_1\cong V_2\cong \cp^{a_2{-}2}$.
%
%
It is not difficult to see that conditions of
Lemma~\ref{3-grad} are satisfied and therefore
$M_*(\gt p)=M_*(\gt q)=\gt{sp}_{a_1{-}1}{\oplus}\gt{sp}_{a_2{-}2}$.
\end{proof}

\section{Exceptional Lie algebras}
\label{S:exc}

In this section $\gt g$ is a simple exceptional Lie algebra
and $\Pi$ a set of simple roots of $\gt g$ as in Section~\ref{S:sea}.
For each quasi-reductive standard parabolic subalgebra $\gt p\subset\gt g$
we explicitly describe the Lie algebra $M_*(\gt p)$ of its maximal reductive stabiliser.
We use the Vinberg-Onishchik numberings of simple roots.
For convenience  of the reader, it is presented in Table~\ref{VO-numb}.

\begin{table}
\begin{tabular}{|ll|}
\hline
{\sc Type E$_6$}: \qquad \enskip
\begin{picture}(92,42)(0,6)
  \put(5,12){\circle{6}}
  \put(25,12){\circle{6}}
  \put(45,12){\circle{6}}
  \put(65,12){\circle{6}}
  \put(85,12){\circle{6}}    
  \put(45,32){\circle{6}} 
  \put(8,12){\line(1,0){14}}
  \put(28,12){\line(1,0){14}}
  \put(48,12){\line(1,0){14}}
  \put(68,12){\line(1,0){14}}
  \put(45,15){\line(0,1){14}}
\put(0,0){$\alpha_1$}
\put(20,0){$\alpha_2$}
\put(40,0){$\alpha_3$}
\put(60,0){$\alpha_4$}
\put(80,0){$\alpha_5$}
\put(40,38){$\alpha_6$}
\end{picture}
&
{\sc Type E$_7$}: \qquad \enskip
\begin{picture}(112,42)(0,6)
  \put(5,12){\circle{6}}
  \put(25,12){\circle{6}}
  \put(45,12){\circle{6}}
  \put(65,12){\circle{6}}
  \put(85,12){\circle{6}}
  \put(105,12){\circle{6}}
  \put(65,32){\circle{6}} 
  \put(8,12){\line(1,0){14}}
  \put(28,12){\line(1,0){14}}
  \put(48,12){\line(1,0){14}}
  \put(68,12){\line(1,0){14}}
  \put(88,12){\line(1,0){14}}
  \put(65,15){\line(0,1){14}}
\put(0,0){$\alpha_1$}
\put(20,0){$\alpha_2$}
\put(40,0){$\alpha_3$}
\put(60,0){$\alpha_4$}
\put(80,0){$\alpha_5$}
\put(100,0){$\alpha_6$}
\put(60,38){$\alpha_7$}
\end{picture} \\[2ex]
{\sc Type E$_8$}: \qquad \enskip
\begin{picture}(132,42)(0,6)
  \put(5,12){\circle{6}}
  \put(25,12){\circle{6}}
  \put(45,12){\circle{6}}
  \put(65,12){\circle{6}}
  \put(85,12){\circle{6}}
  \put(105,12){\circle{6}}
  \put(125,12){\circle{6}}
  \put(85,32){\circle{6}} 
  \put(8,12){\line(1,0){14}}
  \put(28,12){\line(1,0){14}}
  \put(48,12){\line(1,0){14}}
  \put(68,12){\line(1,0){14}}
  \put(88,12){\line(1,0){14}}
  \put(108,12){\line(1,0){14}}
  \put(85,15){\line(0,1){14}}
\put(0,0){$\alpha_1$}
\put(20,0){$\alpha_2$}
\put(40,0){$\alpha_3$}
\put(60,0){$\alpha_4$}
\put(80,0){$\alpha_5$}
\put(100,0){$\alpha_6$}
\put(120,0){$\alpha_7$}
\put(80,38){$\alpha_8$}
\end{picture}
& \\
{\sc Type F$_4$}: \qquad \enskip
\begin{picture}(105,30)(0,3)
  \put(10,4){\circle{6}}
  \put(40,4){\circle{6}}
  \put(70,4){\circle{6}}
  \put(100,4){\circle{6}}
    \put(13,4){\line(1,0){24}}
  \put(73,4){\line(1,0){24}}
  \put(42,6){\line(1,0){26}}
  \put(42,2){\line(1,0){26}}
  \put(50,0){\Large $<$}
\put(5,12){$\alpha_1$}
\put(35,12){$\alpha_2$}
\put(65,12){$\alpha_3$}
\put(95,12){$\alpha_4$}
\end{picture}
&
{\sc Type G$_2$}:\qquad\enskip
\begin{picture}(40,17)(0,4)
  \put(5,5){\circle{6}}
  \put(35,5){\circle{6}}
  \put(5,2){\line(1,0){30}}
  \put(8,5){\line(1,0){24}}
  \put(5,8){\line(1,0){30}}
  \put(15,1){\Large $<$}
  \put(0,11){{$\alpha_1$}}
  \put(31,11){{$\alpha_2$}}
\end{picture}
\\[0.5ex]
\hline
\end{tabular}
\vspace{.3cm}
\caption{Vinberg-Onishchik numbering of simple roots}
\label{VO-numb}
\end{table}

\begin{table}
\begin{tabular}{|ll|}
\hline
$\widetilde{{\rm E}_6}$: \qquad
\begin{picture}(90,65)
  \put(3,5){\circle{6}}
  \put(23,5){\circle{6}}
  \put(43,5){\circle{6}}
  \put(63,5){\circle{6}}
  \put(83,5){\circle{6}}
  \put(43,25){\circle{6}}
  \put(43,45){\circle{6}}
  \put(6,5){\line(1,0){14}}
  \put(26,5){\line(1,0){14}}
  \put(46,5){\line(1,0){14}}
  \put(66,5){\line(1,0){14}}
  \put(43,8){\line(0,1){14}}
  \put(43,28){\line(0,1){14}}
\put(0,-8){{\small $1$}}
\put(20,-8){{\small $2$}}
\put(40,-8){{\small $3$}}
\put(61,-8){{\small $2$}}
\put(81,-8){{\small $1$}}
\put(47,22){{\small $2$}}
\put(47,42){{\small $(1)$}}
\end{picture}
& \\
$\widetilde{{\rm E}_7}$:\qquad
\begin{picture}(130,45)
  \put(3,5){\circle{6}}
  \put(23,5){\circle{6}}
  \put(43,5){\circle{6}}
  \put(63,5){\circle{6}}
  \put(83,5){\circle{6}}
  \put(63,25){\circle{6}}
  \put(103,5){\circle{6}}
  \put(123,5){\circle{6}}

  \put(6,5){\line(1,0){14}}
  \put(26,5){\line(1,0){14}}
  \put(46,5){\line(1,0){14}}
  \put(66,5){\line(1,0){14}}
  \put(63,8){\line(0,1){14}}
\put(86,5){\line(1,0){14}}
\put(106,5){\line(1,0){14}}
\put(0,-8){{\small $1$}}
\put(20,-8){{\small $2$}}
\put(40,-8){{\small $3$}}
\put(61,-8){{\small $4$}}
\put(81,-8){{\small $3$}}
\put(101,-8){{\small $2$}}
\put(118,-8){{\small $(1)$}}
\put(67,22){{\small $2$}}
\end{picture} \qquad\qquad\qquad
&
$\widetilde{{\rm E}_8}$:\qquad
\begin{picture}(150,45)
  \put(3,5){\circle{6}}
  \put(23,5){\circle{6}}
  \put(43,5){\circle{6}}
  \put(63,5){\circle{6}}
  \put(83,5){\circle{6}}
  \put(103,25){\circle{6}}
  \put(103,5){\circle{6}}
  \put(123,5){\circle{6}}
  \put(143,5){\circle{6}}
 \put(6,5){\line(1,0){14}}
  \put(26,5){\line(1,0){14}}
  \put(46,5){\line(1,0){14}}
  \put(66,5){\line(1,0){14}}
  \put(103,8){\line(0,1){14}}
\put(86,5){\line(1,0){14}}
\put(106,5){\line(1,0){14}}
\put(126,5){\line(1,0){14}}
\put(-2.5,-8){{\small $(1)$}}
\put(20,-8){{\small $2$}}
\put(40,-8){{\small $3$}}
\put(60,-8){{\small $4$}}
\put(80,-8){{\small $5$}}
\put(100,-8){{\small $6$}}
\put(121,-8){{\small $4$}}
\put(141,-8){{\small $2$}}
\put(107,22){{\small $3$}}
\end{picture}
\\[3ex]
$\widetilde{{\rm F}_4}$:\qquad
\begin{picture}(135,15)(0,4)
  \put(10,5){\circle{6}}
  \put(40,5){\circle{6}}
  \put(70,5){\circle{6}}
  \put(100,5){\circle{6}}
  \put(130,5){\circle{6}}
  \put(13,5){\line(1,0){24}}
  \put(73,5){\line(1,0){24}}
  \put(42,7){\line(1,0){26}}
  \put(42,3){\line(1,0){26}}
  \put(50,1){\Large $<$}
  \put(103,5){\line(1,0){24}}
  \put(7,11){{\small $2$}}
   \put(37,11){{\small $4$}}
     \put(67,11){{\small $3$}}
       \put(97,11){{\small $2$}}
         \put(123.5,11){{\small $(1)$}}
\end{picture}
&
$\widetilde{{\rm G}_2}$:\qquad
\begin{picture}(70,17)(0,4)
  \put(5,5){\circle{6}}
  \put(35,5){\circle{6}}
  \put(65,5){\circle{6}}
  \put(38,5){\line(1,0){24}}
  \put(5,2){\line(1,0){30}}
  \put(8,5){\line(1,0){24}}
  \put(5,8){\line(1,0){30}}
  \put(15,1){\Large $<$}
  \put(2,11){{\small $3$}}
  \put(32,11){{\small $2$}}
  \put(59,11){{\small $(1)$}}
\end{picture}
\\[0.5ex]
\hline
\end{tabular}
\vspace{.3cm}
\caption{Extended Dynkin diagrams with the highest root coefficients}
\label{extended}
\end{table}

In the exceptional case,
there is a unique
simple root, say $\widetilde{\alpha}$, which is not orthogonal to
the highest positive root (see e.g. Table~\ref{extended}).
Set $\widetilde{\Pi}=\Pi\setminus\{\widetilde{\alpha}\}$. This is in fact the
same $\widetilde{\Pi}$ as in Section~\ref{S:sea}.
A subset $\pi \subset \Pi$ 
defines a standard parabolic subalgebra
$\p:=\gt p(\pi):=\p_{\pi}^+$ in $\gt g$, and all parabolic subalgebras arise in this way.
Let $P\subset G$ denote the corresponding (connected) parabolic
subgroup.
If $\pi \subset \widetilde{\Pi}$, then for the description of $M_*(\p)$
the highest root reduction (Prop.~\ref{min-n-red}) can be applied.
It reduces parabolics in E$_8$ to parabolics in E$_7$,
in E$_7$ to D$_6$, in E$_6$ to A$_5$,
in F$_4$ to C$_3$, and, finally, in G$_2$ to A$_1$.
Therefore we will assume that
$\widetilde{\alpha}\in \pi$.
Outside the E$_6$ type, the additivity property holds and therefore
we will consider only connected subsets $\pi$.

Unfortunately,
some explicit computations are needed,
especially in type E$_6$, where additivity does not work.
Here reductions are done by cutting ideals in the nilpotent radical of
$\gt p$.
The Lie algebra structure of $\gt p$ can be read from the
root system of $\gt g$. For example,
$\gt p$ always has an $m$-grading where
$\gt p(0)$ is the Levi part, lowest weight vectors of
$\gt p(1)$ (w.r.t. $P(0)$) correspond to
simple roots in $\Pi\setminus\pi$, and $\gt p(m{-}1)$ is an
irreducible $P(0)$-module with the highest weight vector being
the highest root vector of $\gt g$.
If $\pi=\Pi\setminus\alpha_i$, then
$m$ is the coefficient of $\alpha_i$
in the decomposition of the highest root
(these coefficients can be found in Table~\ref{extended}).

\vskip1ex

\noindent
{\bf Explanations concerning tables}: we let $\varpi_r$ denote the fundamental
weights and the $R(\varpi_r)$ corresponding irreducible representations;
embedding $M_*(\gt p)\subset\gt p$ is described in terms of the restriction to
$M_*(\gt p)$ of
the defining representation of the Levi (usually $R(\varpi_1)$);
$\id$ stands for the $1$-dimensional trivial
representation. If $\gt g$ is of type E$_6$, then MRS of $\gt p$ is not always
semisimple. We give the index of $\gt p$ and indicate generators of the
centre  using the Chevalley basis $(e_i,h_j,f_i)$
and elements $h_i\!\!\!\!^{^\vee}$ of
the 
Cartan subalgebra $\gt h\subset \gt g$
such that $\alpha_i(h_j\!\!\!\!^{^\vee})=\delta_{i,j}$.
%
%
A method for computing the index of a parabolic is given in \cite{Jo2}.

\subsection{E$_6$}
\label{sub:MRS-E6}

\vskip0.1ex

\noindent
According to \cite[Subsection 8.2]{BM},
$\gt p$ is not quasi-reductive if either
$\{\alpha_6\}$ is a connected component of
$\pi$ or $\pi=\{\alpha_1,\alpha_2,\alpha_3,\alpha_5,\alpha_6\}$
up to the diagram automorphism.
%
%
%

%
Now we describe (up to the diagram automorphism) all subsets
$\pi$ leading to quasi-reductive parabolic
subalgebras,
separating them by cardinality.
By our assumptions, $\alpha_6 \in\pi$.
Since $\alpha_6$ cannot be a connected component
of $\pi$, it must contain $\alpha_3$.
Therefore for $|\pi|=2$, there is just one
possibility, $\pi=\{\alpha_3,\alpha_6\}$.
This is item (4) in Table~\ref{tE6}.

Let $|\pi|=3$. To $\{\alpha_3,\alpha_6\}$
we can add any other root. This leads to
items (3) and (6). If $|\pi|=4$, then
again nothing is forbidden and we get
items (7), (2), (8), and (5).
The last case is $|\pi|=5$. If
$\alpha_2\not\in\pi$, the parabolic is not quasi-reductive.
Thus the only possibility is item (1).

\begin{table}
\begin{tabular}{|cc|c|c|c|}
\hline
& $\pi$ & $\ind\p_{\pi}^+$ & $M_*(\p_{\pi}^{+})$ &  Embedding \\
\hline
(1) &
\begin{picture}(90,30)
  \put(3,0){\circle*{6}}
  \put(23,0){\circle*{6}}
  \put(43,0){\circle*{6}}
  \put(63,0){\circle*{6}}
  \put(83,0){\circle{6}}    
  \put(43,20){\circle*{6}} 
  \put(6,0){\line(1,0){14}}
  \put(26,0){\line(1,0){14}}
  \put(46,0){\line(1,0){14}}
  \put(66,0){\line(1,0){14}}
  \put(43,3){\line(0,1){14}}
\end{picture}

& 2 & $\mathrm{G}_2$ & $R(\varpi_1)+3\id$    \\
(2) &

\begin{picture}(90,30)
  \put(3,0){\circle{6}}
  \put(23,0){\circle*{6}}
  \put(43,0){\circle*{6}}
  \put(63,0){\circle*{6}}
  \put(83,0){\circle{6}}    
  \put(43,20){\circle*{6}} 
  \put(6,0){\line(1,0){14}}
  \put(26,0){\line(1,0){14}}
  \put(46,0){\line(1,0){14}}
  \put(66,0){\line(1,0){14}}
  \put(43,3){\line(0,1){14}}
\end{picture}

& 4 & B$_3\oplus\cp$ & $R(\varpi_3);h_1\!\!\!\!^{^\vee}-h_5\!\!\!\!^{^\vee}$  \\
(3) &

\begin{picture}(90,30)
  \put(3,0){\circle{6}}
  \put(23,0){\circle*{6}}
  \put(43,0){\circle*{6}}
  \put(63,0){\circle{6}}
  \put(83,0){\circle{6}}    
  \put(43,20){\circle*{6}} 
  \put(6,0){\line(1,0){14}}
  \put(26,0){\line(1,0){14}}
  \put(46,0){\line(1,0){14}}
  \put(66,0){\line(1,0){14}}
  \put(43,3){\line(0,1){14}}
\end{picture}

& 2 & A$_1\oplus\cp$ &  $R(\varpi_1)+2\id;h_1\!\!\!\!^{^\vee}-h_5\!\!\!\!^{^\vee}$\\
(4) &

\begin{picture}(90,30)
  \put(3,0){\circle{6}}
  \put(23,0){\circle{6}}
  \put(43,0){\circle*{6}}
  \put(63,0){\circle{6}}
  \put(83,0){\circle{6}}    
  \put(43,20){\circle*{6}} 
  \put(6,0){\line(1,0){14}}
  \put(26,0){\line(1,0){14}}
  \put(46,0){\line(1,0){14}}
  \put(66,0){\line(1,0){14}}
  \put(43,3){\line(0,1){14}}
\end{picture}

& 3 & A$_1\oplus 2\cp$  & $R(\varpi_1)+\id;h_1\!\!\!\!^{^\vee}-h_5\!\!\!\!^{^\vee}\, ,
                    h_2-h_4$ \\
(5) &

\begin{picture}(90,30)
  \put(3,0){\circle*{6}}
  \put(23,0){\circle{6}}
  \put(43,0){\circle*{6}}
  \put(63,0){\circle*{6}}
  \put(83,0){\circle{6}}    
  \put(43,20){\circle*{6}} 
  \put(6,0){\line(1,0){14}}
  \put(26,0){\line(1,0){14}}
  \put(46,0){\line(1,0){14}}
  \put(66,0){\line(1,0){14}}
  \put(43,3){\line(0,1){14}}
\end{picture}

& 1 & A$_1$ & $2\id\otimes(R(\varpi_1)+2\id)$ \\
(6) &

\begin{picture}(90,30)
  \put(3,0){\circle*{6}}
  \put(23,0){\circle{6}}
  \put(43,0){\circle*{6}}
  \put(63,0){\circle{6}}
  \put(83,0){\circle{6}}    
  \put(43,20){\circle*{6}} 
  \put(6,0){\line(1,0){14}}
  \put(26,0){\line(1,0){14}}
  \put(46,0){\line(1,0){14}}
  \put(66,0){\line(1,0){14}}
  \put(43,3){\line(0,1){14}}
\end{picture}

&2& A$_1\oplus\cp$  &  $2\id\otimes(R(\varpi_1)+\id);h_2-h_4$ \\
(7) &

\begin{picture}(90,30)
  \put(3,0){\circle*{6}}
  \put(23,0){\circle*{6}}
  \put(43,0){\circle*{6}}
  \put(63,0){\circle{6}}
  \put(83,0){\circle{6}}    
  \put(43,20){\circle*{6}} 
  \put(6,0){\line(1,0){14}}
  \put(26,0){\line(1,0){14}}
  \put(46,0){\line(1,0){14}}
  \put(66,0){\line(1,0){14}}
  \put(43,3){\line(0,1){14}}
\end{picture}

& 0 & 0 & trivial \\
(8) &

\begin{picture}(90,30)
  \put(3,0){\circle*{6}}
  \put(23,0){\circle{6}}
  \put(43,0){\circle*{6}}
  \put(63,0){\circle{6}}
  \put(83,0){\circle*{6}}    
  \put(43,20){\circle*{6}} 
  \put(6,0){\line(1,0){14}}
  \put(26,0){\line(1,0){14}}
  \put(46,0){\line(1,0){14}}
  \put(66,0){\line(1,0){14}}
  \put(43,3){\line(0,1){14}}
\end{picture}
& 3 & $2\mathrm{A}_1\oplus\cp$ &
$R(\varpi_1)\otimes(R(\varpi_1')+\id)\otimes R(\varpi_1);
h_2\!\!\!\!^{^\vee}-h_4\!\!\!\!^{^\vee}$  \\[0.8ex]
\hline
\end{tabular}
\vspace{.7cm}
\caption{$M_*(\p_{\pi}^{+})$ in E$_6$}
\label{tE6}
\end{table}

All cases are treated separately.

\noindent
\underline{Case (1)}:\ We have
$P=(\cp^{^\times}\!{\times}Spin_{10})\ltimes\cp^{16}$.
This parabolic is $2$-graded.
The reductive part, $P(0)$, acts on the dual space of the Abelian nilpotent radical,
$\gt p(1)=\cp^{16}$, with an open orbit and
the corresponding stabiliser is
$P''=(\cp^{^\times}\!{\times}Spin_7)\ltimes \cp^8$, see e.g.
\cite{elashvili}. Due to Lemma~\ref{2-grad}, we are reduced to computing
$M_*(\gt p'')$. Here again $\cp^{^\times}\!{\times}Spin_7$ acts
on $(\cp^8)^*$ with an open orbit and the corresponding stabiliser is
of type G$_2$.

The next three parabolics are treated in a unified way.
More precisely,
parabolics (3) and (4) are contained in the second one.
We start with the largest and then apply some reductions.
In case (2), the Levi part of $P$ is
$L=P(0)=\cp^{^\times}{\times}\cp^{^\times}{\times}Spin_8$,
the nilpotent radical of $\gt p$ is two-step nilpotent
and can be decomposed as $\gt p(1){\oplus}\gt p(2)$, where
$L$ acts on
$\gt p(1)$ via $R(\varpi_1)+R(\varpi_4)$ and on $\gt p(2)$ via $R(\varpi_3)$.
All representations of the group $Spin_8$ are self dual.
The group $P(0)$ and its Borel subgroup have open orbits
in $\gt p(2)^*$. When restricted to a generic
stabiliser $L_*(\gt p(2))=\cp^{^\times}{\times}Spin_7$, both representations
$R(\varpi_1), R(\varpi_4)$ stay irreducible. Hence $\hat\alpha$
is non-degenerate on $\gt p(1)$ for generic $\alpha\in\gt p(2)^*$.
By virtue of Lemma~\ref{3-grad},
MRS of $\gt p$ is equal to $\cp^{^\times}{\times}Spin_7$.
Central part of $M_*(\gt p)$ coincides with the kernel
of the highest root of E$_6$ in the centre of $\gt p(0)$.

Set $\gt l=\Lie L$, where $L$ is the Levi part of $P(\pi_2)$ as above. 
The parabolics $\gt p(\pi_3)$, $\gt p(\pi_4)$ in lines (3) and (4)
both are semi-direct sums of
the parabolics $\gt l\cap\gt p(\pi_i)$ for $i=3,4$, respectively, and the nilpotent radical
$\gt p(1){\oplus}\gt p(2)$ of $\gt p(\pi_2)$.
Therefore, for them MRS is equal to a
maximal reductive stabiliser of the intersection
$P(\pi_i)\cap L_*(\gt p(2))$, where $L_*(\gt p(2))$ is chosen
to be transversal to $P(\pi_i)$.

\noindent
\underline{Case (3)}:\ Here the
Levi part of $P(\pi_3)$ is $\cp^{^\times}{\times}\cp^{^\times}{\times}GL_4$ and
$Q:=P(\pi_3)\cap L_*(\gt p(2))$ equals
$(\cp^{^\times}\!\times GL_3){\ltimes}\exp(\cp^3{\oplus}\Lambda^2\cp^3)$,
where the first $\cp^{^\times}$ lies in the centre of $Q$.
We can disregard the semisimple central elements and assume that
the reductive part of $Q$ is just $GL_3$.
One readily sees that $\gt q=\Lie Q$ is $3$-graded with
$\gt q(2)=\Lambda^2\cp^3$, $\gt q(1)=\cp^3$.

Take a non-zero
$\alpha\in\gt q(2)^*$, considered just as a linear function
on $\gt q(2)$, not on $\gt q$.
Then $\gt q(\alpha)=\gt q_{\cp\alpha}/\ker\alpha$ is a Lie algebra
$((\gt{gl}_2\oplus\cp)\ltimes\cp^2)\ltimes(\cp{\oplus}\cp^2\oplus\cp)$.
We apply Lemma~\ref{ideal}.
Let $V_0$ be the first $\cp^2$ (a subset of $\gt q(0)$);
$V_1=\cp$ a subset of $\gt q(1)$, $V_1'=\cp^2$ the second part
of $\gt q(1)$, and, finally, $V_2=\gt q(2)/\ker\alpha$.
Then non-zero commutators in the nilpotent radical of $\gt q(\alpha)$ are
$[V_0,V_1']=V_1$ and $[V_1',V_1']=V_2$. The centre of
$\gt{gl}_2$ acts on $V_1$ by a non-trivial character.
Thus  $\gt q(\alpha)=(\gt{gl}_2{\oplus}\cp)\oplus(V_0{\oplus} V_1')\oplus
 (V_1{\oplus}V_2)$ is $3$-graded and conditions of
Lemma~\ref{3-grad} are satisfied.
Hence $M_*(\gt q(\alpha))=\gt{sl}_2=M_*(\gt q)$.
In order to get $M_*(\gt p(\pi_3))$ we have to add the same central element as
in case (2).

\noindent
\underline{Case (4)}:\ The Levi subgroup of $P(\pi_4)$ is $(\cp^{^\times})^3{\times}GL_3$.
Hence for the description of $P(\pi_4)\cap L_*(\gt p(2))$
we can use Lemma~\ref{wings} with $k=2$.
This gives a reduction to a parabolic subalgebra $\gt p_1\subset\gt{so}_6$, which in this situation is equal to
$(\cp{\oplus}\gt{gl}_2)\ltimes (2\cp^2{\oplus}\Lambda^2\cp^2)$.
For $\gt p_1$ the maximal reductive stabiliser is $\cp^{^\times}{\times}SL_2$.
Remembering the central element in $L_*(\gt p(2))$ we get that
$M_*(\gt p(\pi_4))$ is of type $A_1\oplus 2\cp$.

Parabolics $\gt p(\pi_5)$ and $\gt p(\pi_6)$ are contained
in $\gt p(\pi_1)$ and the Abelian nilpotent radical
of $\gt p(\pi_1)$ is also an ideal of both of them.
By Lemma~\ref{2-grad},
$M_*(\gt p)$ in both these cases is equal to
$M_*(\gt p\cap\gt p'')$, where $\gt p''=(\cp{\oplus}\gt{so}_7)\ltimes\cp^8$
is a subalgebra of  $\gt p(\pi_1)$ transversal to $\gt p$.
More precisely, $\gt{so}_7$ is embedded
(via the spin-representation) into a regular
subalgebra $\gt{so}_8\subset\gt{so}_{10}\subset\gt p(\pi_1)$ and
$\cp^8$ is an $\gt{so}_8$-invariant subspace in
the complement of $\gt{so}_8$ in $\gt{so}_{10}$.

\noindent
\underline{Case (5)}:\ First we compute the intersection
$\gt q:=\gt p(\pi_5)\cap\gt p''$.
Outside of $\gt{so}_7$ it equals
$\cp\oplus\cp$ with the first $\cp$ consisting of semisimple
elements and the second of nilpotent.
To compute the intersection inside $\gt{so}_7$ we pass
to a different realisation of the embedding $Spin_7\subset SO_8$,
namely to $SO_7\subset SO_8$. Then one can easily see that the intersection
$\gt{so}_{8}\cap\gt p(\pi_5)$  is
a maximal parabolic subalgebra $\gt{gl}_4\ltimes\Lambda^2\cp^4$ of $\gt{so}_8$.
Thus $\gt q$ is a direct sum of two subalgebras, the first of which is
of index zero and the second, $\gt{gl}_3\ltimes(\cp^3{\oplus}\Lambda^2\cp^3)$,
is the same as the one we came across
in case~(3). Hence, the maximal reductive stabiliser is $\gt{sl}_2$.

\noindent
\underline{Case (6)}:\ We proceed as in case (5) and compute
$\gt p(\pi_6)\cap\gt p''$.
The only difference here is that the parabolic $\gt p(\pi_6)$
is slightly smaller than $\gt p(\pi_5)$. Namely, $\gt{gl}_4$ is replaced by
its maximal parabolic subalgebra with the reductive part $\gt{gl}_3{\oplus}\cp$.
This means that $\gt p(\pi_6)\cap\gt p''$
is a direct sum of
an index zero subalgebra and a Lie algebra $\gt p_1$ arising in case~(4).
Therefore $M_*(\gt p(\pi_6))$ is of type $A_1\oplus\cp$ and is a subalgebra
of $M_*(\gt p(\pi_4))$.

\noindent
\underline{Case (7)}:\ The parabolic is of index zero. Hence
$M_*(\gt p)=0$ and there is nothing to describe.

\noindent
\underline{Case (8)}:\
In the last line $\pi=\{\alpha_1,\alpha_3,\alpha_5,\alpha_6\}$
and it is a union of two disconected subsets
$\{\alpha_1,\alpha_5\}$ and $\{\alpha_3,\alpha_6\}$.
Let $\gt p_1\subset\gt g$ be the standard parabolic subalgebra corresponding to the first subset, and $\gt p_2$ to the second.
Since $\rk\gt g\ne \kg{\g}$, we cannot apply Theorem~\ref{t:add}.
However, we intend to show that the additivity property works for
the derived algebras of $M_*(\gt p_1)$ and of $M_*(\gt p_2)$.

By the proofs of~\cite[Corollary 2.10, Theorem 3.6, and Lemma 5.5]{BM},
there is  $(a',a'',a) \in (\C^{^\times})^3$ such that $\varphi_{v_1}$
(resp.~$\varphi_{v_2}$, $\varphi_{v}$) is regular and has reductive type for
$\gt p_1$ (resp.~$\gt p_2$, $\gt p$)
with $v_1:= a'x_{\alpha_1} + a''x_{\alpha_5} + u_{\varnothing}^- $,
$v_2:=a x_{\alpha_3+\alpha_6}+ u_{\varnothing}^-$ and $v:=v_1+v_2 - u_{\varnothing}^-$.
Let $\gt t$ (resp.~$\gt r$) be the stabiliser of
$\varphi_v$ (resp.~$\varphi_{v^{\t^{\perp}}}$) in $\p$ and
$\gt t_i$ (resp.~$\gt r_i$) the stabilisers of
$\varphi_{v_i}$ (resp.~$\varphi_{v_{i,\t_i^{\perp}}}$)
in $\p_i$, for $i=1,2$.
Let $\gt g$ stand for the derived algebra of the Levi $\gt l$ of
$\gt p$ and $\gt m$ for the nilpotent radical of $\gt p$.
The same applies to $\gt p_i$ with adding an index $i$.
Let also
$\gt s_i$ denote the image of $\gt t_i$ under the projection map from $\p_i$
to $\gt g_i\oplus\m_i$, and by $\gt k_i$ the intersection of $\gt t_i$ with
the centre of $\gt l_i$.

Let $i=1,2$.
We have, $\t_i=\gt k_i \oplus \gt s_i$ and $\t=\gt k_1 \oplus \gt s_1 \oplus \gt s_2$.
Note that $\dim \gt k_1=1$, $\dim \gt k_2=2$, and $\dim \gt s_i=1$ for $i=1,2$.
Moreover, $\gt s_i$ is contained in $\gt g_i\oplus\m$.
All this comes from~\cite{BM} (mostly from Lemma~2.9).

By the $\gt {sl}_6$ case and the highest root reduction 
(Prop.~\ref{min-n-red}),
one knows that $\gt r_1$ has type $\mathrm{A}_1\oplus \C$
and, by the case (4) of E$_6$, $\gt r_2$ has type $\mathrm{A}_1\oplus 2\C$.
The centre of $\r_i$ is $\gt k_i$ and we can
write $\r_i=\gt k_i \oplus \gt a_i$, where $\gt a_i$ is a reductive Lie subalgebra
complementary to $\gt k_i$ in $\gt r_i$ with $\gt s_i$ as a maximal torus.
Further, $\gt r = \gt k_1 \oplus \gt a$, where $\gt a$ is a reductive Lie subalgebra
complementary to $\gt k_1$ in $\gt r$ with $\gt s_1\oplus \gt s_2$ as a maximal torus.
We have $\gt a_i=\gt s_i \oplus [\gt s_i,\gt a_i]$ and
$\gt a_i \subseteq \gt g_i\oplus\gt m$.
Also, $\gt a =\gt s \oplus [\gt s,\gt a]$ and
$\gt a \subseteq \gt g_1\oplus \gt g_2\oplus \gt m$.
Our goal is to show that $\gt a= \gt a_1 \oplus \gt a_2$.

Prove first the inclusion $\gt a_i \subseteq \gt a$.
We argue as in the proof of Theorem~\ref{t:add} (inclusion $\gt r_i \subseteq \gt r$).
The same line of arguments gives
$\kappa(v_{i,\t_i^{\perp}},[\gt a_i,\p_i])=0$
since $\gt a_i$ is contained in $\gt g_i\oplus\gt m$.
Next, setting $w_1= a' x_{\alpha_1} + a'' x_{\alpha_5}$ and
$w_2=a x_{\alpha_3+\alpha_6}$, we also show
that $\kappa(v_{\t^{\perp}}-v_{i,\t_i^{\perp}},[\gt a_i,\p_i])=\{0\}$.
Indeed, it suffices to observe that $w_i \in \gt \t_j^{\perp}$ for $j \ne i$.
Therefore, we get $\gt a_i \subset \gt r$; so $\gt a_i \subset \gt a$.

Prove now the inclusion $\gt a \subseteq \gt a_1\oplus \gt a_2$.
Again, we argue as in the proof of Theorem~\ref{t:add}
(inclusion $\gt r \subseteq \gt r_1\oplus \gt r_2$).
Here, Step~1 gives us $\gt r \cap (\gt g_i \oplus \gt m) \subseteq \gt r_i$;
so $\gt r \cap (\gt g_i \oplus \gt m) \subseteq \gt a_i$.
We resume now Step 2 with $\gt a$, $\gt a_1$, and $\gt a_2$ instead of $\gt r$,
$\gt r_1$, and $\gt r_2$, respectively.
We obtain here that $\gt a=\gt a_1 \oplus \gt a_2$.

Since $\gt a_1$ and $\gt a_2$ have both type A$_1$,
we deduce that $\gt a$ has type $2{\rm A}_1$.
By the index considerations,
$M_*(\p)=\gt{sl}_2{\oplus}\gt{sl}_2{\oplus}\cp$.
We also obtain that the first A$_1$ is a diagonal
in $\gt{sl}_2{\oplus}\gt{sl}_2\subset\gt{sl}_5\subset{\rm E}_6$ and
the second is embedded  into A$_2$-subalgebra corresponding to
$\{\alpha_3,\alpha_6\}$ in the same way as in case~(4).

Now we briefly described an alternative approach to the same
parabolic $\gt p$, which makes use of the $5$-grading:
$$
\gt p=\gt{gl}_2{\oplus}\gt{gl}_2{\oplus}\gt{sl}_3\ltimes\left(
   (\cp_{{\rm I}}^{2}{\otimes}\cp^3{\oplus}\cp^3{\otimes}\cp_{{\rm II}}^2)
   \oplus
   \cp_{\rm I}{\otimes}\cp^2_{\rm II}{\otimes}(\cp^3)^*
   \oplus
   (\cp_{\rm I}^2{\oplus}\cp_{\rm II}^2)
   \oplus \cp^3\right),
$$
where subscripts ${\rm I}$, ${\rm II}$ indicate non-trivial actions of
the first and the second $\gt{sl}_2$-factors.
Note also that $[\gt p(i),\p(j)]=\p(i{+}j)$ whenever $i+j\le 4$ and
the kernel of each commutator is trivial.
Lemma~\ref{ideal} applied to $\gt a=\p(3){\oplus}\gt p(4)$ and
$\alpha\in\gt a^*$, which is zero on $\gt p(3)$ and non-zero on $\gt p(4)$,
reduces calculation of $M_*(\gt p)$ to
$$
\gt p(\alpha)=\gt{gl}_2{\oplus}\gt{gl}_2{\oplus}(\gt{gl}_2{\oplus}\cp^2)\ltimes
 V_1\oplus V_2\oplus \cp,
 $$
where $[V_1,V_1]\subset V_2$, $[V_1,V_2]=0$, $[V_2,V_2]=\cp$, and
$V_1\cong \cp^2_{\rm I}\otimes\cp^2{\oplus}\cp^2{\otimes}\cp^2_{\rm II}$,
$V_2\cong \cp^2_{\rm I}{\otimes}\cp^2_{\rm II}\oplus
 \cp^2_{\rm I}{\otimes}\cp^2_{\rm II}{\otimes}\cp^2$.
Note that the nilpotent part of $\gt p(0)_\alpha$ acts on $V_2$ non-trivially.
The central element $h$ of $M_*(\gt p)$ is defined by the conditions
that it commutes with the semisimple part of $M_*(\gt p)$,
stabilises $\alpha$, and also that $\varpi_2(h)=-\varpi_4(h)$.
Together they amount to $h=h_2\!\!\!\!^{^\vee}-h_4\!\!\!\!^{^\vee}$\ ,
up to a scalar.

\subsection{E$_7$}
\label{sub:MRS-E7}

\noindent
Recall that if $\pi$ has as a connected component one of the following sets
\begin{eqnarray*}
\{\alpha_6\}, \{\alpha_6,\alpha_5,\alpha_4\},
\{\alpha_6,\alpha_5,\alpha_4,\alpha_3,\alpha_2\},
\end{eqnarray*}
then $\p$ is not quasi-reductive \cite[Table 6]{BM}.
Quasi-reductive parabolics corresponding to connected subsets of roots
$\pi$ are listed in Table~\ref{tE7} and for each of them we find a suitable
reduction.

\begin{table}
\vspace{-.3cm}
\begin{tabular}{|cc|c|c|c|}
\hline
& $\pi$ & $\ind\p_{\pi}^+$ & $M_*(\p_{\pi}^{+})$ &  Embedding \\
\hline
(1) &

\begin{picture}(110,30)
  \put(3,0){\circle*{6}}
  \put(23,0){\circle*{6}}
  \put(43,0){\circle*{6}}
  \put(63,0){\circle*{6}}
  \put(83,0){\circle*{6}}
  \put(63,20){\circle{6}} 
  \put(103,0){\circle*{6}} 

  \put(6,0){\line(1,0){14}}
  \put(26,0){\line(1,0){14}}
  \put(46,0){\line(1,0){14}}
  \put(66,0){\line(1,0){14}}
  \put(63,3){\line(0,1){14}}
\put(86,0){\line(1,0){14}}
\end{picture}
&3& C$_3$ & $R(\varpi_1)+\id$ \\
(2) &
\begin{picture}(110,30)
  \put(3,0){\circle{6}}
  \put(23,0){\circle*{6}}
  \put(43,0){\circle*{6}}
  \put(63,0){\circle*{6}}
  \put(83,0){\circle*{6}}
  \put(63,20){\circle*{6}} 
  \put(103,0){\circle*{6}} 

  \put(6,0){\line(1,0){14}}
  \put(26,0){\line(1,0){14}}
  \put(46,0){\line(1,0){14}}
  \put(66,0){\line(1,0){14}}
  \put(63,3){\line(0,1){14}}
\put(86,0){\line(1,0){14}}
\end{picture}
&4& F$_4$ &  $R(\varpi_1)+\id$ \\
(3) &
\begin{picture}(110,30)
  \put(3,0){\circle{6}}
  \put(23,0){\circle{6}}
  \put(43,0){\circle*{6}}
  \put(63,0){\circle*{6}}
  \put(83,0){\circle*{6}}
  \put(63,20){\circle*{6}} 
  \put(103,0){\circle*{6}} 

  \put(6,0){\line(1,0){14}}
  \put(26,0){\line(1,0){14}}
  \put(46,0){\line(1,0){14}}
  \put(66,0){\line(1,0){14}}
  \put(63,3){\line(0,1){14}}
\put(86,0){\line(1,0){14}}
\end{picture}
&4& B$_4$ & $R(\varpi_1)+\id$ \\
(4) &
\begin{picture}(110,30)
  \put(3,0){\circle{6}}
  \put(23,0){\circle{6}}
  \put(43,0){\circle*{6}}
  \put(63,0){\circle*{6}}
  \put(83,0){\circle*{6}}
  \put(63,20){\circle{6}} 
  \put(103,0){\circle*{6}} 

  \put(6,0){\line(1,0){14}}
  \put(26,0){\line(1,0){14}}
  \put(46,0){\line(1,0){14}}
  \put(66,0){\line(1,0){14}}
  \put(63,3){\line(0,1){14}}
\put(86,0){\line(1,0){14}}
\end{picture}
&2& C$_2$ &  $R(\varpi_1)+\id$  \\
(5) &
\begin{picture}(110,30)
  \put(3,0){\circle{6}}
  \put(23,0){\circle{6}}
  \put(43,0){\circle{6}}
  \put(63,0){\circle*{6}}
  \put(83,0){\circle*{6}}
  \put(63,20){\circle*{6}} 
  \put(103,0){\circle*{6}} 

  \put(6,0){\line(1,0){14}}
  \put(26,0){\line(1,0){14}}
  \put(46,0){\line(1,0){14}}
  \put(66,0){\line(1,0){14}}
  \put(63,3){\line(0,1){14}}
\put(86,0){\line(1,0){14}}
\end{picture}
&2& C$_2$ &  $R(\varpi_1)+\id$  \\
(6) &
\begin{picture}(110,30)
  \put(3,0){\circle{6}}
  \put(23,0){\circle{6}}
  \put(43,0){\circle{6}}
  \put(63,0){\circle{6}}
  \put(83,0){\circle*{6}}
  \put(63,20){\circle{6}} 
  \put(103,0){\circle*{6}} 

  \put(6,0){\line(1,0){14}}
  \put(26,0){\line(1,0){14}}
  \put(46,0){\line(1,0){14}}
  \put(66,0){\line(1,0){14}}
  \put(63,3){\line(0,1){14}}
\put(86,0){\line(1,0){14}}
\end{picture}
&1& A$_1$ & $R(\varpi_1)+\id$ \\[1.4ex]
\hline
\end{tabular}
\vspace{.3cm}
\caption{$M_*(\p_{\pi}^{+})$ in E$_7$}
\label{tE7}
\end{table}

The four last items are dealt with in a unified way.
Take first a maximal parabolic subgroup $P$ with the Levi
$\cp^{^\times}{\times}Spin_{10}{\times}SL_2$, so all simple roots except
$\alpha_2$ are contained in it.
Then $\gt p$ is $3$-graded with
$\gt p(1)=\cp^{16}{\otimes}\cp^2$ and
$\gt p(2)=\cp^{10}$. There is an open $P(0)$-orbit in
$\gt p(2)^*$ with the stabiliser $Spin_9{\times}SL_2$,
$p(1)$ remains irreducible under restriction to this stabiliser,
hence the skew-symmetric form $\hat\alpha$ is non-degenerate
for generic $\alpha\in\gt p(2)^*$ and Lemma~\ref{3-grad} yields
that MRS equals $Spin_9{\times}SL_2$.

If now $P(\pi_i)$ is one of the parabolics (3)--(6) in the table,
then $P(\pi_i)\subset P$ and
$M_*(\gt p(\pi_i))$ is equal to a maximal reductive stabiliser
of the intersection
$\gt p(\pi_i)\cap\gt{so}_9$ for the complementary to $\gt p(\pi_i)$
subalgebra $\gt{so}_9$. This readily shows that
$M_*(\gt p(\pi_3))=\gt{so}_9$. In the remaining three cases
the Levi part of $P(\pi_i)$ is  either $GL_3$ or $GL_5$.
Thus Lemma~\ref{wings} provides a reduction to a parabolic $\gt p_1\subset\gt{so}_8$.
Here the answer follows from Theorem~\ref{t:so}.

Items (1) and (2) are treated separately in a different way.

\vskip0.5ex

\noindent
\underline{Case (1)}:\ Here $\gt p$ is $3$-graded:
$$
\gt p=\gt{gl}_7\oplus\Lambda^4\cp^7\oplus\cp^7.
$$
Unfortunately, a skew-symmetric form $\hat\alpha$ is
degenerate on $\gt p(1)$ for all $\alpha\in\gt p(2)^*$ by a
 simple reason  that $\gt p(1)$ is odd dimensional.
We apply Lemma~\ref{ideal} to $\gt p(2)$ and a non-zero $\alpha\in\gt p(2)^*$.
As a result we get
$$
\gt q=\gt p(\alpha)=(\gt{gl}_6\oplus\cp\oplus\cp^6)\oplus(\Lambda^3\cp^6\oplus\Lambda^4\cp^6)\oplus\cp,
$$
where $\Lambda^4\cp^6$ is the kernel of $\hat\alpha$ and
the first $\cp$ consists of semisimple elements.
This is again a $3$-grading with
$\gt q(0)=\gt{gl}_6\oplus\cp$,
$\gt q(1)=\cp^6\oplus\Lambda^3\cp^6$, and
$\gt q(2)=\Lambda^4\cp^6\oplus\cp$.
The commutators are $[\Lambda^3\cp^6,\Lambda^3\cp^6]=\cp$
(the form $\hat\alpha$ is non-degenerate on $\Lambda^3\cp^6$);
and $[\cp^6,\Lambda^3\cp^6]=\Lambda^4\cp^6$.
Let us check that for generic $\beta\in\gt q(2)^*$ the skew-symmetric
form $\hat\beta$ is non-degenerate on $\gt q(1)$.
Take $\beta=\gamma+\alpha$ with generic
$\gamma\in(\Lambda^4\cp^6)^*$.
The stabiliser of $\beta$ in $\gt q(0)$ is $\gt{sp}_6$.
As a representation of $Sp_6$, the subspace $\gt q(1)$ decomposes
as $\cp^6\oplus V(\varpi_1)\oplus V(\varpi_3)$, where
$V(\varpi_j)$ is a space of the irreducible representation
with the highest weight $\varpi_j$. Note that $\hat\alpha$ is non-degenerate on
$V(\varpi_3)$ and
$\hat\gamma$ defines a non-degenerate pairing  between $\cp^6$ and
$V(\varpi_1)$. Since $\hat\gamma(\cp^6,V(\varpi_3))=0$,
we conclude that
$\hat\beta$ is indeed non-degenerate on $\gt q(1)$ and
the maximal reductive stabiliser is $Sp_6$ by Lemma~\ref{3-grad}.

\noindent
\underline{Case (2)}:\ Here $\gt p$ is $2$-graded
with $\gt p(0)=\mathrm{E}_6\oplus\cp$ being reductive
and $\gt p(1)=V(\varpi_1)=\cp^{27}$ an Abelian ideal.
The reductive past acts on $\gt p(1)$ with an open orbit and
a generic stabiliser  F$_4$. This is a maximal reductive stabiliser.

\subsection{E$_8$}
\label{sub:MRS-E8}
Recall that if $\pi$ has as a connected component one of the following sets
\begin{eqnarray*}
&&\{\alpha_1\}, \{\alpha_3,\alpha_2,\alpha_1\},
\{\alpha_5,\alpha_4,\alpha_3,\alpha_2,\alpha_1\},\\
&&\{\alpha_7,\alpha_6,\alpha_5,\alpha_4,\alpha_3,\alpha_2,\alpha_1,\alpha_8\},
\end{eqnarray*}
then $\p$ is not quasi-reductive \cite[Table 6]{BM}.
We list all connected subsets $\pi$ corresponding to
quasi-reductive parabolics in Table~\ref{tE8}.
By the assumptions, $\alpha_1\in\pi$.
Since $\pi$ is connected, either
$\alpha_7\in\pi$, then  $\pi$ is the A$_7$-line and
$\gt p(\pi)$ is not quasi-reductive,
or $\alpha_7\not\in\pi$. As one readily sees,
no subset in Table~\ref{tE8} contains $\alpha_7$.

The item in the line (1) was considered in section~\ref{sb-EF}
and $M_*(\gt p(\pi_1))=\gt{so}_{12}$.
For all other lines Corollary~\ref{red-other-end}
provides a reduction to a parabolic in $\gt{so}_{12}$.

\hspace{-2cm}\begin{table}
\begin{tabular}{|cc|c|c|c|}
\hline
& $\pi$ & $\ind\p_{\pi}^+$& $M_*(\p_{\pi}^{+})$ &  Embedding \\
\hline
(1) &
\begin{picture}(127,30)
  \put(3,0){\circle*{6}}
  \put(23,0){\circle*{6}}
  \put(43,0){\circle*{6}}
  \put(63,0){\circle*{6}}
  \put(83,0){\circle*{6}}
  \put(83,20){\circle*{6}}  
  \put(103,0){\circle*{6}}
  \put(123,0){\circle{6}} 

  \put(6,0){\line(1,0){14}}
  \put(26,0){\line(1,0){14}}
  \put(46,0){\line(1,0){14}}
  \put(66,0){\line(1,0){14}}
  \put(83,3){\line(0,1){14}}
\put(86,0){\line(1,0){14}}
\put(106,0){\line(1,0){14}}
\end{picture}
&6& B$_6$ & $R(\varpi_1)+\id$ \\
(2) &
\begin{picture}(127,30)
  \put(3,0){\circle*{6}}
  \put(23,0){\circle*{6}}
  \put(43,0){\circle*{6}}
  \put(63,0){\circle*{6}}
  \put(83,0){\circle*{6}}
  \put(83,20){\circle{6}}  
  \put(103,0){\circle*{6}}
  \put(123,0){\circle{6}} 

  \put(6,0){\line(1,0){14}}
  \put(26,0){\line(1,0){14}}
  \put(46,0){\line(1,0){14}}
  \put(66,0){\line(1,0){14}}
  \put(83,3){\line(0,1){14}}
\put(86,0){\line(1,0){14}}
\put(106,0){\line(1,0){14}}
\end{picture}
&3& C$_3$ & $R(\varpi_1)+\id$  \\
(3) &
\begin{picture}(127,30)
  \put(3,0){\circle*{6}}
  \put(23,0){\circle*{6}}
  \put(43,0){\circle*{6}}
  \put(63,0){\circle*{6}}
  \put(83,0){\circle*{6}}
  \put(83,20){\circle*{6}}  
  \put(103,0){\circle{6}}
  \put(123,0){\circle{6}} 

  \put(6,0){\line(1,0){14}}
  \put(26,0){\line(1,0){14}}
  \put(46,0){\line(1,0){14}}
  \put(66,0){\line(1,0){14}}
  \put(83,3){\line(0,1){14}}
\put(86,0){\line(1,0){14}}
\put(106,0){\line(1,0){14}}
\end{picture}
&3& C$_3$ & $R(\varpi_1)+\id$  \\
(4) &
\begin{picture}(127,30)
  \put(3,0){\circle*{6}}
  \put(23,0){\circle*{6}}
  \put(43,0){\circle*{6}}
  \put(63,0){\circle*{6}}
  \put(83,0){\circle{6}}
  \put(83,20){\circle{6}}  
  \put(103,0){\circle{6}}
  \put(123,0){\circle{6}} 

  \put(6,0){\line(1,0){14}}
  \put(26,0){\line(1,0){14}}
  \put(46,0){\line(1,0){14}}
  \put(66,0){\line(1,0){14}}
  \put(83,3){\line(0,1){14}}
\put(86,0){\line(1,0){14}}
\put(106,0){\line(1,0){14}}
\end{picture}
&2& C$_2$ &  $R(\varpi_1)+\id$ \\
(5) &
\begin{picture}(127,30)
  \put(3,0){\circle*{6}}
  \put(23,0){\circle*{6}}
  \put(43,0){\circle{6}}
  \put(63,0){\circle{6}}
  \put(83,0){\circle{6}}
  \put(83,20){\circle{6}}  
  \put(103,0){\circle{6}}
  \put(123,0){\circle{6}} 

  \put(6,0){\line(1,0){14}}
  \put(26,0){\line(1,0){14}}
  \put(46,0){\line(1,0){14}}
  \put(66,0){\line(1,0){14}}
  \put(83,3){\line(0,1){14}}
\put(86,0){\line(1,0){14}}
\put(106,0){\line(1,0){14}}
\end{picture}
&1& A$_1$ & $R(\varpi_1)+\id$ \\[1.4ex]
\hline
\end{tabular}
\vspace{.3cm}
\caption{$M_*(\p_{\pi}^{+})$ in E$_8$}
\label{tE8}
\end{table}

\subsection{F$_4$}
\label{sub:MRS-F4}

Recall that if one of the connected components
of $\pi$ is $\alpha_{4}=\alpha_{\iota}$, then $\p(\pi)$ is
not quasi-reductive \cite[Table 6]{BM}.
By our assumption, $\alpha_4\in\pi$.
Since it cannot be a connected component,
$\alpha_3\in\pi$.
Because also $\pi$ is assumed to be connected,
we need to consider only two parabolics.

The item in the first line was treated in section~\ref{sb-EF} and here
$M_*(\gt p)=\gt{so}_6$.
For line (2) Corollary~\ref{red-other-end} provides a reduction
to a parabolic $\gt p_0$ in $\gt{so}_5$.

\begin{table}
\begin{tabular}{|cc|c|c|c|}
\hline
& $\pi$ & $\ind\p_{\pi}^+$ & $M_*(\p_{\pi}^{+})$ &  Embedding \\
\hline
(1) &
\begin{picture}(105,20)
  \put(10,4){\circle{6}}
  \put(40,4){\circle*{6}}
  \put(70,4){\circle*{6}}
  \put(100,4){\circle*{6}}
    \put(13,4){\line(1,0){24}}
  \put(73,4){\line(1,0){24}}
  \put(42,6){\line(1,0){26}}
  \put(42,2){\line(1,0){26}}
  \put(50,0){\Large $<$}
\end{picture}
&3& $\mathrm{D}_3 (\simeq \mathrm{A}_3)$  &
$R(\varpi_1)+ \id $  \\

(2) &
\begin{picture}(105,20)
  \put(10,4){\circle{6}}
  \put(40,4){\circle{6}}
  \put(70,4){\circle*{6}}
  \put(100,4){\circle*{6}}
    \put(13,4){\line(1,0){24}}
  \put(73,4){\line(1,0){24}}
  \put(42,6){\line(1,0){26}}
  \put(42,2){\line(1,0){26}}
  \put(50,0){\Large $<$}
\end{picture}
&1& A$_1$ &  $R(\varpi_1)+\id$\\[0.5ex]
\hline
\end{tabular}
\vspace{.3cm}
\caption{$M_*(\p_{\pi}^{+})$ in F$_4$}
\label{tF4}
\end{table}

\subsection{G$_2$}
The only parabolic that could be of interest
corresponds to $\pi=\{\alpha_2\}$, but it is not
quasi-reductive.

\begin{rmk} \label{r:gap}
The above computations give more elegant proofs for the fact that the parabolic
subalgebras appearing in tables of Section 6 are quasi-reductive than what is done
in~\cite{BM} (where \texttt{GAP} was used to check the result in cases with
$|\pi| \geq 3$.)
\end{rmk}

\section{Conclusion}

As was already noticed, if $\gt q\subset\gt{gl}_n$ is a seaweed,
then MRS is connected.
In general, MRS are not always connected.
Take for example the parabolic in line (1) of Table~\ref{tF4}.
A maximal reductive stabiliser is equal to
$\mathbb Z/2\mathbb Z\times SL_{4}$ and the component group
acts on $M_*(\gt p)$ non-trivially.

In this paper we have studied maximal reductive stabilisers
on the Lie algebra level.
Given any particular quasi-reductive parabolic subalgebra it is
not difficult to compute the component group of its MRS, for instance,
along the lines indicated in the Introduction.
Also
our explicit calculations presented in Section~\ref{S:exc}
allow to get MRS on the group level for all considered Lie algebras.
Nevertheless, an intriguing question remains, does
the ``additivity'' property hold on the group level?
We believe it does and would like to find a general proof,
not a case-by-case verification.

Contemplating results of Section~\ref{S:cla}, one may notice that
``additivity'' works for all parabolics in type C.
If $\gt p(\pi)\subset\gt{sl}_n$ is a parabolic and
$\pi=\pi_1\cup\pi_2$, where $\pi_1,\pi_2$ are not connected to each other
and both are invariant
under the diagram automorphism, then ``additivity'' holds for the derived algebras of
$M_*(\gt p)$, $M_*(\gt p(\pi_1)), M_*(\gt p(\pi_2))$.
A more careful analysis shows that it holds
in type B and in type D fails only
for D$_{2m{+}1}$ and either $r=2m{+}1$, $a_t$ even
(case (5) of Theorem~\ref{t:so}) or $r=n-4$, $a_t$ odd.
 As is well known, a Weyl involution is inner in type D$_{2m}$
 and outer in type D$_{2m{+}1}$.
In the exceptional case, the  ``additivity'' does not work only in E$_6$ and here
holds in the same sense as for type A.
All these observations 
indicate
a relation with the existence of outer automorphisms
or, more precisely, with the condition that a Weyl involution is inner.
It would be nice to clarify this dependence.

Two  further questions in the same spirit were suggested to us by Michel Duflo.
What is the dimension of the centre of a maximal reductive stabilser
and what is the number of its simple factors?
Analysing
Section~\ref{S:exc} and Theorem~\ref{t:so}, one might conclude that
these questions are also 
related to the existence of outer automorphisms.

One of  the most interesting direction for further investigation is:
what happens for real Lie algebras? What kind of unitary representations
come out as results of Duflo's construction for $\gamma\in\gt q^*_{\rm red}$
of nilpotent type or other linear forms of reductive type?

Stabilisers $\gt q_\gamma$ play an important r{\^ole} in representation theory,
in all characteristics. They are of particular interest if
$\gamma$ is of nilpotent type, i.e., vanishes on $\gt q_{\gamma}$
(for relevant results in finite characteristic see e.g. \cite{Pr}, \cite{PS}).
Studying the nilpotent cone in $\gt q^*$ and linear forms of nilpotent type
on $\gt q$ for strongly quasi-reductive Lie algebras looks like a very promising
project.


\end{document}